\input amstex
%

\def\next{AMS-SEKR}\ifx\styname\next \endinput\fi
\catcode`\@=11
\def\styname{AMS-SEKR}
\def\styversion{2.0}
{\W@{}\W@{\styname.STY - Version \styversion}\W@{}}
\hyphenation{acad-e-my acad-e-mies af-ter-thought anom-aly anom-alies
an-ti-deriv-a-tive an-tin-o-my an-tin-o-mies apoth-e-o-ses apoth-e-o-sis
ap-pen-dix ar-che-typ-al as-sign-a-ble as-sist-ant-ship as-ymp-tot-ic
asyn-chro-nous at-trib-uted at-trib-ut-able bank-rupt bank-rupt-cy
bi-dif-fer-en-tial blue-print busier busiest cat-a-stroph-ic
cat-a-stroph-i-cally con-gress cross-hatched data-base de-fin-i-tive
de-riv-a-tive dis-trib-ute dri-ver dri-vers eco-nom-ics econ-o-mist
elit-ist equi-vari-ant ex-quis-ite ex-tra-or-di-nary flow-chart
for-mi-da-ble forth-right friv-o-lous ge-o-des-ic ge-o-det-ic geo-met-ric
griev-ance griev-ous griev-ous-ly hexa-dec-i-mal ho-lo-no-my ho-mo-thetic
ideals idio-syn-crasy in-fin-ite-ly in-fin-i-tes-i-mal ir-rev-o-ca-ble
key-stroke lam-en-ta-ble light-weight mal-a-prop-ism man-u-script
mar-gin-al meta-bol-ic me-tab-o-lism meta-lan-guage me-trop-o-lis
met-ro-pol-i-tan mi-nut-est mol-e-cule mono-chrome mono-pole mo-nop-oly
mono-spline mo-not-o-nous mul-ti-fac-eted mul-ti-plic-able non-euclid-ean
non-iso-mor-phic non-smooth par-a-digm par-a-bol-ic pa-rab-o-loid
pa-ram-e-trize para-mount pen-ta-gon phe-nom-e-non post-script pre-am-ble
pro-ce-dur-al pro-hib-i-tive pro-hib-i-tive-ly pseu-do-dif-fer-en-tial
pseu-do-fi-nite pseu-do-nym qua-drat-ics quad-ra-ture qua-si-smooth
qua-si-sta-tion-ary qua-si-tri-an-gu-lar quin-tes-sence quin-tes-sen-tial
re-arrange-ment rec-tan-gle ret-ri-bu-tion retro-fit retro-fit-ted
right-eous right-eous-ness ro-bot ro-bot-ics sched-ul-ing se-mes-ter
semi-def-i-nite semi-ho-mo-thet-ic set-up se-vere-ly side-step sov-er-eign
spe-cious spher-oid spher-oid-al star-tling star-tling-ly
sta-tis-tics sto-chas-tic straight-est strange-ness strat-a-gem strong-hold
sum-ma-ble symp-to-matic syn-chro-nous topo-graph-i-cal tra-vers-a-ble
tra-ver-sal tra-ver-sals treach-ery turn-around un-at-tached un-err-ing-ly
white-space wide-spread wing-spread wretch-ed wretch-ed-ly Brown-ian
Eng-lish Euler-ian Feb-ru-ary Gauss-ian Grothen-dieck Hamil-ton-ian
Her-mit-ian Jan-u-ary Japan-ese Kor-te-weg Le-gendre Lip-schitz
Lip-schitz-ian Mar-kov-ian Noe-ther-ian No-vem-ber Rie-mann-ian
Schwarz-schild Sep-tem-ber
form per-iods Uni-ver-si-ty cri-ti-sism for-ma-lism}
\Invalid@\nofrills
\Invalid@\usualspace
\newif\ifnofrills@
\def\nofrills@#1#2{\relaxnext@
  \DN@{\ifx\next\nofrills
    \nofrills@true\let#2\relax\DN@\nofrills{\nextii@}%
  \else
    \nofrills@false\def#2{#1}\let\next@\nextii@\fi
\next@}}
\def\usualspace@#1{\ifnofrills@\def\usualspace{#1}\fi}
\def\addto#1#2{\csname \expandafter\eat@\string#1@\endcsname
  \expandafter{\the\csname \expandafter\eat@\string#1@\endcsname#2}}
\newdimen\bigsize@
\def\big@#1#2{{\hbox{$\left#2\vcenter to#1\bigsize@{}%
  \right.\nulldelimiterspace\z@\m@th$}}}
\def\big{\big@\@ne}
\def\Big{\big@{1.5}}
\def\bigg{\big@\tw@}
\def\Bigg{\big@{2.5}}
\def\raggedcenter@{\leftskip\z@ plus.4\hsize \rightskip\leftskip
 \parfillskip\z@ \parindent\z@ \spaceskip.3333em \xspaceskip.5em
 \pretolerance9999\tolerance9999 \exhyphenpenalty\@M
 \hyphenpenalty\@M \let\\\linebreak}
\def\upperspecialchars{\def\ss{SS}\let\i=I\let\j=J\let\ae\AE\let\oe\OE
  \let\o\O\let\aa\AA\let\l\L}
\def\uppercasetext@#1{%
  {\spaceskip1.2\fontdimen2\the\font plus1.2\fontdimen3\the\font
   \upperspecialchars\uctext@#1$\m@th\aftergroup\eat@$}}
\def\uctext@#1$#2${\endash@#1-\endash@$#2$\uctext@}
\def\endash@#1-#2\endash@{\uppercase{#1}\if\notempty{#2}--\endash@#2\endash@\fi}
\def\runaway@#1{\DN@{#1}\ifx\envir@\next@
  \Err@{You seem to have a missing or misspelled \string\end#1 ...}%
  \let\envir@\empty\fi}
\newif\iftemp@
\def\notempty#1{TT\fi\def\test@{#1}\ifx\test@\empty\temp@false
  \else\temp@true\fi \iftemp@}
\font@\tensmc=cmcsc10
\font@\sevenex=cmex7
\font@\sevenit=cmti7
\font@\eightrm=cmr8 
\font@\sixrm=cmr6 
\font@\eighti=cmmi8     \skewchar\eighti='177 
\font@\sixi=cmmi6       \skewchar\sixi='177   
\font@\eightsy=cmsy8    \skewchar\eightsy='60 
\font@\sixsy=cmsy6      \skewchar\sixsy='60   
\font@\eightex=cmex8
\font@\eightbf=cmbx8 
\font@\sixbf=cmbx6   
\font@\eightit=cmti8 
\font@\eightsl=cmsl8 
\font@\eightsmc=cmcsc8
\font@\eighttt=cmtt8 

\loadmsam
\loadmsbm
\loadeufm
\UseAMSsymbols
\newtoks\tenpoint@
\def\tenpoint{\normalbaselineskip12\p@
 \abovedisplayskip12\p@ plus3\p@ minus9\p@
 \belowdisplayskip\abovedisplayskip
 \abovedisplayshortskip\z@ plus3\p@
 \belowdisplayshortskip7\p@ plus3\p@ minus4\p@
 \textonlyfont@\rm\tenrm \textonlyfont@\it\tenit
 \textonlyfont@\sl\tensl \textonlyfont@\bf\tenbf
 \textonlyfont@\smc\tensmc \textonlyfont@\tt\tentt
 \textonlyfont@\bsmc\tenbsmc
 \ifsyntax@ \def\big##1{{\hbox{$\left##1\right.$}}}%
  \let\Big\big \let\bigg\big \let\Bigg\big
 \else
  \textfont\z@=\tenrm  \scriptfont\z@=\sevenrm  \scriptscriptfont\z@=\fiverm
  \textfont\@ne=\teni  \scriptfont\@ne=\seveni  \scriptscriptfont\@ne=\fivei
  \textfont\tw@=\tensy \scriptfont\tw@=\sevensy \scriptscriptfont\tw@=\fivesy
  \textfont\thr@@=\tenex \scriptfont\thr@@=\sevenex
        \scriptscriptfont\thr@@=\sevenex
  \textfont\itfam=\tenit \scriptfont\itfam=\sevenit
        \scriptscriptfont\itfam=\sevenit
  \textfont\bffam=\tenbf \scriptfont\bffam=\sevenbf
        \scriptscriptfont\bffam=\fivebf
  \setbox\strutbox\hbox{\vrule height8.5\p@ depth3.5\p@ width\z@}%
  \setbox\strutbox@\hbox{\lower.5\normallineskiplimit\vbox{%
        \kern-\normallineskiplimit\copy\strutbox}}%
 \setbox\z@\vbox{\hbox{$($}\kern\z@}\bigsize@=1.2\ht\z@
 \fi
 \normalbaselines\rm\ex@.2326ex\jot3\ex@\the\tenpoint@}
\newtoks\eightpoint@
\def\eightpoint{\normalbaselineskip10\p@
 \abovedisplayskip10\p@ plus2.4\p@ minus7.2\p@
 \belowdisplayskip\abovedisplayskip
 \abovedisplayshortskip\z@ plus2.4\p@
 \belowdisplayshortskip5.6\p@ plus2.4\p@ minus3.2\p@
 \textonlyfont@\rm\eightrm \textonlyfont@\it\eightit
 \textonlyfont@\sl\eightsl \textonlyfont@\bf\eightbf
 \textonlyfont@\smc\eightsmc \textonlyfont@\tt\eighttt
 \textonlyfont@\bsmc\eightbsmc
 \ifsyntax@\def\big##1{{\hbox{$\left##1\right.$}}}%
  \let\Big\big \let\bigg\big \let\Bigg\big
 \else
  \textfont\z@=\eightrm \scriptfont\z@=\sixrm \scriptscriptfont\z@=\fiverm
  \textfont\@ne=\eighti \scriptfont\@ne=\sixi \scriptscriptfont\@ne=\fivei
  \textfont\tw@=\eightsy \scriptfont\tw@=\sixsy \scriptscriptfont\tw@=\fivesy
  \textfont\thr@@=\eightex \scriptfont\thr@@=\sevenex
   \scriptscriptfont\thr@@=\sevenex
  \textfont\itfam=\eightit \scriptfont\itfam=\sevenit
   \scriptscriptfont\itfam=\sevenit
  \textfont\bffam=\eightbf \scriptfont\bffam=\sixbf
   \scriptscriptfont\bffam=\fivebf
 \setbox\strutbox\hbox{\vrule height7\p@ depth3\p@ width\z@}%
 \setbox\strutbox@\hbox{\raise.5\normallineskiplimit\vbox{%
   \kern-\normallineskiplimit\copy\strutbox}}%
 \setbox\z@\vbox{\hbox{$($}\kern\z@}\bigsize@=1.2\ht\z@
 \fi
 \normalbaselines\eightrm\ex@.2326ex\jot3\ex@\the\eightpoint@}

\font@\twelverm=cmr10 scaled\magstep1
\font@\twelveit=cmti10 scaled\magstep1
\font@\twelvesl=cmsl10 scaled\magstep1
\font@\twelvesmc=cmcsc10 scaled\magstep1
\font@\twelvett=cmtt10 scaled\magstep1
\font@\twelvebf=cmbx10 scaled\magstep1
\font@\twelvei=cmmi10 scaled\magstep1
\font@\twelvesy=cmsy10 scaled\magstep1
\font@\twelveex=cmex10 scaled\magstep1
\font@\twelvemsa=msam10 scaled\magstep1
\font@\twelveeufm=eufm10 scaled\magstep1
\font@\twelvemsb=msbm10 scaled\magstep1
\newtoks\twelvepoint@
\def\twelvepoint{\normalbaselineskip15\p@
 \abovedisplayskip15\p@ plus3.6\p@ minus10.8\p@
 \belowdisplayskip\abovedisplayskip
 \abovedisplayshortskip\z@ plus3.6\p@
 \belowdisplayshortskip8.4\p@ plus3.6\p@ minus4.8\p@
 \textonlyfont@\rm\twelverm \textonlyfont@\it\twelveit
 \textonlyfont@\sl\twelvesl \textonlyfont@\bf\twelvebf
 \textonlyfont@\smc\twelvesmc \textonlyfont@\tt\twelvett
 \textonlyfont@\bsmc\twelvebsmc
 \ifsyntax@ \def\big##1{{\hbox{$\left##1\right.$}}}%
  \let\Big\big \let\bigg\big \let\Bigg\big
 \else
  \textfont\z@=\twelverm  \scriptfont\z@=\tenrm  \scriptscriptfont\z@=\sevenrm
  \textfont\@ne=\twelvei  \scriptfont\@ne=\teni  \scriptscriptfont\@ne=\seveni
  \textfont\tw@=\twelvesy \scriptfont\tw@=\tensy \scriptscriptfont\tw@=\sevensy
  \textfont\thr@@=\twelveex \scriptfont\thr@@=\tenex
        \scriptscriptfont\thr@@=\tenex
  \textfont\itfam=\twelveit \scriptfont\itfam=\tenit
        \scriptscriptfont\itfam=\tenit
  \textfont\bffam=\twelvebf \scriptfont\bffam=\tenbf
        \scriptscriptfont\bffam=\sevenbf
  \setbox\strutbox\hbox{\vrule height10.2\p@ depth4.2\p@ width\z@}%
  \setbox\strutbox@\hbox{\lower.6\normallineskiplimit\vbox{%
        \kern-\normallineskiplimit\copy\strutbox}}%
 \setbox\z@\vbox{\hbox{$($}\kern\z@}\bigsize@=1.4\ht\z@
 \fi
 \normalbaselines\rm\ex@.2326ex\jot3.6\ex@\the\twelvepoint@}

\def\headfonts{\twelvepoint\bf}

\font@\fourteenrm=cmr10 scaled\magstep2
\font@\fourteenit=cmti10 scaled\magstep2
\font@\fourteensl=cmsl10 scaled\magstep2
\font@\fourteensmc=cmcsc10 scaled\magstep2
\font@\fourteentt=cmtt10 scaled\magstep2
\font@\fourteenbf=cmbx10 scaled\magstep2
\font@\fourteeni=cmmi10 scaled\magstep2
\font@\fourteensy=cmsy10 scaled\magstep2
\font@\fourteenex=cmex10 scaled\magstep2
\font@\fourteenmsa=msam10 scaled\magstep2
\font@\fourteeneufm=eufm10 scaled\magstep2
\font@\fourteenmsb=msbm10 scaled\magstep2
\newtoks\fourteenpoint@
\def\fourteenpoint{\normalbaselineskip15\p@
 \abovedisplayskip18\p@ plus4.3\p@ minus12.9\p@
 \belowdisplayskip\abovedisplayskip
 \abovedisplayshortskip\z@ plus4.3\p@
 \belowdisplayshortskip10.1\p@ plus4.3\p@ minus5.8\p@
 \textonlyfont@\rm\fourteenrm \textonlyfont@\it\fourteenit
 \textonlyfont@\sl\fourteensl \textonlyfont@\bf\fourteenbf
 \textonlyfont@\smc\fourteensmc \textonlyfont@\tt\fourteentt
 \textonlyfont@\bsmc\fourteenbsmc
 \ifsyntax@ \def\big##1{{\hbox{$\left##1\right.$}}}%
  \let\Big\big \let\bigg\big \let\Bigg\big
 \else
  \textfont\z@=\fourteenrm  \scriptfont\z@=\twelverm  \scriptscriptfont\z@=\tenrm
  \textfont\@ne=\fourteeni  \scriptfont\@ne=\twelvei  \scriptscriptfont\@ne=\teni
  \textfont\tw@=\fourteensy \scriptfont\tw@=\twelvesy \scriptscriptfont\tw@=\tensy
  \textfont\thr@@=\fourteenex \scriptfont\thr@@=\twelveex
        \scriptscriptfont\thr@@=\twelveex
  \textfont\itfam=\fourteenit \scriptfont\itfam=\twelveit
        \scriptscriptfont\itfam=\twelveit
  \textfont\bffam=\fourteenbf \scriptfont\bffam=\twelvebf
        \scriptscriptfont\bffam=\tenbf
  \setbox\strutbox\hbox{\vrule height12.2\p@ depth5\p@ width\z@}%
  \setbox\strutbox@\hbox{\lower.72\normallineskiplimit\vbox{%
        \kern-\normallineskiplimit\copy\strutbox}}%
 \setbox\z@\vbox{\hbox{$($}\kern\z@}\bigsize@=1.7\ht\z@
 \fi
 \normalbaselines\rm\ex@.2326ex\jot4.3\ex@\the\fourteenpoint@}

\def\chapheadfonts{\fourteenpoint\bf}

\font@\seventeenrm=cmr10 scaled\magstep3
\font@\seventeenit=cmti10 scaled\magstep3
\font@\seventeensl=cmsl10 scaled\magstep3
\font@\seventeensmc=cmcsc10 scaled\magstep3
\font@\seventeentt=cmtt10 scaled\magstep3
\font@\seventeenbf=cmbx10 scaled\magstep3
\font@\seventeeni=cmmi10 scaled\magstep3
\font@\seventeensy=cmsy10 scaled\magstep3
\font@\seventeenex=cmex10 scaled\magstep3
\font@\seventeenmsa=msam10 scaled\magstep3
\font@\seventeeneufm=eufm10 scaled\magstep3
\font@\seventeenmsb=msbm10 scaled\magstep3
\newtoks\seventeenpoint@
\def\seventeenpoint{\normalbaselineskip18\p@
 \abovedisplayskip21.6\p@ plus5.2\p@ minus15.4\p@
 \belowdisplayskip\abovedisplayskip
 \abovedisplayshortskip\z@ plus5.2\p@
 \belowdisplayshortskip12.1\p@ plus5.2\p@ minus7\p@
 \textonlyfont@\rm\seventeenrm \textonlyfont@\it\seventeenit
 \textonlyfont@\sl\seventeensl \textonlyfont@\bf\seventeenbf
 \textonlyfont@\smc\seventeensmc \textonlyfont@\tt\seventeentt
 \textonlyfont@\bsmc\seventeenbsmc
 \ifsyntax@ \def\big##1{{\hbox{$\left##1\right.$}}}%
  \let\Big\big \let\bigg\big \let\Bigg\big
 \else
  \textfont\z@=\seventeenrm  \scriptfont\z@=\fourteenrm  \scriptscriptfont\z@=\twelverm
  \textfont\@ne=\seventeeni  \scriptfont\@ne=\fourteeni  \scriptscriptfont\@ne=\twelvei
  \textfont\tw@=\seventeensy \scriptfont\tw@=\fourteensy \scriptscriptfont\tw@=\twelvesy
  \textfont\thr@@=\seventeenex \scriptfont\thr@@=\fourteenex
        \scriptscriptfont\thr@@=\fourteenex
  \textfont\itfam=\seventeenit \scriptfont\itfam=\fourteenit
        \scriptscriptfont\itfam=\fourteenit
  \textfont\bffam=\seventeenbf \scriptfont\bffam=\fourteenbf
        \scriptscriptfont\bffam=\twelvebf
  \setbox\strutbox\hbox{\vrule height14.6\p@ depth6\p@ width\z@}%
  \setbox\strutbox@\hbox{\lower.86\normallineskiplimit\vbox{%
        \kern-\normallineskiplimit\copy\strutbox}}%
 \setbox\z@\vbox{\hbox{$($}\kern\z@}\bigsize@=2\ht\z@
 \fi
 \normalbaselines\rm\ex@.2326ex\jot5.2\ex@\the\seventeenpoint@}

\font@\rrrrrm=cmr10 scaled\magstep4
\font@\bigtitlefont=cmbx10 scaled\magstep4

\parindent1pc
\normallineskiplimit\p@
\newdimen\indenti \indenti=2pc
\def\pageheight#1{\vsize#1}
\def\pagewidth#1{\hsize#1%
   \captionwidth@\hsize \advance\captionwidth@-2\indenti}
\pagewidth{30pc} \pageheight{47pc}
\def\topmatter{%
 \ifx\undefined\msafam
 \else\font@\eightmsa=msam8 \font@\sixmsa=msam6
   \ifsyntax@\else \addto\tenpoint{\textfont\msafam=\tenmsa
              \scriptfont\msafam=\sevenmsa \scriptscriptfont\msafam=\fivemsa}%
     \addto\eightpoint{\textfont\msafam=\eightmsa \scriptfont\msafam=\sixmsa
              \scriptscriptfont\msafam=\fivemsa}%
   \fi
 \fi
 \ifx\undefined\msbfam
 \else\font@\eightmsb=msbm8 \font@\sixmsb=msbm6
   \ifsyntax@\else \addto\tenpoint{\textfont\msbfam=\tenmsb
         \scriptfont\msbfam=\sevenmsb \scriptscriptfont\msbfam=\fivemsb}%
     \addto\eightpoint{\textfont\msbfam=\eightmsb \scriptfont\msbfam=\sixmsb
         \scriptscriptfont\msbfam=\fivemsb}%
   \fi
 \fi
 \ifx\undefined\eufmfam
 \else \font@\eighteufm=eufm8 \font@\sixeufm=eufm6
   \ifsyntax@\else \addto\tenpoint{\textfont\eufmfam=\teneufm
       \scriptfont\eufmfam=\seveneufm \scriptscriptfont\eufmfam=\fiveeufm}%
     \addto\eightpoint{\textfont\eufmfam=\eighteufm
       \scriptfont\eufmfam=\sixeufm \scriptscriptfont\eufmfam=\fiveeufm}%
   \fi
 \fi
 \ifx\undefined\eufbfam
 \else \font@\eighteufb=eufb8 \font@\sixeufb=eufb6
   \ifsyntax@\else \addto\tenpoint{\textfont\eufbfam=\teneufb
      \scriptfont\eufbfam=\seveneufb \scriptscriptfont\eufbfam=\fiveeufb}%
    \addto\eightpoint{\textfont\eufbfam=\eighteufb
      \scriptfont\eufbfam=\sixeufb \scriptscriptfont\eufbfam=\fiveeufb}%
   \fi
 \fi
 \ifx\undefined\eusmfam
 \else \font@\eighteusm=eusm8 \font@\sixeusm=eusm6
   \ifsyntax@\else \addto\tenpoint{\textfont\eusmfam=\teneusm
       \scriptfont\eusmfam=\seveneusm \scriptscriptfont\eusmfam=\fiveeusm}%
     \addto\eightpoint{\textfont\eusmfam=\eighteusm
       \scriptfont\eusmfam=\sixeusm \scriptscriptfont\eusmfam=\fiveeusm}%
   \fi
 \fi
 \ifx\undefined\eusbfam
 \else \font@\eighteusb=eusb8 \font@\sixeusb=eusb6
   \ifsyntax@\else \addto\tenpoint{\textfont\eusbfam=\teneusb
       \scriptfont\eusbfam=\seveneusb \scriptscriptfont\eusbfam=\fiveeusb}%
     \addto\eightpoint{\textfont\eusbfam=\eighteusb
       \scriptfont\eusbfam=\sixeusb \scriptscriptfont\eusbfam=\fiveeusb}%
   \fi
 \fi
 \ifx\undefined\eurmfam
 \else \font@\eighteurm=eurm8 \font@\sixeurm=eurm6
   \ifsyntax@\else \addto\tenpoint{\textfont\eurmfam=\teneurm
       \scriptfont\eurmfam=\seveneurm \scriptscriptfont\eurmfam=\fiveeurm}%
     \addto\eightpoint{\textfont\eurmfam=\eighteurm
       \scriptfont\eurmfam=\sixeurm \scriptscriptfont\eurmfam=\fiveeurm}%
   \fi
 \fi
 \ifx\undefined\eurbfam
 \else \font@\eighteurb=eurb8 \font@\sixeurb=eurb6
   \ifsyntax@\else \addto\tenpoint{\textfont\eurbfam=\teneurb
       \scriptfont\eurbfam=\seveneurb \scriptscriptfont\eurbfam=\fiveeurb}%
    \addto\eightpoint{\textfont\eurbfam=\eighteurb
       \scriptfont\eurbfam=\sixeurb \scriptscriptfont\eurbfam=\fiveeurb}%
   \fi
 \fi
 \ifx\undefined\cmmibfam
 \else \font@\eightcmmib=cmmib8 \font@\sixcmmib=cmmib6
   \ifsyntax@\else \addto\tenpoint{\textfont\cmmibfam=\tencmmib
       \scriptfont\cmmibfam=\sevencmmib \scriptscriptfont\cmmibfam=\fivecmmib}%
    \addto\eightpoint{\textfont\cmmibfam=\eightcmmib
       \scriptfont\cmmibfam=\sixcmmib \scriptscriptfont\cmmibfam=\fivecmmib}%
   \fi
 \fi
 \ifx\undefined\cmbsyfam
 \else \font@\eightcmbsy=cmbsy8 \font@\sixcmbsy=cmbsy6
   \ifsyntax@\else \addto\tenpoint{\textfont\cmbsyfam=\tencmbsy
      \scriptfont\cmbsyfam=\sevencmbsy \scriptscriptfont\cmbsyfam=\fivecmbsy}%
    \addto\eightpoint{\textfont\cmbsyfam=\eightcmbsy
      \scriptfont\cmbsyfam=\sixcmbsy \scriptscriptfont\cmbsyfam=\fivecmbsy}%
   \fi
 \fi
 \let\topmatter\relax}
\def\chapterno@{\uppercase\expandafter{\romannumeral\chaptercount@}}
\newcount\chaptercount@
\def\chapter{\nofrills@{\afterassignment\chapterno@
                        CHAPTER \global\chaptercount@=}\chapter@
 \DNii@##1{\leavevmode\hskip-\leftskip
   \rlap{\vbox to\z@{\vss\centerline{\eightpoint
   \chapter@##1\unskip}\baselineskip2pc\null}}\hskip\leftskip
   \nofrills@false}%
 \FN@\next@}
\newbox\titlebox@

\def\title{\nofrills@{\relax}\title@%
 \DNii@##1\endtitle{\global\setbox\titlebox@\vtop{\tenpoint\bf
 \raggedcenter@\ignorespaces
 \baselineskip1.3\baselineskip\title@{##1}\endgraf}%
 \ifmonograph@ \edef\next{\the\leftheadtoks}\ifx\next\empty
    \leftheadtext{##1}\fi
 \fi
 \edef\next{\the\rightheadtoks}\ifx\next\empty \rightheadtext{##1}\fi
 }\FN@\next@}
\newbox\authorbox@
\def\author#1\endauthor{\global\setbox\authorbox@
 \vbox{\tenpoint\smc\raggedcenter@\ignorespaces
 #1\endgraf}\relaxnext@ \edef\next{\the\leftheadtoks}%
 \ifx\next\empty\leftheadtext{#1}\fi}
\newbox\affilbox@
\def\affil#1\endaffil{\global\setbox\affilbox@
 \vbox{\tenpoint\raggedcenter@\ignorespaces#1\endgraf}}
\newcount\addresscount@
\addresscount@\z@
\def\address#1\endaddress{\global\advance\addresscount@\@ne
  \expandafter\gdef\csname address\number\addresscount@\endcsname
  {\vskip12\p@ minus6\p@\noindent\eightpoint\smc\ignorespaces#1\par}}
\def\email{\nofrills@{\eightpoint{\it E-mail\/}:\enspace}\email@
  \DNii@##1\endemail{%
  \expandafter\gdef\csname email\number\addresscount@\endcsname
  {\def\usualspace{{\it\enspace}}\smallskip\noindent\eightpoint\email@
  \ignorespaces##1\par}}%
 \FN@\next@}
\def\thedate@{}
\def\date#1\enddate{\gdef\thedate@{\tenpoint\ignorespaces#1\unskip}}
\def\thethanks@{}
\def\thanks#1\endthanks{\gdef\thethanks@{\eightpoint\ignorespaces#1.\unskip}}
\def\thekeywords@{}
\def\keywords{\nofrills@{{\it Key words and phrases.\enspace}}\keywords@
 \DNii@##1\endkeywords{\def\thekeywords@{\def\usualspace{{\it\enspace}}%
 \eightpoint\keywords@\ignorespaces##1\unskip.}}%
 \FN@\next@}
\def\thesubjclass@{}
\def\subjclass{\nofrills@{{\rm2020 {\it Mathematics Subject
   Classification\/}.\enspace}}\subjclass@
 \DNii@##1\endsubjclass{\def\thesubjclass@{\def\usualspace
  {{\rm\enspace}}\eightpoint\subjclass@\ignorespaces##1\unskip.}}%
 \FN@\next@}
\newbox\abstractbox@
\def\abstract{\nofrills@{{\smc Abstract.\enspace}}\abstract@
 \DNii@{\setbox\abstractbox@\vbox\bgroup\noindent$$\vbox\bgroup
  \def\envir@{abstract}\advance\hsize-2\indenti
  \usualspace@{{\enspace}}\eightpoint \noindent\abstract@\ignorespaces}%
 \FN@\next@}
\def\endabstract{\par\unskip\egroup$$\egroup}
\def\widestnumber#1#2{\begingroup\let\head\null\let\subhead\empty
   \let\subsubhead\subhead
   \ifx#1\head\global\setbox\tocheadbox@\hbox{#2.\enspace}%
   \else\ifx#1\subhead\global\setbox\tocsubheadbox@\hbox{#2.\enspace}%
   \else\ifx#1\key\bgroup\let\endrefitem@\egroup
        \key#2\endrefitem@\global\refindentwd\wd\keybox@
   \else\ifx#1\no\bgroup\let\endrefitem@\egroup
        \no#2\endrefitem@\global\refindentwd\wd\nobox@
   \else\ifx#1\page\global\setbox\pagesbox@\hbox{\quad\bf#2}%
   \else\ifx#1\item\setboxz@h{#2}\global\rosteritemwd\wdz@
        \global\advance\rosteritemwd by.5\parindent
   \else\message{\string\widestnumber is not defined for this option
   (\string#1)}%
\fi\fi\fi\fi\fi\fi\endgroup}
\newif\ifmonograph@
\def\Monograph{\monograph@true \let\headmark\rightheadtext
  \let\varindent@\indent \def\headfont@{\bf}\def\proclaimheadfont@{\smc}%
  \def\demofont@{\smc}}
\let\varindent@\indent

\newbox\tocheadbox@    \newbox\tocsubheadbox@
\newbox\tocbox@
\def\toc{\toc@{Contents}}
\def\newtocdefs{%
   \def \title##1\endtitle
       {\penaltyandskip@\z@\smallskipamount
        \hangindent\wd\tocheadbox@\noindent{\bf##1}}%
   \def \chapter##1{%
        Chapter \uppercase\expandafter{\romannumeral##1.\unskip}\enspace}%
   \def \specialhead##1\endspecialhead
       {\par\hangindent\wd\tocheadbox@ \noindent##1\par}%
   \def \head##1 ##2\endhead
       {\par\hangindent\wd\tocheadbox@ \noindent
        \if\notempty{##1}\hbox to\wd\tocheadbox@{\hfil##1\enspace}\fi
        ##2\par}%
   \def \subhead##1 ##2\endsubhead
       {\par\vskip-\parskip {\normalbaselines
        \advance\leftskip\wd\tocheadbox@
        \hangindent\wd\tocsubheadbox@ \noindent
        \if\notempty{##1}\hbox to\wd\tocsubheadbox@{##1\unskip\hfil}\fi
         ##2\par}}%
   \def \subsubhead##1 ##2\endsubsubhead
       {\par\vskip-\parskip {\normalbaselines
        \advance\leftskip\wd\tocheadbox@
        \hangindent\wd\tocsubheadbox@ \noindent
        \if\notempty{##1}\hbox to\wd\tocsubheadbox@{##1\unskip\hfil}\fi
        ##2\par}}}
\def\toc@#1{\relaxnext@
   \def\page##1%
       {\unskip\penalty0\null\hfil
        \rlap{\hbox to\wd\pagesbox@{\quad\hfil##1}}\hfilneg\penalty\@M}%
 \DN@{\ifx\next\nofrills\DN@\nofrills{\nextii@}%
      \else\DN@{\nextii@{{#1}}}\fi
      \next@}%
 \DNii@##1{%
\ifmonograph@\bgroup\else\setbox\tocbox@\vbox\bgroup
   \centerline{\headfont@\ignorespaces##1\unskip}\nobreak
   \vskip\belowheadskip \fi
   \setbox\tocheadbox@\hbox{0.\enspace}%
   \setbox\tocsubheadbox@\hbox{0.0.\enspace}%
   \leftskip\indenti \rightskip\leftskip
   \setbox\pagesbox@\hbox{\bf\quad000}\advance\rightskip\wd\pagesbox@
   \newtocdefs
 }%
 \FN@\next@}
\def\endtoc{\par\egroup}
\let\pretitle\relax
\let\preauthor\relax
\let\preaffil\relax
\let\predate\relax
\let\preabstract\relax
\let\prepaper\relax
\def\dedicatory #1\enddedicatory{\def\preabstract{{\medskip
  \eightpoint\it \raggedcenter@#1\endgraf}}}
\def\thetranslator@{}
\def\translator#1\endtranslator{\def\thetranslator@{\nobreak\medskip
 \line{\eightpoint\hfil Translated by \uppercase{#1}\qquad\qquad}\nobreak}}
\outer\def\endtopmatter{\runaway@{abstract}%
 \edef\next{\the\leftheadtoks}\ifx\next\empty
  \expandafter\leftheadtext\expandafter{\the\rightheadtoks}\fi
 \ifmonograph@\else
   \ifx\thesubjclass@\empty\else \makefootnote@{}{\thesubjclass@}\fi
   \ifx\thekeywords@\empty\else \makefootnote@{}{\thekeywords@}\fi
   \ifx\thethanks@\empty\else \makefootnote@{}{\thethanks@}\fi
 \fi
  \pretitle
  \ifmonograph@ \topskip7pc \else \topskip4pc \fi
  \box\titlebox@
  \topskip10pt
  \preauthor
  \ifvoid\authorbox@\else \vskip2.5pc plus1pc \unvbox\authorbox@\fi
  \preaffil
  \ifvoid\affilbox@\else \vskip1pc plus.5pc \unvbox\affilbox@\fi
  \predate
  \ifx\thedate@\empty\else \vskip1pc plus.5pc \line{\hfil\thedate@\hfil}\fi
  \preabstract
  \ifvoid\abstractbox@\else \vskip1.5pc plus.5pc \unvbox\abstractbox@ \fi
  \ifvoid\tocbox@\else\vskip1.5pc plus.5pc \unvbox\tocbox@\fi
  \prepaper
  \vskip2pc plus1pc
}
\def\document{\let\fontlist@\relax\let\alloclist@\relax
  \tenpoint}

\newskip\aboveheadskip       \aboveheadskip1.8\bigskipamount
\newdimen\belowheadskip      \belowheadskip1.8\medskipamount

\def\headfont@{\smc}
\def\penaltyandskip@#1#2{\relax\ifdim\lastskip<#2\relax\removelastskip
      \ifnum#1=\z@\else\penalty@#1\relax\fi\vskip#2%
  \else\ifnum#1=\z@\else\penalty@#1\relax\fi\fi}
\def\nobreak{\penalty\@M
  \ifvmode\def\penalty@{\let\penalty@\penalty\count@@@}%
  \everypar{\let\penalty@\penalty\everypar{}}\fi}
\let\penalty@\penalty
\def\heading#1\endheading{\head#1\endhead}

\def\specialheadfont@{\bf}
\outer\def\specialhead{\par\penaltyandskip@{-200}\aboveheadskip
  \begingroup\interlinepenalty\@M\rightskip\z@ plus\hsize \let\\\linebreak
  \specialheadfont@\noindent\ignorespaces}
\def\endspecialhead{\par\endgroup\nobreak\vskip\belowheadskip}
\let\headmark\eat@
\newskip\subheadskip       \subheadskip\medskipamount
\def\subheadfont@{\bf}
\outer\def\subhead{\nofrills@{.\enspace}\subhead@
 \DNii@##1\endsubhead{\par\penaltyandskip@{-100}\subheadskip
  \varindent@{\usualspace@{{\subheadfont@\enspace}}%
 \subheadfont@\ignorespaces##1\unskip\subhead@}\ignorespaces}%
 \FN@\next@}
\outer\def\subsubhead{\nofrills@{.\enspace}\subsubhead@
 \DNii@##1\endsubsubhead{\par\penaltyandskip@{-50}\medskipamount
      {\usualspace@{{\it\enspace}}%
  \it\ignorespaces##1\unskip\subsubhead@}\ignorespaces}%
 \FN@\next@}
\def\proclaimheadfont@{\bf}
\outer\def\proclaim{\runaway@{proclaim}\def\envir@{proclaim}%
  \nofrills@{.\enspace}\proclaim@
 \DNii@##1{\penaltyandskip@{-100}\medskipamount\varindent@
   \usualspace@{{\proclaimheadfont@\enspace}}\proclaimheadfont@
   \ignorespaces##1\unskip\proclaim@
  \sl\ignorespaces}%
 \FN@\next@}
\outer\def\endproclaim{\let\envir@\relax\par\rm
  \penaltyandskip@{55}\medskipamount}
\def\demoheadfont@{\it}
\def\demo{\runaway@{proclaim}\nofrills@{.\enspace}\demo@
     \DNii@##1{\par\penaltyandskip@\z@\medskipamount
  {\usualspace@{{\demoheadfont@\enspace}}%
  \varindent@\demoheadfont@\ignorespaces##1\unskip\demo@}\rm
  \ignorespaces}\FN@\next@}
\def\enddemo{\par\medskip}
\def\qed{\ifhmode\unskip\nobreak\fi\quad\ifmmode\square\else$\m@th\square$\fi}
\let\remark\demo
\let\endremark\enddemo
\def\definition{\runaway@{proclaim}%
  \nofrills@{.\demoheadfont@\enspace}\definition@
        \DNii@##1{\penaltyandskip@{-100}\medskipamount
        {\usualspace@{{\demoheadfont@\enspace}}%
        \varindent@\demoheadfont@\ignorespaces##1\unskip\definition@}%
        \rm \ignorespaces}\FN@\next@}


\newdimen\rosteritemwd
\newcount\rostercount@
\newif\iffirstitem@
\let\plainitem@\item
\newtoks\everypartoks@
\def\par@{\everypartoks@\expandafter{\the\everypar}\everypar{}}
\def\roster{\edef\leftskip@{\leftskip\the\leftskip}%
 \relaxnext@
 \rostercount@\z@  
 \def\item{\FN@\rosteritem@}%
 \DN@{\ifx\next\runinitem\let\next@\nextii@\else
  \let\next@\nextiii@\fi\next@}%
 \DNii@\runinitem  
  {\unskip  
   \DN@{\ifx\next[\let\next@\nextii@\else
    \ifx\next"\let\next@\nextiii@\else\let\next@\nextiv@\fi\fi\next@}%
   \DNii@[####1]{\rostercount@####1\relax
    \enspace{\rm(\number\rostercount@)}~\ignorespaces}%
   \def\nextiii@"####1"{\enspace{\rm####1}~\ignorespaces}%
   \def\nextiv@{\enspace{\rm(1)}\rostercount@\@ne~}%
   \par@\firstitem@false  
   \FN@\next@}%
 \def\nextiii@{\par\par@  
  \penalty\@m\smallskip\vskip-\parskip
  \firstitem@true}%
 \FN@\next@}
\def\rosteritem@{\iffirstitem@\firstitem@false\else\par\vskip-\parskip\fi
 \leftskip3\parindent\noindent  
 \DNii@[##1]{\rostercount@##1\relax
  \llap{\hbox to2.5\parindent{\hss\rm(\number\rostercount@)}%
   \hskip.5\parindent}\ignorespaces}%
 \def\nextiii@"##1"{%
  \llap{\hbox to2.5\parindent{\hss\rm##1}\hskip.5\parindent}\ignorespaces}%
 \def\nextiv@{\advance\rostercount@\@ne
  \llap{\hbox to2.5\parindent{\hss\rm(\number\rostercount@)}%
   \hskip.5\parindent}}%
 \ifx\next[\let\next@\nextii@\else\ifx\next"\let\next@\nextiii@\else
  \let\next@\nextiv@\fi\fi\next@}

\newif\ifnextRunin@
\def\endroster{\relaxnext@
 \par\leftskip@  
 \penalty-50 \vskip-\parskip\smallskip  
 \DN@{\ifx\next\Runinitem\let\next@\relax
  \else\nextRunin@false\let\item\plainitem@  
   \ifx\next\par 
    \DN@\par{\everypar\expandafter{\the\everypartoks@}}%
   \else  
    \DN@{\noindent\everypar\expandafter{\the\everypartoks@}}%
  \fi\fi\next@}%
 \FN@\next@}
\newcount\rosterhangafter@
\def\Runinitem#1\roster\runinitem{\relaxnext@
 \rostercount@\z@ 
 \def\item{\FN@\rosteritem@}%
 \def\runinitem@{#1}%
 \DN@{\ifx\next[\let\next\nextii@\else\ifx\next"\let\next\nextiii@
  \else\let\next\nextiv@\fi\fi\next}%
 \DNii@[##1]{\rostercount@##1\relax
  \def\item@{{\rm(\number\rostercount@)}}\nextv@}%
 \def\nextiii@"##1"{\def\item@{{\rm##1}}\nextv@}%
 \def\nextiv@{\advance\rostercount@\@ne
  \def\item@{{\rm(\number\rostercount@)}}\nextv@}%
 \def\nextv@{\setbox\z@\vbox  
  {\ifnextRunin@\noindent\fi  
  \runinitem@\unskip\enspace\item@~\par  
  \global\rosterhangafter@\prevgraf}%
  \firstitem@false  
  \ifnextRunin@\else\par\fi  
  \hangafter\rosterhangafter@\hangindent3\parindent
  \ifnextRunin@\noindent\fi  
  \runinitem@\unskip\enspace 
  \item@~\ifnextRunin@\else\par@\fi  
  \nextRunin@true\ignorespaces}%
 \FN@\next@}
\def\footmarkform@#1{$\m@th^{#1}$}
\let\thefootnotemark\footmarkform@
\def\makefootnote@#1#2{\insert\footins
 {\interlinepenalty\interfootnotelinepenalty
 \eightpoint\splittopskip\ht\strutbox\splitmaxdepth\dp\strutbox
 \floatingpenalty\@MM\leftskip\z@\rightskip\z@\spaceskip\z@\xspaceskip\z@
 \leavevmode{#1}\footstrut\ignorespaces#2\unskip\lower\dp\strutbox
 \vbox to\dp\strutbox{}}}
\newcount\footmarkcount@
\footmarkcount@\z@
\def\footnotemark{\let\@sf\empty\relaxnext@
 \ifhmode\edef\@sf{\spacefactor\the\spacefactor}\/\fi
 \DN@{\ifx[\next\let\next@\nextii@\else
  \ifx"\next\let\next@\nextiii@\else
  \let\next@\nextiv@\fi\fi\next@}%
 \DNii@[##1]{\footmarkform@{##1}\@sf}%
 \def\nextiii@"##1"{{##1}\@sf}%
 \def\nextiv@{\iffirstchoice@\global\advance\footmarkcount@\@ne\fi
  \footmarkform@{\number\footmarkcount@}\@sf}%
 \FN@\next@}
\def\footnotetext{\relaxnext@
 \DN@{\ifx[\next\let\next@\nextii@\else
  \ifx"\next\let\next@\nextiii@\else
  \let\next@\nextiv@\fi\fi\next@}%
 \DNii@[##1]##2{\makefootnote@{\footmarkform@{##1}}{##2}}%
 \def\nextiii@"##1"##2{\makefootnote@{##1}{##2}}%
 \def\nextiv@##1{\makefootnote@{\footmarkform@{\number\footmarkcount@}}{##1}}%
 \FN@\next@}
\def\footnote{\let\@sf\empty\relaxnext@
 \ifhmode\edef\@sf{\spacefactor\the\spacefactor}\/\fi
 \DN@{\ifx[\next\let\next@\nextii@\else
  \ifx"\next\let\next@\nextiii@\else
  \let\next@\nextiv@\fi\fi\next@}%
 \DNii@[##1]##2{\footnotemark[##1]\footnotetext[##1]{##2}}%
 \def\nextiii@"##1"##2{\footnotemark"##1"\footnotetext"##1"{##2}}%
 \def\nextiv@##1{\footnotemark\footnotetext{##1}}%
 \FN@\next@}
\def\adjustfootnotemark#1{\advance\footmarkcount@#1\relax}
\def\footnoterule{\kern-3\p@
  \hrule width 5pc\kern 2.6\p@} 
\def\captionfont@{\smc}
\def\topcaption#1#2\endcaption{%
  {\dimen@\hsize \advance\dimen@-\captionwidth@
   \rm\raggedcenter@ \advance\leftskip.5\dimen@ \rightskip\leftskip
  {\captionfont@#1}%
  \if\notempty{#2}.\enspace\ignorespaces#2\fi
  \endgraf}\nobreak\bigskip}
\def\botcaption#1#2\endcaption{%
  \nobreak\bigskip
  \setboxz@h{\captionfont@#1\if\notempty{#2}.\enspace\rm#2\fi}%
  {\dimen@\hsize \advance\dimen@-\captionwidth@
   \leftskip.5\dimen@ \rightskip\leftskip
   \noindent \ifdim\wdz@>\captionwidth@ 
   \else\hfil\fi 
  {\captionfont@#1}\if\notempty{#2}.\enspace\rm#2\fi\endgraf}}
\def\@ins{\par\begingroup\def\vspace##1{\vskip##1\relax}%
  \def\captionwidth##1{\captionwidth@##1\relax}%
  \setbox\z@\vbox\bgroup} 
\def\block{\RIfMIfI@\nondmatherr@\block\fi
       \else\ifvmode\vskip\abovedisplayskip\noindent\fi
        $$\def\endblock{\par\egroup$$}\fi
  \vbox\bgroup\advance\hsize-2\indenti\noindent}
\def\endblock{\par\egroup}
\def\cite#1{{\rm[{\citefont@\m@th#1}]}}
\def\citefont@{\rm}
\def\refsfont@{\eightpoint}
\outer\def\Refs{\runaway@{proclaim}%
 \relaxnext@ \DN@{\ifx\next\nofrills\DN@\nofrills{\nextii@}\else
  \DN@{\nextii@{References}}\fi\next@}%
 \DNii@##1{\penaltyandskip@{-200}\aboveheadskip
  \line{\hfil\headfont@\ignorespaces##1\unskip\hfil}\nobreak
  \vskip\belowheadskip
  \begingroup\refsfont@\sfcode`.=\@m}%
 \FN@\next@}
\def\endRefs{\par\endgroup}
\newbox\nobox@            \newbox\keybox@           \newbox\bybox@
\newbox\paperbox@         \newbox\paperinfobox@     \newbox\jourbox@
\newbox\volbox@           \newbox\issuebox@         \newbox\yrbox@
\newbox\pagesbox@         \newbox\bookbox@          \newbox\bookinfobox@
\newbox\publbox@          \newbox\publaddrbox@      \newbox\finalinfobox@
\newbox\edsbox@           \newbox\langbox@
\newif\iffirstref@        \newif\iflastref@
\newif\ifprevjour@        \newif\ifbook@            \newif\ifprevinbook@
\newif\ifquotes@          \newif\ifbookquotes@      \newif\ifpaperquotes@
\newdimen\bysamerulewd@
\setboxz@h{\refsfont@\kern3em}
\bysamerulewd@\wdz@
\newdimen\refindentwd
\setboxz@h{\refsfont@ 00. }
\refindentwd\wdz@
\outer\def\ref{\begingroup \noindent\hangindent\refindentwd
 \firstref@true \def\nofrills{\def\refkern@{\kern3sp}}%
 \ref@}
\def\ref@{\book@false \bgroup\let\endrefitem@\egroup \ignorespaces}
\def\moreref{\endrefitem@\endref@\firstref@false\ref@}%
\def\transl{\endrefitem@\endref@\firstref@false
  \book@false
  \prepunct@
  \setboxz@h\bgroup \aftergroup\unhbox\aftergroup\z@
    \def\endrefitem@{\unskip\refkern@\egroup}\ignorespaces}%
\def\emptyifempty@{\dimen@\wd\currbox@
  \advance\dimen@-\wd\z@ \advance\dimen@-.1\p@
  \ifdim\dimen@<\z@ \setbox\currbox@\copy\voidb@x \fi}
\let\refkern@\relax
\def\endrefitem@{\unskip\refkern@\egroup
  \setboxz@h{\refkern@}\emptyifempty@}\ignorespaces
\def\refdef@#1#2#3{\edef\next@{\noexpand\endrefitem@
  \let\noexpand\currbox@\csname\expandafter\eat@\string#1box@\endcsname
    \noexpand\setbox\noexpand\currbox@\hbox\bgroup}%
  \toks@\expandafter{\next@}%
  \if\notempty{#2#3}\toks@\expandafter{\the\toks@
  \def\endrefitem@{\unskip#3\refkern@\egroup
  \setboxz@h{#2#3\refkern@}\emptyifempty@}#2}\fi
  \toks@\expandafter{\the\toks@\ignorespaces}%
  \edef#1{\the\toks@}}
\refdef@\no{}{. }
\refdef@\key{[\m@th}{] }
\refdef@\by{}{}
\def\bysame{\by\hbox to\bysamerulewd@{\hrulefill}\thinspace
   \kern0sp}
\def\manyby{\message{\string\manyby is no longer necessary; \string\by
  can be used instead, starting with version 2.0 of \styname.STY}\by}
\refdef@\paper{\ifpaperquotes@``\fi\it}{}
\refdef@\paperinfo{}{}
\def\jour{\endrefitem@\let\currbox@\jourbox@
  \setbox\currbox@\hbox\bgroup
  \def\endrefitem@{\unskip\refkern@\egroup
    \setboxz@h{\refkern@}\emptyifempty@
    \ifvoid\jourbox@\else\prevjour@true\fi}%
\ignorespaces}
\refdef@\vol{\ifbook@\else\bf\fi}{}
\refdef@\issue{no. }{}
\refdef@\yr{}{}
\refdef@\pages{}{}
\def\page{\endrefitem@\def\pp@{\def\pp@{pp.~}p.~}\let\currbox@\pagesbox@
  \setbox\currbox@\hbox\bgroup\ignorespaces}
\def\pp@{pp.~}
\def\book{\endrefitem@ \let\currbox@\bookbox@
 \setbox\currbox@\hbox\bgroup\def\endrefitem@{\unskip\refkern@\egroup
  \setboxz@h{\ifbookquotes@``\fi}\emptyifempty@
  \ifvoid\bookbox@\else\book@true\fi}%
  \ifbookquotes@``\fi\it\ignorespaces}
\def\inbook{\endrefitem@
  \let\currbox@\bookbox@\setbox\currbox@\hbox\bgroup
  \def\endrefitem@{\unskip\refkern@\egroup
  \setboxz@h{\ifbookquotes@``\fi}\emptyifempty@
  \ifvoid\bookbox@\else\book@true\previnbook@true\fi}%
  \ifbookquotes@``\fi\ignorespaces}
\refdef@\eds{(}{, eds.)}
\def\ed{\endrefitem@\let\currbox@\edsbox@
 \setbox\currbox@\hbox\bgroup
 \def\endrefitem@{\unskip, ed.)\refkern@\egroup
  \setboxz@h{(, ed.)}\emptyifempty@}(\ignorespaces}
\refdef@\bookinfo{}{}
\refdef@\publ{}{}
\refdef@\publaddr{}{}
\refdef@\finalinfo{}{}
\refdef@\lang{(}{)}

\let\refdef@\relax 
\def\ppunbox@#1{\ifvoid#1\else\prepunct@\unhbox#1\fi}
\def\nocomma@#1{\ifvoid#1\else\changepunct@3\prepunct@\unhbox#1\fi}
\def\changepunct@#1{\ifnum\lastkern<3 \unkern\kern#1sp\fi}
\def\prepunct@{\count@\lastkern\unkern
  \ifnum\lastpenalty=0
    \let\penalty@\relax
  \else
    \edef\penalty@{\penalty\the\lastpenalty\relax}%
  \fi
  \unpenalty
  \let\refspace@\ \ifcase\count@,
\or;\or.\or 
  \or\let\refspace@\relax
  \else,\fi
  \ifquotes@''\quotes@false\fi \penalty@ \refspace@
}
\def\transferpenalty@#1{\dimen@\lastkern\unkern
  \ifnum\lastpenalty=0\unpenalty\let\penalty@\relax
  \else\edef\penalty@{\penalty\the\lastpenalty\relax}\unpenalty\fi
  #1\penalty@\kern\dimen@}
\def\endref{\endrefitem@\lastref@true\endref@
  \par\endgroup \prevjour@false \previnbook@false }
\def\endref@{%
\iffirstref@
  \ifvoid\nobox@\ifvoid\keybox@\indent\fi
  \else\hbox to\refindentwd{\hss\unhbox\nobox@}\fi
  \ifvoid\keybox@
  \else\ifdim\wd\keybox@>\refindentwd
         \box\keybox@
       \else\hbox to\refindentwd{\unhbox\keybox@\hfil}\fi\fi
  \kern4sp\ppunbox@\bybox@
\fi 
  \ifvoid\paperbox@
  \else\prepunct@\unhbox\paperbox@
    \ifpaperquotes@\quotes@true\fi\fi
  \ppunbox@\paperinfobox@
  \ifvoid\jourbox@
    \ifprevjour@ \nocomma@\volbox@
      \nocomma@\issuebox@
      \ifvoid\yrbox@\else\changepunct@3\prepunct@(\unhbox\yrbox@
        \transferpenalty@)\fi
      \ppunbox@\pagesbox@
    \fi 
  \else \prepunct@\unhbox\jourbox@
    \nocomma@\volbox@
    \nocomma@\issuebox@
    \ifvoid\yrbox@\else\changepunct@3\prepunct@(\unhbox\yrbox@
      \transferpenalty@)\fi
    \ppunbox@\pagesbox@
  \fi 
  \ifbook@\prepunct@\unhbox\bookbox@ \ifbookquotes@\quotes@true\fi \fi
  \nocomma@\edsbox@
  \ppunbox@\bookinfobox@
  \ifbook@\ifvoid\volbox@\else\prepunct@ vol.~\unhbox\volbox@
  \fi\fi
  \ppunbox@\publbox@ \ppunbox@\publaddrbox@
  \ifbook@ \ppunbox@\yrbox@
    \ifvoid\pagesbox@
    \else\prepunct@\pp@\unhbox\pagesbox@\fi
  \else
    \ifprevinbook@ \ppunbox@\yrbox@
      \ifvoid\pagesbox@\else\prepunct@\pp@\unhbox\pagesbox@\fi
    \fi \fi
  \ppunbox@\finalinfobox@
  \iflastref@
    \ifvoid\langbox@.\ifquotes@''\fi
    \else\changepunct@2\prepunct@\unhbox\langbox@\fi
  \else
    \ifvoid\langbox@\changepunct@1%
    \else\changepunct@3\prepunct@\unhbox\langbox@
      \changepunct@1\fi
  \fi
}
\outer\def\enddocument{%
 \runaway@{proclaim}%
\ifmonograph@ 
\else
 \nobreak
 \thetranslator@
 \count@\z@ \loop\ifnum\count@<\addresscount@\advance\count@\@ne
 \csname address\number\count@\endcsname
 \csname email\number\count@\endcsname
 \repeat
\fi
 \vfill\supereject\end}

\def\headfont@{\headfonts}
\def\proclaimheadfont@{\bf}
\def\specialheadfont@{\bf}
\def\subheadfont@{\bf}
\def\demoheadfont@{\smc}

\newif\ifThisToToc \ThisToTocfalse
\newif\iftocloaded \tocloadedfalse

\def\C@L{\noexpand\Cal}\def\B@B{\noexpand\Bbb}\def\fR@K{\noexpand\frak}
\def\S@{\noexpand\S}\def\P@P{\noexpand\"}
\def\xpar{\\}

\def\writetoc#1{\iftocloaded\ifThisToToc\begingroup\def\totoc{}
  \def\Cal{\noexpand\C@L}\def\Bbb{\noexpand\B@B}
  \def\frak{\noexpand\fR@K}\def\goth{\frak}\def\S{\noexpand\S@}
  \def\"{\noexpand\P@P}
  \def\xpar{\par\penalty100000 }\def\idx##1{##1}\def\\{\xpar}
  \edef\next@{\write\toc{\noindent#1\leaderfill\noexpand\folio\par}}%
  \next@\endgroup\global\ThisToTocfalse\fi\fi}
\def\leaderfill{\leaders\hbox to 1em{\hss.\hss}\hfill}

\newif\ifindexloaded \indexloadedfalse
\def\idx#1{\ifindexloaded\begingroup\def\ign{}\def\it{}\def\/{}%
 \def\smc{}\def\bf{}\def\tt{}%
 \def\Cal{\noexpand\C@L}\def\Bbb{\noexpand\B@B}%
 \def\frak{\noexpand\fR@K}\def\goth{\frak}\def\S{\noexpand\S@}%
  \def\"{\noexpand\P@P}%
 {\edef\next@{\write\index{#1, \noexpand\folio}}\next@}%
 \endgroup\fi{#1}}
\def\ign#1{}

\def\totoc{\global\ThisToToctrue}

\outer\def\head#1\endhead{\par\penaltyandskip@{-200}\aboveheadskip
 {\headfont@\raggedcenter@\interlinepenalty\@M
 \ignorespaces#1\endgraf}\nobreak
 \vskip\belowheadskip
 \headmark{#1}\writetoc{#1}}

\outer\def\chaphead#1\endchaphead{\par\penaltyandskip@{-200}\aboveheadskip
 {\chapheadfonts\raggedcenter@\interlinepenalty\@M
 \ignorespaces#1\endgraf}\nobreak
 \vskip3\belowheadskip
 \headmark{#1}\writetoc{#1}}

\def\folio{{\foliofont@\ifnum\pageno<\z@ \romannumeral-\pageno
 \else\number\pageno \fi}}
\newtoks\leftheadtoks
\newtoks\rightheadtoks

\def\leftheadtext{\nofrills@{\relax}\lht@
  \DNii@##1{\leftheadtoks\expandafter{\lht@{##1}}%
    \mark{\the\leftheadtoks\noexpand\else\the\rightheadtoks}
    \ifsyntax@\setboxz@h{\def\\{\unskip\space\ignorespaces}%
        \headlinefont@##1}\fi}%
  \FN@\next@}
\def\rightheadtext{\nofrills@{\relax}\rht@
  \DNii@##1{\rightheadtoks\expandafter{\rht@{##1}}%
    \mark{\the\leftheadtoks\noexpand\else\the\rightheadtoks}%
    \ifsyntax@\setboxz@h{\def\\{\unskip\space\ignorespaces}%
        \headlinefont@##1}\fi}%
  \FN@\next@}
\def\NoRunningHeads{\global\runheads@false\global\let\headmark\eat@}

\newif\iffirstpage@     \firstpage@true
\newif\ifrunheads@      \runheads@true

\newdimen\fullhsize \fullhsize=\hsize
\newdimen\fullvsize \fullvsize=\vsize
\def\fullline{\hbox to\fullhsize}

\def\pagenumbers{\gdef\folio{\folio@}}

\let\norunningheads\NoRunningHeads
\def\userunningheads{\global\runheads@true}
\norunningheads

\headline={\def\chapter#1{\chapterno@. }%
  \def\\{\unskip\space\ignorespaces}\ifrunheads@\headlinefont@
    \ifodd\pageno\rightheadline \else\leftheadline\fi
   \else\hfil\fi\ifNoRunHeadline\global\NoRunHeadlinefalse\fi}
\let\folio@\folio
\def\foliofont@{\foliofont}
\def\foliofont{\eightrm}
\def\headlinefont@{\headlinefont}
\def\headlinefont{\eightpoint\smc}
\def\leftheadline{\rlap{\folio}\hfill
   \ifNoRunHeadline\else\iftrue\topmark\fi\fi \hfill}
\def\rightheadline{\hfill\ifNoRunHeadline
   \else \expandafter\fi
  \hfill \llap{\folio}}
\footline={{\eightpoint\bottremark}%
   \ifrunheads@\else\hfil{\let\foliofont\tenrm\folio}\fi\hfil}
\def\bottremark{}
 
\newif\ifNoRunHeadline      
\def\norunninghead{\global\NoRunHeadlinetrue}
\norunninghead

\output={\output@}
%
\newif\ifoffset\offsetfalse
\output={\output@}
\def\output@{%
 \ifoffset 
  \ifodd\count0\advance\hoffset by0.5truecm
   \else\advance\hoffset by-0.5truecm\fi\fi
 \shipout\vbox{%
  \makeheadline \pagebody \makefootline }%
 \advancepageno \ifnum\outputpenalty>-\@MM\else\dosupereject\fi}

\def\indexoutput#1{%
  \ifoffset 
   \ifodd\count0\advance\hoffset by0.5truecm
    \else\advance\hoffset by-0.5truecm\fi\fi
  \shipout\vbox{\makeheadline
  \vbox to\fullvsize{\boxmaxdepth\maxdepth%
  \ifvoid\topins\else\unvbox\topins\fi%
  #1 %
  \ifvoid\footins\else 
    \vskip\skip\footins
    \footnoterule
    \unvbox\footins\fi
  \ifr@ggedbottom \kern-\dimen@ \vfil \fi}%
  \baselineskip2pc
  \makefootline}%
 \global\advance\pageno\@ne
 \ifnum\outputpenalty>-\@MM\else\dosupereject\fi}
 
 \newbox\partialpage \newdimen\halfsize \halfsize=0.5\fullhsize
 \advance\halfsize by-0.5em

 \def\begindoublecolumns{\output={\indexoutput{\unvbox255}}%
   \begingroup \def\line{\fullline}
   \output={\global\setbox\partialpage=\vbox{\unvbox255\bigskip}}\eject
   \output={\doublecolumnout}\hsize=\halfsize \vsize=2\fullvsize}
 \def\enddoublecolumns{\output={\balancecolumns}\eject
  \endgroup \pagegoal=\fullvsize%
  \output={\output@}}
\def\doublecolumnout{\splittopskip=\topskip \splitmaxdepth=\maxdepth
  \dimen@=\fullvsize \advance\dimen@ by-\ht\partialpage
  \setbox0=\vsplit255 to \dimen@ \setbox2=\vsplit255 to \dimen@
  \indexoutput{\pagesofar} \unvbox255 \penalty\outputpenalty}
\def\pagesofar{\unvbox\partialpage
  \wd0=\hsize \wd2=\hsize \hbox to\fullhsize{\box0\hfil\box2}}
\def\balancecolumns{\setbox0=\vbox{\unvbox255} \dimen@=\ht0
  \advance\dimen@ by\topskip \advance\dimen@ by-\baselineskip
  \divide\dimen@ by2 \splittopskip=\topskip
  {\vbadness=10000 \loop \global\setbox3=\copy0
    \global\setbox1=\vsplit3 to\dimen@
    \ifdim\ht3>\dimen@ \global\advance\dimen@ by1pt \repeat}
  \setbox0=\vbox to\dimen@{\unvbox1} \setbox2=\vbox to\dimen@{\unvbox3}
  \pagesofar}

\tenpoint
\catcode`\@=\active

\def\smallheadings{\let\chapheadfonts\tenpoint\let\headfonts\tenpoint}

\tenpoint
\catcode`\@=\active

\def\LL{\leavevmode\setbox0=\hbox{L}\hbox to\wd0{\hss\char'40L}}
\def\al{\alpha}

\def\la{\lambda}
\def\rh{\rho}


\def\Z{{\Bbb Z}}

\def\today{\ifcase\month\or
 January\or February\or March\or April\or May\or June\or
 July\or August\or September\or October\or November\or December\fi
 \space\number\day, \number\year}

\def\({\left(}
\def\){\right)}
\def\[{\left[}
\def\]{\right]}

\def\3{\ss}
\catcode`\@=11
\def\dddot#1{\vbox{\ialign{##\crcr
      .\hskip-.5pt.\hskip-.5pt.\crcr\noalign{\kern1.5\p@\nointerlineskip}
      $\hfil\displaystyle{#1}\hfil$\crcr}}}

\newif\iftab@\tab@false
\newif\ifvtab@\vtab@false
\def\tab{\bgroup\tab@true\vtab@false\vst@bfalse\Strich@false%
   \def\\{\global\hline@@false%
     \ifhline@\global\hline@false\global\hline@@true\fi\cr}
   \edef\l@{\the\leftskip}\ialign\bgroup\hskip\l@##\hfil&&##\hfil\cr}
\def\endtab{\cr\egroup\egroup}
\def\vtab{\vtop\bgroup\vst@bfalse\vtab@true\tab@true\Strich@false%
   \bgroup\def\\{\cr}\ialign\bgroup&##\hfil\cr}
\def\endvtab{\cr\egroup\egroup\egroup}
\def\stab{\D@cke0.5pt\null 
 \bgroup\tab@true\vtab@false\vst@bfalse\Strich@true\Let@@\vspace@
 \normalbaselines\offinterlineskip
  \openup\spreadmlines@
 \edef\l@{\the\leftskip}\ialign
 \bgroup\hskip\l@##\hfil&&##\hfil\crcr}
\def\endstab{\crcr\egroup
 \egroup}
\newif\ifvst@b\vst@bfalse
\def\vstab{\D@cke0.5pt\null
 \vtop\bgroup\tab@true\vtab@false\vst@btrue\Strich@true\bgroup\Let@@\vspace@
 \normalbaselines\offinterlineskip
  \openup\spreadmlines@\bgroup}
\def\endvstab{\crcr\egroup\egroup
 \egroup\tab@false\Strich@false}

\newdimen\htstrut@
\htstrut@8.5\p@
\newdimen\htStrut@
\htStrut@12\p@
\newdimen\dpstrut@
\dpstrut@3.5\p@
\newdimen\dpStrut@
\dpStrut@3.5\p@
\def\openup{\afterassignment\@penup\dimen@=}
\def\@penup{\advance\lineskip\dimen@
  \advance\baselineskip\dimen@
  \advance\lineskiplimit\dimen@
  \divide\dimen@ by2
  \advance\htstrut@\dimen@
  \advance\htStrut@\dimen@
  \advance\dpstrut@\dimen@
  \advance\dpStrut@\dimen@}
\def\Let@@{\relax%
    \def\\{\global\hline@@false%
     \ifhline@\global\hline@false\global\hline@@true\fi\cr}%
    \iffalse}\fi}
\def\matrix{\null\,\vcenter\bgroup
 \tab@false\vtab@false\vst@bfalse\Strich@false\Let@@\vspace@
 \normalbaselines\openup\spreadmlines@\ialign
 \bgroup\hfil$\m@th##$\hfil&&\quad\hfil$\m@th##$\hfil\crcr
 \Mathstrut@\crcr\noalign{\kern-\baselineskip}}
\def\endmatrix{\crcr\Mathstrut@\crcr\noalign{\kern-\baselineskip}\egroup
 \egroup\,}
\def\smatrix{\D@cke0.5pt\null\,
 \vcenter\bgroup\tab@false\vtab@false\vst@bfalse\Strich@true\Let@@\vspace@
 \normalbaselines\offinterlineskip
  \openup\spreadmlines@\ialign
 \bgroup\hfil$\m@th##$\hfil&&\quad\hfil$\m@th##$\hfil\crcr}
\def\endsmatrix{\crcr\egroup
 \egroup\,\Strich@false}
\newdimen\D@cke
\def\Dicke#1{\global\D@cke#1}
\newtoks\tabs@\tabs@{&}
\newif\ifStrich@\Strich@false
\newif\iff@rst

\def\Stricherr@{\iftab@\ifvtab@\errmessage{\noexpand\s not allowed
     here. Use \noexpand\vstab!}%
  \else\errmessage{\noexpand\s not allowed here. Use \noexpand\stab!}%
  \fi\else\errmessage{\noexpand\s not allowed
     here. Use \noexpand\smatrix!}\fi}
\def\format{\ifvst@b\else\crcr\fi\egroup\iffalse{\fi\ifnum`}=0 \fi\format@}
\def\format@#1\\{\def\preamble@{#1}%
 \def\Str@chfehlt##1{\ifx##1\s\Stricherr@\fi\ifx##1\\\let\Next\relax%
   \else\let\Next\Str@chfehlt\fi\Next}%
 \def\c{\hfil\noexpand\ifhline@@\hbox{\vrule height\htStrut@%
   depth\dpstrut@ width\z@}\noexpand\fi%
   \ifStrich@\hbox{\vrule height\htstrut@ depth\dpstrut@ width\z@}%
   \fi\iftab@\else$\m@th\fi\the\hashtoks@\iftab@\else$\fi\hfil}%
 \def\r{\hfil\noexpand\ifhline@@\hbox{\vrule height\htStrut@%
   depth\dpstrut@ width\z@}\noexpand\fi%
   \ifStrich@\hbox{\vrule height\htstrut@ depth\dpstrut@ width\z@}%
   \fi\iftab@\else$\m@th\fi\the\hashtoks@\iftab@\else$\fi}%
 \def\l{\noexpand\ifhline@@\hbox{\vrule height\htStrut@%
   depth\dpstrut@ width\z@}\noexpand\fi%
   \ifStrich@\hbox{\vrule height\htstrut@ depth\dpstrut@ width\z@}%
   \fi\iftab@\else$\m@th\fi\the\hashtoks@\iftab@\else$\fi\hfil}%
 \def\s{\ifStrich@\ \the\tabs@\vrule width\D@cke\the\hashtoks@%
          \fi\the\tabs@\ }%
 \def\sa{\ifStrich@\vrule width\D@cke\the\hashtoks@%
            \the\tabs@\ %
            \fi}%
 \def\se{\ifStrich@\ \the\tabs@\vrule width\D@cke\the\hashtoks@\fi}%
 \def\cd{\hfil\noexpand\ifhline@@\hbox{\vrule height\htStrut@%
   depth\dpstrut@ width\z@}\noexpand\fi%
   \ifStrich@\hbox{\vrule height\htstrut@ depth\dpstrut@ width\z@}%
   \fi$\dsize\m@th\the\hashtoks@$\hfil}%
 \def\rd{\hfil\noexpand\ifhline@@\hbox{\vrule height\htStrut@%
   depth\dpstrut@ width\z@}\noexpand\fi%
   \ifStrich@\hbox{\vrule height\htstrut@ depth\dpstrut@ width\z@}%
   \fi$\dsize\m@th\the\hashtoks@$}%
 \def\ld{\noexpand\ifhline@@\hbox{\vrule height\htStrut@%
   depth\dpstrut@ width\z@}\noexpand\fi%
   \ifStrich@\hbox{\vrule height\htstrut@ depth\dpstrut@ width\z@}%
   \fi$\dsize\m@th\the\hashtoks@$\hfil}%
 \ifStrich@\else\Str@chfehlt#1\\\fi%
 \setbox\z@\hbox{\xdef\Preamble@{\preamble@}}\ifnum`{=0 \fi\iffalse}\fi
 \ialign\bgroup\span\Preamble@\crcr}
\newif\ifhline@\hline@false
\newif\ifhline@@\hline@@false
\def\hlinefor#1{\multispan@{\strip@#1 }\leaders\hrule height\D@cke\hfill%
    \global\hline@true\ignorespaces}
\def\Item "#1"{\par\noindent\hangindent2\parindent%
  \hangafter1\setbox0\hbox{\rm#1\enspace}\ifdim\wd0>2\parindent%
  \box0\else\hbox to 2\parindent{\rm#1\hfil}\fi\ignorespaces}
\def\ITEM #1"#2"{\par\noindent\hangafter1\hangindent#1%
  \setbox0\hbox{\rm#2\enspace}\ifdim\wd0>#1%
  \box0\else\hbox to 0pt{\rm#2\hss}\hskip#1\fi\ignorespaces}
\def\item"#1"{\par\noindent\hang%
  \setbox0=\hbox{\rm#1\enspace}\ifdim\wd0>\the\parindent%
  \box0\else\hbox to \parindent{\rm#1\hfil}\enspace\fi\ignorespaces}
\let\plainitem@\item
\catcode`\@=13

\hsize13cm
\vsize19cm
\newdimen\fullhsize
\newdimen\fullvsize
\newdimen\halfsize
\fullhsize13cm
\fullvsize19cm
\halfsize=0.5\fullhsize
\advance\halfsize by-0.5em

\magnification1200

\catcode`\@=11
\font\tenln    = line10
\font\tenlnw   = linew10

\newskip\Einheit \Einheit=0.5cm
\newcount\xcoord \newcount\ycoord
\newdimen\xdim \newdimen\ydim \newdimen\PfadD@cke \newdimen\Pfadd@cke

\newcount\@tempcnta
\newcount\@tempcntb

\newdimen\@tempdima
\newdimen\@tempdimb

\newdimen\@wholewidth
\newdimen\@halfwidth

\newcount\@xarg
\newcount\@yarg
\newcount\@yyarg
\newbox\@linechar
\newbox\@tempboxa
\newdimen\@linelen
\newdimen\@clnwd
\newdimen\@clnht

\newif\if@negarg

\def\@whilenoop#1{}
\def\@whiledim#1\do #2{\ifdim #1\relax#2\@iwhiledim{#1\relax#2}\fi}
\def\@iwhiledim#1{\ifdim #1\let\@nextwhile=\@iwhiledim
        \else\let\@nextwhile=\@whilenoop\fi\@nextwhile{#1}}

\def\@whileswnoop#1\fi{}
\def\@whilesw#1\fi#2{#1#2\@iwhilesw{#1#2}\fi\fi}
\def\@iwhilesw#1\fi{#1\let\@nextwhile=\@iwhilesw
         \else\let\@nextwhile=\@whileswnoop\fi\@nextwhile{#1}\fi}

\def\thinlines{\let\@linefnt\tenln \let\@circlefnt\tencirc
  \@wholewidth\fontdimen8\tenln \@halfwidth .5\@wholewidth}
\def\thicklines{\let\@linefnt\tenlnw \let\@circlefnt\tencircw
  \@wholewidth\fontdimen8\tenlnw \@halfwidth .5\@wholewidth}
\thinlines

\PfadD@cke1pt \Pfadd@cke0.5pt
\def\PfadDicke#1{\PfadD@cke#1 \divide\PfadD@cke by2 \Pfadd@cke\PfadD@cke \multiply\PfadD@cke by2}
\long\def\LOOP#1\REPEAT{\def\BODY{#1}\ITERATE}
\def\ITERATE{\BODY \let\next\ITERATE \else\let\next\relax\fi \next}
\let\REPEAT=\fi
\def\Punkt{\hbox{\raise-2pt\hbox to0pt{\hss$\ssize\bullet$\hss}}}
\def\DuennPunkt(#1,#2){\unskip
  \raise#2 \Einheit\hbox to0pt{\hskip#1 \Einheit
          \raise-2.5pt\hbox to0pt{\hss$\bullet$\hss}\hss}}
\def\NormalPunkt(#1,#2){\unskip
  \raise#2 \Einheit\hbox to0pt{\hskip#1 \Einheit
          \raise-3pt\hbox to0pt{\hss\twelvepoint$\bullet$\hss}\hss}}
\def\DickPunkt(#1,#2){\unskip
  \raise#2 \Einheit\hbox to0pt{\hskip#1 \Einheit
          \raise-4pt\hbox to0pt{\hss\fourteenpoint$\bullet$\hss}\hss}}
\def\Kreis(#1,#2){\unskip
  \raise#2 \Einheit\hbox to0pt{\hskip#1 \Einheit
          \raise-4pt\hbox to0pt{\hss\fourteenpoint$\circ$\hss}\hss}}

\def\Line@(#1,#2)#3{\@xarg #1\relax \@yarg #2\relax
\@linelen=#3\Einheit
\ifnum\@xarg =0 \@vline
  \else \ifnum\@yarg =0 \@hline \else \@sline\fi
\fi}

\def\@sline{\ifnum\@xarg< 0 \@negargtrue \@xarg -\@xarg \@yyarg -\@yarg
  \else \@negargfalse \@yyarg \@yarg \fi
\ifnum \@yyarg >0 \@tempcnta\@yyarg \else \@tempcnta -\@yyarg \fi
\ifnum\@tempcnta>6 \@badlinearg\@tempcnta0 \fi
\ifnum\@xarg>6 \@badlinearg\@xarg 1 \fi
\setbox\@linechar\hbox{\@linefnt\@getlinechar(\@xarg,\@yyarg)}%
\ifnum \@yarg >0 \let\@upordown\raise \@clnht\z@
   \else\let\@upordown\lower \@clnht \ht\@linechar\fi
\@clnwd=\wd\@linechar
\if@negarg \hskip -\wd\@linechar \def\@tempa{\hskip -2\wd\@linechar}\else
     \let\@tempa\relax \fi
\@whiledim \@clnwd <\@linelen \do
  {\@upordown\@clnht\copy\@linechar
   \@tempa
   \advance\@clnht \ht\@linechar
   \advance\@clnwd \wd\@linechar}%
\advance\@clnht -\ht\@linechar
\advance\@clnwd -\wd\@linechar
\@tempdima\@linelen\advance\@tempdima -\@clnwd
\@tempdimb\@tempdima\advance\@tempdimb -\wd\@linechar
\if@negarg \hskip -\@tempdimb \else \hskip \@tempdimb \fi
\multiply\@tempdima \@m
\@tempcnta \@tempdima \@tempdima \wd\@linechar \divide\@tempcnta \@tempdima
\@tempdima \ht\@linechar \multiply\@tempdima \@tempcnta
\divide\@tempdima \@m
\advance\@clnht \@tempdima
\ifdim \@linelen <\wd\@linechar
   \hskip \wd\@linechar
  \else\@upordown\@clnht\copy\@linechar\fi}

\def\@hline{\ifnum \@xarg <0 \hskip -\@linelen \fi
\vrule height\Pfadd@cke width \@linelen depth\Pfadd@cke
\ifnum \@xarg <0 \hskip -\@linelen \fi}

\def\@getlinechar(#1,#2){\@tempcnta#1\relax\multiply\@tempcnta 8
\advance\@tempcnta -9 \ifnum #2>0 \advance\@tempcnta #2\relax\else
\advance\@tempcnta -#2\relax\advance\@tempcnta 64 \fi
\char\@tempcnta}

\def\Vektor(#1,#2)#3(#4,#5){\unskip\leavevmode
  \xcoord#4\relax \ycoord#5\relax
      \raise\ycoord \Einheit\hbox to0pt{\hskip\xcoord \Einheit
         \Vector@(#1,#2){#3}\hss}}

\def\Vector@(#1,#2)#3{\@xarg #1\relax \@yarg #2\relax
\@tempcnta \ifnum\@xarg<0 -\@xarg\else\@xarg\fi
\ifnum\@tempcnta<5\relax
\@linelen=#3\Einheit
\ifnum\@xarg =0 \@vvector
  \else \ifnum\@yarg =0 \@hvector \else \@svector\fi
\fi
\else\@badlinearg\fi}

\def\@hvector{\@hline\hbox to 0pt{\@linefnt
\ifnum \@xarg <0 \@getlarrow(1,0)\hss\else
    \hss\@getrarrow(1,0)\fi}}

\def\@vvector{\ifnum \@yarg <0 \@downvector \else \@upvector \fi}

\def\@svector{\@sline
\@tempcnta\@yarg \ifnum\@tempcnta <0 \@tempcnta=-\@tempcnta\fi
\ifnum\@tempcnta <5
  \hskip -\wd\@linechar
  \@upordown\@clnht \hbox{\@linefnt  \if@negarg
  \@getlarrow(\@xarg,\@yyarg) \else \@getrarrow(\@xarg,\@yyarg) \fi}%
\else\@badlinearg\fi}

\def\@upline{\hbox to \z@{\hskip -.5\Pfadd@cke \vrule width \Pfadd@cke
   height \@linelen depth \z@\hss}}

\def\@downline{\hbox to \z@{\hskip -.5\Pfadd@cke \vrule width \Pfadd@cke
   height \z@ depth \@linelen \hss}}

\def\@upvector{\@upline\setbox\@tempboxa\hbox{\@linefnt\char'66}\raise
     \@linelen \hbox to\z@{\lower \ht\@tempboxa\box\@tempboxa\hss}}

\def\@downvector{\@downline\lower \@linelen
      \hbox to \z@{\@linefnt\char'77\hss}}

\def\@getlarrow(#1,#2){\ifnum #2 =\z@ \@tempcnta='33\else
\@tempcnta=#1\relax\multiply\@tempcnta \sixt@@n \advance\@tempcnta
-9 \@tempcntb=#2\relax\multiply\@tempcntb \tw@
\ifnum \@tempcntb >0 \advance\@tempcnta \@tempcntb\relax
\else\advance\@tempcnta -\@tempcntb\advance\@tempcnta 64
\fi\fi\char\@tempcnta}

\def\@getrarrow(#1,#2){\@tempcntb=#2\relax
\ifnum\@tempcntb < 0 \@tempcntb=-\@tempcntb\relax\fi
\ifcase \@tempcntb\relax \@tempcnta='55 \or
\ifnum #1<3 \@tempcnta=#1\relax\multiply\@tempcnta
24 \advance\@tempcnta -6 \else \ifnum #1=3 \@tempcnta=49
\else\@tempcnta=58 \fi\fi\or
\ifnum #1<3 \@tempcnta=#1\relax\multiply\@tempcnta
24 \advance\@tempcnta -3 \else \@tempcnta=51\fi\or
\@tempcnta=#1\relax\multiply\@tempcnta
\sixt@@n \advance\@tempcnta -\tw@ \else
\@tempcnta=#1\relax\multiply\@tempcnta
\sixt@@n \advance\@tempcnta 7 \fi\ifnum #2<0 \advance\@tempcnta 64 \fi
\char\@tempcnta}

\def\Diagonale(#1,#2)#3{\unskip\leavevmode
  \xcoord#1\relax \ycoord#2\relax
      \raise\ycoord \Einheit\hbox to0pt{\hskip\xcoord \Einheit
         \Line@(1,1){#3}\hss}}
\def\AntiDiagonale(#1,#2)#3{\unskip\leavevmode
  \xcoord#1\relax \ycoord#2\relax 
      \raise\ycoord \Einheit\hbox to0pt{\hskip\xcoord \Einheit
         \Line@(1,-1){#3}\hss}}
\def\Pfad(#1,#2),#3\endPfad{\unskip\leavevmode
  \xcoord#1 \ycoord#2 \thicklines\ZeichnePfad#3\endPfad\thinlines}
\def\ZeichnePfad#1{\ifx#1\endPfad\let\next\relax
  \else\let\next\ZeichnePfad
    \ifnum#1=1
      \raise\ycoord \Einheit\hbox to0pt{\hskip\xcoord \Einheit
         \vrule height\Pfadd@cke width1 \Einheit depth\Pfadd@cke\hss}%
      \advance\xcoord by 1
    \else\ifnum#1=2
      \raise\ycoord \Einheit\hbox to0pt{\hskip\xcoord \Einheit
        \hbox{\hskip-\PfadD@cke\vrule height1 \Einheit width\PfadD@cke depth0pt}\hss}%
      \advance\ycoord by 1
    \else\ifnum#1=3
      \raise\ycoord \Einheit\hbox to0pt{\hskip\xcoord \Einheit
         \Line@(1,1){1}\hss}
      \advance\xcoord by 1
      \advance\ycoord by 1
    \else\ifnum#1=4
      \raise\ycoord \Einheit\hbox to0pt{\hskip\xcoord \Einheit
         \Line@(1,-1){1}\hss}
      \advance\xcoord by 1
      \advance\ycoord by -1
    \else\ifnum#1=5
      \advance\xcoord by -1
      \raise\ycoord \Einheit\hbox to0pt{\hskip\xcoord \Einheit
         \vrule height\Pfadd@cke width1 \Einheit depth\Pfadd@cke\hss}%
    \else\ifnum#1=6
      \advance\ycoord by -1
      \raise\ycoord \Einheit\hbox to0pt{\hskip\xcoord \Einheit
        \hbox{\hskip-\PfadD@cke\vrule height1 \Einheit width\PfadD@cke depth0pt}\hss}%
    \else\ifnum#1=7
      \advance\xcoord by -1
      \advance\ycoord by -1
      \raise\ycoord \Einheit\hbox to0pt{\hskip\xcoord \Einheit
         \Line@(1,1){1}\hss}
    \else\ifnum#1=8
      \advance\xcoord by -1
      \advance\ycoord by +1
      \raise\ycoord \Einheit\hbox to0pt{\hskip\xcoord \Einheit
         \Line@(1,-1){1}\hss}
    \fi\fi\fi\fi
    \fi\fi\fi\fi
  \fi\next}
\def\hSSchritt{\leavevmode\raise-.4pt\hbox to0pt{\hss.\hss}\hskip.2\Einheit
  \raise-.4pt\hbox to0pt{\hss.\hss}\hskip.2\Einheit
  \raise-.4pt\hbox to0pt{\hss.\hss}\hskip.2\Einheit
  \raise-.4pt\hbox to0pt{\hss.\hss}\hskip.2\Einheit
  \raise-.4pt\hbox to0pt{\hss.\hss}\hskip.2\Einheit}
\def\vSSchritt{\vbox{\baselineskip.2\Einheit\lineskiplimit0pt
\hbox{.}\hbox{.}\hbox{.}\hbox{.}\hbox{.}}}
\def\DSSchritt{\leavevmode\raise-.4pt\hbox to0pt{%
  \hbox to0pt{\hss.\hss}\hskip.2\Einheit
  \raise.2\Einheit\hbox to0pt{\hss.\hss}\hskip.2\Einheit
  \raise.4\Einheit\hbox to0pt{\hss.\hss}\hskip.2\Einheit
  \raise.6\Einheit\hbox to0pt{\hss.\hss}\hskip.2\Einheit
  \raise.8\Einheit\hbox to0pt{\hss.\hss}\hss}}
\def\dSSchritt{\leavevmode\raise-.4pt\hbox to0pt{%
  \hbox to0pt{\hss.\hss}\hskip.2\Einheit
  \raise-.2\Einheit\hbox to0pt{\hss.\hss}\hskip.2\Einheit
  \raise-.4\Einheit\hbox to0pt{\hss.\hss}\hskip.2\Einheit
  \raise-.6\Einheit\hbox to0pt{\hss.\hss}\hskip.2\Einheit
  \raise-.8\Einheit\hbox to0pt{\hss.\hss}\hss}}
\def\SPfad(#1,#2),#3\endSPfad{\unskip\leavevmode
  \xcoord#1 \ycoord#2 \ZeichneSPfad#3\endSPfad}
\def\ZeichneSPfad#1{\ifx#1\endSPfad\let\next\relax
  \else\let\next\ZeichneSPfad
    \ifnum#1=1
      \raise\ycoord \Einheit\hbox to0pt{\hskip\xcoord \Einheit
         \hSSchritt\hss}%
      \advance\xcoord by 1
    \else\ifnum#1=2
      \raise\ycoord \Einheit\hbox to0pt{\hskip\xcoord \Einheit
        \hbox{\hskip-2pt \vSSchritt}\hss}%
      \advance\ycoord by 1
    \else\ifnum#1=3
      \raise\ycoord \Einheit\hbox to0pt{\hskip\xcoord \Einheit
         \DSSchritt\hss}
      \advance\xcoord by 1
      \advance\ycoord by 1
    \else\ifnum#1=4
      \raise\ycoord \Einheit\hbox to0pt{\hskip\xcoord \Einheit
         \dSSchritt\hss}
      \advance\xcoord by 1
      \advance\ycoord by -1
    \else\ifnum#1=5
      \advance\xcoord by -1
      \raise\ycoord \Einheit\hbox to0pt{\hskip\xcoord \Einheit
         \hSSchritt\hss}%
    \else\ifnum#1=6
      \advance\ycoord by -1
      \raise\ycoord \Einheit\hbox to0pt{\hskip\xcoord \Einheit
        \hbox{\hskip-2pt \vSSchritt}\hss}%
    \else\ifnum#1=7
      \advance\xcoord by -1
      \advance\ycoord by -1
      \raise\ycoord \Einheit\hbox to0pt{\hskip\xcoord \Einheit
         \DSSchritt\hss}
    \else\ifnum#1=8
      \advance\xcoord by -1
      \advance\ycoord by 1
      \raise\ycoord \Einheit\hbox to0pt{\hskip\xcoord \Einheit
         \dSSchritt\hss}
    \fi\fi\fi\fi
    \fi\fi\fi\fi
  \fi\next}
\def\Koordinatenachsen(#1,#2){\unskip
 \hbox to0pt{\hskip-.5pt\vrule height#2 \Einheit width.5pt depth1 \Einheit}%
 \hbox to0pt{\hskip-1 \Einheit \xcoord#1 \advance\xcoord by1
    \vrule height0.25pt width\xcoord \Einheit depth0.25pt\hss}}
\def\Koordinatenachsen(#1,#2)(#3,#4){\unskip
 \hbox to0pt{\hskip-.5pt \ycoord-#4 \advance\ycoord by1
    \vrule height#2 \Einheit width.5pt depth\ycoord \Einheit}%
 \hbox to0pt{\hskip-1 \Einheit \hskip#3\Einheit 
    \xcoord#1 \advance\xcoord by1 \advance\xcoord by-#3 
    \vrule height0.25pt width\xcoord \Einheit depth0.25pt\hss}}
\def\Gitter(#1,#2){\unskip \xcoord0 \ycoord0 \leavevmode
  \LOOP\ifnum\ycoord<#2
    \loop\ifnum\xcoord<#1
      \raise\ycoord \Einheit\hbox to0pt{\hskip\xcoord \Einheit\Punkt\hss}%
      \advance\xcoord by1
    \repeat
    \xcoord0
    \advance\ycoord by1
  \REPEAT}
\def\Gitter(#1,#2)(#3,#4){\unskip \xcoord#3 \ycoord#4 \leavevmode
  \LOOP\ifnum\ycoord<#2
    \loop\ifnum\xcoord<#1
      \raise\ycoord \Einheit\hbox to0pt{\hskip\xcoord \Einheit\Punkt\hss}%
      \advance\xcoord by1
    \repeat
    \xcoord#3
    \advance\ycoord by1
  \REPEAT}
\def\Label#1#2(#3,#4){\unskip \xdim#3 \Einheit \ydim#4 \Einheit
  \def\lo{\advance\xdim by-.5 \Einheit \advance\ydim by.5 \Einheit}%
  \def\llo{\advance\xdim by-.25cm \advance\ydim by.5 \Einheit}%
  \def\loo{\advance\xdim by-.5 \Einheit \advance\ydim by.25cm}%
  \def\o{\advance\ydim by.25cm}%
  \def\ro{\advance\xdim by.5 \Einheit \advance\ydim by.5 \Einheit}%
  \def\rro{\advance\xdim by.25cm \advance\ydim by.5 \Einheit}%
  \def\roo{\advance\xdim by.5 \Einheit \advance\ydim by.25cm}%
  \def\l{\advance\xdim by-.30cm}%
  \def\r{\advance\xdim by.30cm}%
  \def\lu{\advance\xdim by-.5 \Einheit \advance\ydim by-.6 \Einheit}%
  \def\llu{\advance\xdim by-.25cm \advance\ydim by-.6 \Einheit}%
  \def\luu{\advance\xdim by-.5 \Einheit \advance\ydim by-.30cm}%
  \def\u{\advance\ydim by-.30cm}%
  \def\ru{\advance\xdim by.5 \Einheit \advance\ydim by-.6 \Einheit}%
  \def\rru{\advance\xdim by.25cm \advance\ydim by-.6 \Einheit}%
  \def\ruu{\advance\xdim by.5 \Einheit \advance\ydim by-.30cm}%
  #1\raise\ydim\hbox to0pt{\hskip\xdim
     \vbox to0pt{\vss\hbox to0pt{\hss$#2$\hss}\vss}\hss}%
}
\catcode`\@=13

\TagsOnRight

\def\AnCTAA{1}
\def\BuchAA{2}
\def\ChPfAA{3}
\def\DrubAA{4}
\def\DrSmAA{5}
\def\FaGuAA{6}
\def\GrKPAA{7}
\def\HawkAA{8}
\def\HoOrAA{9}
\def\HoLLAA{10}
\def\KiSYAA{11}
\def\KnMYAA{12}
\def\KratAC{13}
\def\LaLiAA{14}
\def\LaLiAB{15}
\def\PaYuAA{16}
\def\PeScAA{17}
\def\ReTYAA{18}
\def\SlatAC{19}
\def\YuTiAA{20}

\def\AAA{1.1}
\def\AAB{1.2}
\def\AAC{1.3}
\def\AAD{1.4}
\def\AAE{1.5}
\def\AAF{1.6}
\def\AAG{1.7}
\def\AAH{1.8}
\def\AAI{1.9}
\def\AAJ{1.10}
\def\AAK{3.1}
\def\AA{3.2}
\def\Aaa{3.3}
\def\AAa{3.4}
\def\Abb{3.5}
\def\Acc{3.6}
\def\Add{3.7}
\def\Aee{3.8}
\def\Aff{3.9}
\def\Agg{3.10}
\def\Ahh{4.1}
\def\Aii{4.2}
\def\AAb{5.1}
\def\AB{5.2}
\def\ABa{5.3}
\def\ACc{5.4}
\def\ACa{5.5}
\def\ACb{5.6}
\def\AC{5.7}
\def\AD{5.8}
\def\AE{5.9}
\def\AF{5.10}
\def\AFa{5.11}
\def\AG{6.1}
\def\AJ{6.2}
\def\AL{9.1}
\def\AN{9.2}
\def\AH{10.1}

\def\TA{1}
\def\TB{2}
\def\TC{3}
\def\TD{4}
\def\TE{5}
\def\TF{6}
\def\TG{7}
\def\TH{8}
\def\TI{9}
\def\TJ{10}
\def\TK{11}
\def\TL{12}
\def\TM{13}
\def\TN{14}
\def\TO{15}
\def\TP{16}
\def\TQ{17}
\def\TR{18}
\def\TS{19}
\def\TT{20}
\def\TU{21}
\def\TV{22}
\def\TW{23}
\def\TX{24}
\def\TY{25}
\def\TZ{26}
\def\TTA{27}
\def\TTB{28}
\def\TTC{29}
\def\TTD{30}
\def\TTE{31}
\def\TTF{32}
\def\TTG{33}
\def\TTH{34}
\def\TTI{35}
\def\TTJ{36}

\def\FA{1}
\def\FB{2}
\def\FC{3}

\def\coef#1{\left\langle#1\right\rangle}
\def\po#1#2{(#1)_#2}

\topmatter 
\title Refined enumeration of two-rowed set-valued standard tableaux
via two-coloured Motzkin paths
\endtitle 
\author C.~Krattenthaler$^\dagger$
\endauthor

\affil 
Fakult\"at f\"ur Mathematik, Universit\"at Wien,\\
Oskar-Morgenstern-Platz~1, A-1090 Vienna, Austria.\\
WWW: \tt http://www.mat.univie.ac.at/\~{}kratt
\endaffil
\address Fakult\"at f\"ur Mathematik, Universit\"at Wien,
Oskar-Morgenstern-Platz~1, A-1090 Vienna, Austria.\newline
http://www.mat.univie.ac.at/\~{}kratt
\endaddress

\subjclass Primary 05A15;
 Secondary 05E10 33C20
\endsubjclass
\keywords Set-valued tableau, Motzkin path,
generating function calculus, Lagrange inversion,
hypergeometric series, contiguous relation, Chu--Vandermonde summation
\endkeywords
\thanks $^\dagger$Partially supported by the Austrian
Science Foundation FWF, grant 10.55776/F1002, in the framework
of the Special Research Program ``Discrete Random Structures:
Enumeration and Scaling Limits"\endthanks
\abstract 
We derive formulae for the number of set-valued standard tableaux of
two-rowed shapes, keeping track of the total number of entries, the
number of entries in the first row, and the number of entries in
the second row. Key in the proofs is a bijection with two-coloured
Motzkin paths followed by generating function computations and
coefficient extraction helped by the Lagrange inversion formula.
\endabstract
\endtopmatter
\document

\subhead 1. Introduction and statement of results\endsubhead
The subject matter of this article is {\it set-valued standard tableaux}
of two-rowed shapes.
Given a (straight or skew) Ferrers diagram (with any number of rows),
a set-valued standard
tableau of that shape is a filling of the cells of
that diagram with $1,2,\dots,n$,
for some positive integer~$n$, each being used exactly once, such that
each cell is filled with at least one number, and such that a number
in a cell~$\rh$ is less than {\it all\/} numbers in cells that are
to the right and below~$\rh$. See Section~2 for
precise definitions.

Set-valued tableaux have been introduced by Buch~\cite{\BuchAA}
in order to give a combinatorial description of the Littlewood--Richardson
coefficients for the $K$-theory of Gra\ss mannians. 
Since then set-valued tableaux were used and studied in various
algebraic and combinatorial contexts, see
\cite{\AnCTAA, \ChPfAA, \DrubAA, \DrSmAA, \FaGuAA, \HawkAA, \HoOrAA,
\HoLLAA, \KiSYAA, \KnMYAA, \LaLiAB, \PaYuAA, \PeScAA,  \ReTYAA, \YuTiAA}
and references therein for a highly non-exhaustive list.

In the extended abstract
\cite{\LaLiAA}, Lazar and Linusson embarked on the enumeration
of two-rowed set-valued standard tableaux.
Their approach was mainly experimental;
whenever computer-aided computations seemed to hint at a ``nice" formula,
subsequently they tried to prove that formula combinatorially.
One key ingredient in their approach was a (straightforward)
bijection with two-coloured Motzkin paths. In that manner, they
derived a number of enumeration formulae, and left one open
conjecture~\cite{\LaLiAA, Conj.~16}.
Lazar and Linusson managed to prove this conjecture
in the full version~\cite{\LaLiAB, Theorem~20} of the extended abstract
by applying generating function calculus, following a suggestion
of Pierre Bonnet and Guillaume Chapuy.

The purpose of the present article is to refine and generalise
the results of Lazar and Linusson. The message is that, once
the original enumeration problem has been converted into an
enumeration problem for Motzkin paths by the earlier mentioned
bijection, it is indeed generating
function calculus that allows one to {\it compute} (without guessing)
whatever number one is interested in. In particular,
this leads one to formulae that seem very hard to guess in the first
place (see for example Theorems~\TA\ and~\TF\ below).
As we shall see, this
calculus can only be fully exploited if one puts in the
{\it Lagrange inversion formula} (see Section~4). 

In this paper, we consider the enumeration of set-valued standard tableaux
of two-rowed
straight shapes as well as of skew shapes.
Our main result for straight shapes is
a formula for the number of set-valued standard tableaux of a given
shape and with fixed numbers of entries in the first row and in
the second row.

\proclaim{Theorem \TA}
Let $n$ be a positive integer and $t,c,d,e$ non-negative integers
with $c+d+2e+t=n$.
The number of set-valued standard tableaux of shape $(e+t,e)$ with
$c+e+t$ entries in the first row and
with $d+e$ entries in the second row is given by
$$\multline
\chi(d=e=0)\binom {n-1}{t-1}
+\frac {(n-1)!} {(d+e)\cdot c!\,d!\,(e-1)!\,(e+t-1)!}\\
-\sum_{b=n-c-e+1}^{n-e+1}(-1)^{n-b-c-e-1}
\frac {(n-b)\cdot(n-1)!} {b\cdot d!\,(b-1-d)!\,
(e-1)!\,(n-b-e+1)!},
\endmultline
\tag\AAA$$
where $\chi(\Cal S)=1$ if $\Cal S$ is true and $\chi(\Cal S)=0$
otherwise. Here, and later, an expression containing
a term $m!$ with $m<0$ in the denominator has to be interpreted
as~$0$.
\endproclaim

\remark{Remark}
Clearly, the number $n$ is the total number of entries of the
tableaux in the assertion of the theorem.
\endremark

By summation over all possible $c$ and $d$, we obtain the number
of set-valued tableaux of a given shape and with a fixed number
of entries.

\proclaim{Corollary \TB}
Let $n$ be a positive  integer and $t,e$ non-negative integers.
The number of set-valued standard tableaux of shape $(e+t,e)$
with $n$~entries is given by
$$\multline
\chi(e=0)\binom {n-1}{t-1}
+\sum_{d=0}^{n-2e-t}\frac {(n-1)!}
{(n-d-2e-t)!\,d!\,(e-1)!\,(e+t-1)!}\\
\cdot
\left(\frac {1} {d+e}
-\frac {n-d-e-t-1} {(d+e+t+1)(d+e+t)}\right).
\endmultline
\tag\AAB$$
\endproclaim

We also recover Lazar and Linusson's cumulative results on the
number of set-valued standard tableaux with a given number of
entries.

\proclaim{Corollary \TC}
Let $n$ and $m$ be positive integers and\/ $t$ a non-negative integer.
The number of set-valued standard tableaux of shape
$(e+t,e)$, for some non-negative integer~$e$, 
with $m$~entries in the first row and $n$ entries in total
is given by
$$
\chi(m=n)\binom {n-1}{t-1}
+\chi(m<n)\frac {t} {n-1}\binom nm\binom {n-1}{m-t-1}
+\frac {1} {n-1}\binom {n-1}m\binom {n-1}{m-t-1}.
\tag\AAC
$$
\endproclaim

The special case $t=0$ of this corollary appears in
\cite{\LaLiAA, Prop.~8,~1.} and \cite{\LaLiAB, Cor.~2,~(2)}.

\proclaim{Corollary \TD\
(\cite{\LaLiAA, Theorem~13}, \cite{\LaLiAB, Theorem~13})}
Let $n$ be a positive integer and $t$ a non-negative integer.
The number of set-valued standard tableaux of shape
$(e+t,e)$, for some non-negative integer~$e$, and $n$ entries in total
is given by
$$
\binom {2n-2}{n-t-1}-\binom {2n-2}{n-t-2}+\binom {n-2}{t-2}.
\tag\AAD
$$
\endproclaim

The earlier mentioned conjecture of Lazar and Linusson~\cite{\LaLiAA,
Conj.~16}, subsequently proved 
in~\cite{\LaLiAB, Theorem~20}, concerned the average length of
the second row in all set-valued standard tableaux
of a 2-row rectangular shape with a given number of entries.
Using our approach, we generalise that result to
the set-valued standard tableaux featuring in Corollary~\TD.

\proclaim{Theorem \TE}
Let $n$ be a positive integer and $t$ a non-negative integer.
The expected value of the length of the second row
in set-valued standard tableaux with $n$ entries
with a {\rm(}straight{\rm)}
shape where the first row is by~$t$ longer than the second row
is given by
$$
\frac {\binom {2n-4}{n-t-1}
  + (n - 2) \binom {2 n - 4} {  n - t - 3}
  - (n + 1) \binom {2 n - 4} {n - t - 4} - 
 \binom {n - 3} {t-2}}
 {\binom {2n-2}{n-t-1}-\binom {2n-2}{n-t-2}+\binom {n-2}{t-2}}.
\tag\AAE
$$
\endproclaim

Our main result for the enumeration of set-valued standard tableaux
of two-rowed {\it skew} shapes is the following.

\proclaim{Theorem \TF}
Let $n$ and $f$ be positive integers and $t,c,d,e$ non-negative integers
with $c+d+2e-f+t=n$.

If $t<f$, the number of set-valued standard tableaux of shape
$(e+t,e)/(f,0)$ with
$c+e+t-f$ entries in the first row and
with $d+e$ entries in the second row is given by
$$\align
\chi(&c=0\text{ and }e=f-t)\binom {n-1} {f-t-1}
-\frac {(n-1)!}
{(c+e+t)\cdot c!\,d!\,(e-f-1)!\,(e+t-1)!}\\
&
+\frac {n!}
{(c+e-f+t)(d+e)\cdot c!\,d!\,(e-1)!\,(e-f+t-1)!}\\
&
+\sum_{b=n-c-e+f+1}^{n-e+f+1}(-1)^{n-b-c-e+f}
\frac {(n-b)\cdot(n-1)!} {b\cdot d!\,(b-1-d)!\,
(e-f-1)!\,(n-b-e+f+1)!}.
\\
\tag\AAF
\endalign$$
If $t\ge f$, then this number is given by
$$\allowdisplaybreaks
\align 
\chi(d=e&{}=0)\binom {n-1}{t-f-1}
+
\frac {n!}
{(c+e-f+t)(d+e)\cdot c!\,d!\,(e-1)!\,(e-f+t-1)!}
\\
\\
&-
\frac {(n-1)!}
{(c+e+t)\cdot c!\,d!\,(e-f-1)!\,(e+t-1)!}
\\
&
+\sum_{b=n-c-e+f+1}^{n-e+f+1}(-1)^{n-b-c-e+f}
\frac {(n-b)\cdot(n-1)!} {b\cdot d!\,(b-1-d)!\,
(e-f-1)!\,(n-b-e+f+1)!}.
\\
\tag\AAG
\endalign$$
\endproclaim

\remark{Remark}
Also here, the number $n$ is the total number of entries of the
tableaux in the assertion of the theorem.
\endremark

Our generating function approach also allows us to derive cumulative
formulae for two-rowed skew shapes.

\proclaim{Theorem \TG}
Let $n$ and $f$ be positive integers and $t$ a non-negative integer.

If $t<f$, the number of set-valued standard tableaux of shape
$(e+t,e)/(f,0)$, for some non-negative integer~$e$, with
a total number of $n$~entries is given by
$$\align
&\binom {n-1}{f-t-1}
+\binom {2n-2}{n-f+t-1}
-\binom {2n-2}{n-f-t-2}\\
&+
\sum_{k\ge f-t}(-1)^{k-f+t}\binom {k-1}{f-t-1}
\left(-\binom {2n+k-f+t-3}{n-k-1}
+
\binom {2n+k-f+t-3}{n-k-2}\right)
\\
\tag\AAH
\endalign$$
If $t\ge f$, then this number is given by
$$\align
&\binom {n-1}{t-f-1}
+
2\binom {2n-3}{n+f-t-2}
-
\binom {2n-2}{n-f-t-2}
\\
&+
\sum_{k\ge t-f+1}(-1)^{k-t-f-1}\binom {k-1}{t-f-1}
\left(\binom {2n+k+f-t-3}{n-k-1}-\binom {2n+k+f-t-3}{n-k-2}\right).\\
\tag\AAI
\endalign
$$
\endproclaim

\remark{Remark}
Formula (\AAI) becomes particularly simple if $t=f$ since then the
sum drops out because of the binomial coefficient
$\binom {k-1}{t-f-1}=\binom {k-1}{-1}=0$. We obtain that the
number of set-valued standard tableaux of shape
$(e+t,e)/(t,0)$, for some non-negative integer~$e$, with
a total number of $n$~entries is given by
$$
2\binom {2n-3}{n-2}
-
\binom {2n-2}{n-2t-2}.
\tag\AAJ
$$
For $t=1$ this reduces to \cite{\LaLiAA, Prop.~9,~2.}, respectively
to~\cite{\LaLiAB, Cor.~9, Eq.~(5)}.
\endremark

In the next section we give the precise definition of set-valued
standard tableaux, we introduce the Motzkin paths that we need
in the present context, and we describe the earlier mentioned
bijection between two-rowed set-valued standard tableaux and
these Motzkin paths (see Lemma~\TH).

In Section~3, we derive generating functions
for the Motzkin paths that correspond to our two-rowed set-valued
standard tableaux. As we indicated earlier, in order to extract
coefficients from these generating functions, we need the Lagrange
inversion formula. It is indeed the less frequently used (and less
known) second form of the Lagrange inversion formula that is more
convenient here, which we recall in Section~4. We apply the formula
in Sections~5--9 to derive explicit expressions for the coefficients
of the various pieces of the generating functions that we found
in Section~3. We are then ready for the proofs of the above results.
These are collected in Section~10.

\subhead 2. Set-valued standard tableaux and Motzkin paths\endsubhead
A {\it partition} is a tuple $\la=(\la_1,\la_2,\dots,\la_r)$ of
non-negative integers which are non-increasing. 
The {\it Ferrers
diagram\/} of a partition $\lambda=(\la_1,\la_2,\dots,\la_r)$
is an array of cells with $r$ left-justified rows and $\lambda_i$
cells in row $i$. Figure~\FA.a shows the Ferrers
diagram of the partition $(4,3,3,1)$.
If $\la$ and $\mu$ are two partitions with
$\la\ge\mu$ (i.e\. the Ferrers diagram of $\mu$ is contained in the
Ferrers diagram of $\la$), then the {\it skew shape\/} $\la/\mu$
consists of all cells that are contained in $\la$ but not in $\mu$.
Figure~\FA.b shows the skew Ferrers diagram
corresponding to $(4,3,3,1)/(2,1)$.

\midinsert
\vskip10pt
\vbox{
$$
\gather
\smatrix \format \sa\c\quad \s\c\quad \s\c\quad \s\c\quad \se\\
\hlinefor9\\
&\vphantom{f}&& && && &\\
\hlinefor9\\
&\vphantom{f}&& && &\\
\hlinefor7\\
&\vphantom{f}&& && &\\
\hlinefor7\\
&\vphantom{f} &\\
\hlinefor3
\endsmatrix
\hskip1cm
\smatrix \format\sa\c\s\c\s\c\s\c\se\\
\omit&\omit&\omit&\omit&\hlinefor5\\
\omit&\hbox to10pt{\hss \hss}&\omit&\hbox to10pt{\hss \hss}&&\hbox to10pt{\hss
\hss}&&\hbox to10pt{\hss $$\hss}&\\
\omit&\omit&\hlinefor7\\
\omit& && && &\\
\hlinefor7\\
& && && &\\
\hlinefor7\\
& &\\
\hlinefor3
\endsmatrix\\
 \text{\eightpoint \hskip.5cm a. Ferrers diagram\hskip1cm
b. skew Ferrers diagram}
\endgather
$$
\centerline{\eightpoint Figure \FA}
}
\vskip10pt
\endinsert

Given a (straight or skew) Ferrers diagram,
a {\it set-valued standard
tableau} of that shape is a filling of
the cells of that diagram with $1,2,\dots,n$,
for some positive integer~$n$, each being used exactly once, such that
each cell is filled with at least one number, and such that a number
in a cell~$\rh$ is less than {\it all\/} numbers in cells that are
weakly to the right and weakly below~$\rh$ ($\rh$ excluded).
Figure~\FB\ displays set-valued standard tableaux of shapes
$(4,3,3,1)$ respectively $(4,3,3,1)/(2,1)$ consisting of $n=15$
entries.

\midinsert
\vskip10pt
\vbox{
$$
\gather
\smatrix \format \sa\c \s\c \s\c \s\c \se\\
\hlinefor9\\
&1,2&&5 &&8 &&12 &\\
\hlinefor9\\
&3 &&9,10,11 && 14&\\
\hlinefor7\\
&4 &&13 &&15 &\\
\hlinefor7\\
&6,7 &\\
\hlinefor3
\endsmatrix
\hskip1cm
\smatrix \format\sa\c\s\c\s\c\s\c\se\\
\omit&\omit&\omit&\omit&\hlinefor5\\
\omit&\hbox to10pt{\hss \hss}&\omit&\hbox to10pt{\hss \hss}&&\hbox to10pt{\hss
3,8\hss}&&\hbox to10pt{\hss $12$\hss}&\\
\omit&\omit&\hlinefor7\\
\omit& &&2,9,10,11 &&14 &\\
\hlinefor7\\
&1,4,5 &&13 &&15 &\\
\hlinefor7\\
&6,7 &\\
\hlinefor3
\endsmatrix
\endgather
$$
\centerline{\eightpoint Set-valued standard tableaux}
\centerline{\eightpoint Figure \FB}
}
\vskip10pt
\endinsert

\medskip
A {\it Motzkin path} is a lattice path on
the two-dimensional integer lattice $\Z^2$ with up-steps
$(x,y)\to(x +1,y +1)$, horizontal steps $(x,y)\to(x +1,y)$, and down-steps
$(x,y)\to(x +1,y -1)$, which 
never runs below the $x$-axis.\footnote{We should point out
that we use a slightly extended version of the classical
notion of ``Motzkin path'' here.
Usually, a Motzkin path is required to start and end on the $x$-axis.
We do not make this requirement here.}
The Motzkin paths that we shall need here will have two kinds of
horizontal steps. Following~\cite{\LaLiAA, \LaLiAB}, we shall implement
this by requiring that the horizontal steps come in {\it two colours},
{\it umber} and {\it denim}. (We shall see in Lemma~\TH\ that these
two colours have a mnemonic function, namely for ``up" and ``down",
referring to first respectively second row of our shapes.)
In the sequel, whenever we speak of {\it``two-coloured" Motzkin
paths}, we mean this kind of Motzkin paths.

\medskip
As already said, the key for the enumeration of two-rowed set-valued
standard tableaux is the bijection between these tableaux and
two-coloured Motzkin paths that is the subject of the following lemma.

\proclaim{Lemma \TH}
Set-valued standard tableaux of shape $(e+t,e)/(f,0)$ with
$c+e-f+t$ entries in the first row and
with $d+e$ entries in the second row {\rm(}that is, in total
with $n=c+d+2e-f+t$ entries{\rm)} are in bijection with
two-coloured Motzkin paths from $(0,f)$ to $(n,t)$
with $c$ umber horizontal steps, with $d$ denim horizontal steps,
with $e-f+t$ up-steps,
and with $e$ down-steps {\rm(}that is, in total with $n=c+d+2e-f+t$
steps{\rm)}, where an umber horizontal step does not occur before the first
up-step and not at height~$0$,
and a denim horizontal step does not occur before the first down-step.
\endproclaim

\demo{Proof}
Given a two-rowed set-valued tableau as in the statement of the
lemma, we construct a two-coloured Motzkin path from $(0,f)$
to $(n,t)$ step by step
by letting $i=1,2,\dots, n$, in that
order, and:

\roster 
\item"--"if $i$ is an entry in the first row and the smallest number
in its cell, we add an up-step;
\item"--"if $i$ is an entry in the second row and the smallest number
in its cell, we add a down-step;
\item"--"if $i$ is an entry in the first row but not the smallest number
in its cell, we add an umber horizontal step;
\item"--"if $i$ is an entry in the second row but not the smallest number
in its cell, we add a denim horizontal step.
\endroster
Figure \FC\ shows an example of this mapping for
$f=2$, $t=1$, $c=1$ and $d=3$, $e=3$. In the figure, the colouring of
horizontal steps is indicated by the labels~$u$ (for {\it umber}) and~$d$
(for {\it denim}).
\midinsert
\vskip10pt
\vbox{
$$
\gather
\smatrix \format\sa\c\s\c\s\c\s\c\se\\
\omit&\omit&\omit&\omit&\hlinefor5\\
\omit&\hbox to10pt{\hss \hss}&\omit&\hbox to10pt{\hss \hss}&&\hbox to10pt{\hss
3,4\hss}&&\hbox to10pt{\hss $8$\hss}&\\
\hlinefor9\\
& 1&&2,5,6,7 &&9 &\\
\hlinefor7
\endsmatrix
\hskip2cm
\hbox{
$
\Gitter(10,5)(-1,0)
\Koordinatenachsen(10,5)(-1,0)
\Pfad(0,2),443111134\endPfad
\Label\ro{u}(3,1)
\Label\ro{d}(4,1)
\Label\ro{d}(5,1)
\Label\ro{d}(6,1)
\hskip5cm
$}
\\
 \text{\eightpoint \hskip.5cm a. A set-valued standard tableau\hskip2cm
b. Corresponding Motzkin path\hskip1cm}
\endgather
$$
\centerline{\eightpoint Figure \FC}
}
\vskip10pt
\endinsert
It is straightforward to check that this mapping has the claimed
properties.\quad \quad \qed
\enddemo

\subhead 3. Generating functions\endsubhead
In this section, we derive formulae for the generating functions
$$
\sum_{P\in\Cal P} z^{\#(\text{steps of }P)}
x^{\#(\text{umber horizontal steps of }P)}
y^{\#(\text{denim horizontal steps of }P)}
\al^{\#(\text{down-steps of }P)}
\tag\AAK
$$
for the sets $\Cal P$ of Motzkin paths that correspond to the
set-valued standard tableaux in Theorems~\TA\ and~\TF.
These formulae will be expressed in terms of the generating
function~(\AAK) where $\Cal P$ is the set of two-coloured Motzkin
paths starting at $(0,0)$ and ending on the $x$-axis. We denote
this generating function by~$M(z)$. (For the sake of lighter notation,
we omit to indicate the dependence on $x,y,\al$ in the notation.)
Writing $\Cal M$ for the set of two-coloured Motzkin paths starting at
$(0,0)$ and ending on the $x$-axis, $u$~for an up-step,
$d$~for a down-step, $h_u$~for an umber horizontal step,
$h_d$~for a denim horizontal step,
and $\emptyset$ for the zero-length path from $(0,0)$ to $(0,0)$,
we have the following simple combinatorial decomposition:
$$
\Cal M=\{\emptyset\}\,\sqcup\, h_u\Cal M\,\sqcup\, h_d\Cal M
\,\sqcup\,
u\Cal M d\Cal M,
$$
where $\sqcup$ denotes disjoint union.
In words, a Motzkin path in $\Cal M$ is either the zero-length path,
or an umber horizontal step followed by a Motzkin path in~$\Cal M$,
or a denim horizontal step followed by a Motzkin path in~$\Cal M$,
or an up-step followed by a Motzkin path in~$\Cal M$ elevated by
one unit, followed by the first return to the $x$-axis --- which must
occur via a down-step --- followed by a Motzkin path in~$\Cal M$.

The above decomposition translates into the functional equation
$$
M(z)=1+(x+y)zM(z)+\al z^2M^2(z)
\tag\AA$$
for the generating function $M(z)$.

For later considerations, we need to single out a subset of~$\Cal
M$. Namely, let $\Cal M_0$ denote the set of all two-coloured Motzkin paths
in~$\Cal M$ with no umber horizontal step on the $x$-axis.
The corresponding generating function, that is, the series~(\AAK)
with $\Cal P=\Cal M_0$, is denoted by~$M_0(z)$. Arguing as above,
we obtain the combinatorial decomposition
$$
\Cal M_0=\{\emptyset\}\,\sqcup\, h_d\Cal M_0
\,\sqcup\,
u\Cal M d\Cal M_0.
$$
This decomposition translates into the equation
$$
M_0(z)=1+yzM_0(z)+\al z^2M(z)M_0(z),
$$
or, equivalently,
$$
M_0(z)=\frac {1} {1-yz-\al z^2M(z)}.
\tag\Aaa
$$

We are now in the position to provide a formula for the generating
function for the Motzkin paths that correspond to the
set-valued standard tableaux in Theorem~\TA. More precisely,
via the bijection in Lemma~\TH, these are the two-coloured Motzkin paths
starting at $(0,0)$ and ending at height~$t$, where an umber
horizontal step does not occur at height~0 and a denim horizontal
step does not occur before the first down-step.

\proclaim{Proposition \TI}
Let $t$ be a non-negative integer.
The generating function {\rm(\AAK)} where $\Cal P$ is the above set of
two-coloured Motzkin paths corresponding to
set-valued standard tableaux of shape $(e+t,e)$, for some non-negative
integer~$e$, is given by
$$
\left(\frac {z} {1-xz}\right)^t
+\frac {\al \left(zM(z)\right)^{t+2}} 
{(1+yzM(z))(1+xzM(z))}
+
\frac {\al zM(z)
\left(\left(zM(z)\right)^t-\left(\frac {z} {1-xz}\right)^{t}\right)}
 {\left(y+\al zM(z)\right)
\left(1+yzM(z)\right)}.
\tag\AAa$$
\endproclaim

\demo{Proof}
We start again with a combinatorial decomposition of the set of
Motzkin paths that we are interested in.
For a convenient statement of the decomposition, we need
some notation borrowed from the combinatorics of words. Namely,
for a set $\Cal A$ of paths, we write $\Cal A^*$ for the set of all
concatenations of an arbitrary number of paths from $\Cal A$,
including the zero-length path. (It will always be clear from the
context from which point these concatenations start.)
If $\Cal A=\{s\}$ is a set consisting of a single path with only one
step~$s$, then, in abuse of notation, we shall shortly write $s^*$ instead
of $\{s\}^*$. Using this notation, a combinatorial decomposition of
the set of two-coloured Motzkin paths in the statement of the
proposition is given by
$$\multline
\{uh_u^*\}^t
\,\sqcup\,
\left(
\bigsqcup_{s\ge0}\{uh_u^*\}^{s+1}d\left(\Cal Md\right)^s\Cal M_0
\left(u\Cal M\right)^t\right)\\
\,\sqcup\,
\left(
\bigsqcup_{s\ge0}\{uh_u^*\}^{s+1}d
\bigsqcup_{h=1}^{\min\{s,t\}}
\left(\Cal Md\right)^{s-h}\Cal M
\left(u\Cal M\right)^{t-h}\right).
\endmultline
\tag\Abb$$
In words: a two-coloured Motzkin path from the statement of the
proposition 

\roster 
\item"--"may be without down-steps; in that case only up-steps and
umber horizontal steps are allowed, leading to the first term in~(\Abb);
\item"--"may start with up-steps and umber horizontal steps until
height~$s+1$ before the first down-step occurs; then $s$~times
an elevated Motzkin path in~$\Cal M$ and a down-step follow,
until one reaches height~0; then follows a Motzkin path in~$\Cal M_0$
(the reader must recall that no umber horizontal step is allowed to
occur on height~0); and finally the path continues alternatingly
with an up-step and an elevated Motzkin path in~$\Cal M$ until it
reaches its final point at height~$t$; this explains the second term
in~(\Abb);
\item"--"may start with up-steps and umber horizontal steps until
height~$s+1$ before the first down-step occurs; then $s-h$~times
an elevated Motzkin path in~$\Cal M$ and a down-step follow,
thus reaching height~$h>0$; then follows an elevated Motzkin path in~$\Cal M$;
and finally the path continues alternatingly
with an up-step and an elevated Motzkin path in~$\Cal M$ until it
reaches its final point at height~$t$; this explains the third term
in~(\Abb).
\endroster
The decomposition (\Abb) straightforwardly translates into the
expression
$$\align
&\left(\frac {z} {1-xz}\right)^t
+\sum_{s\ge0}\left(\frac {z} {1-xz}\right)^{s+1}\al z
\left(M(z)\al z\right)^{s}
M_0(z)\left(M(z)z\right)^{t}\\
&\kern2cm
+
\sum_{s=0}^{t-1}\left(\frac {z} {1-xz}\right)^{s+1}\al z
\sum_{h=1}^{s}
\left(M(z)\al z\right)^{s-h}
M(z)\left(M(z)z\right)^{t-h}\\
&\kern2cm
+
\sum_{s\ge t}\left(\frac {z} {1-xz}\right)^{s+1}\al z
\sum_{h=1}^{t}
\left(M(z)\al z\right)^{s-h}
M(z)\left(M(z)z\right)^{t-h}.
\endalign$$
for the generating function for these two-coloured Motzkin paths.

We denote the three sums on the right-hand side
by $S_1,S_2,S_3$, in order.
Then, using (\Aaa), we have
$$
S_1=
\frac {\al z^2M_0(z)\left(zM(z)\right)^{t}} 
{1-xz-\al z^2M(z)}
=\frac {\al z^2\left(zM(z)\right)^{t}} 
{(1-xz-\al z^2M(z))(1-yz-\al z^2M(z))}.
$$
Furthermore, by (\AA) we have
$$
M(z)(1-xz-\al z^2M(z))=1+yzM(z)
\tag\Acc
$$
and
$$
M(z)(1-yz-\al z^2M(z))=1+xzM(z).
\tag\Add
$$
Therefore, our last expression for $S_1$ can be rewritten in the
form
$$
S_1
=\frac {\al \left(zM(z)\right)^{t+2}} 
{(1+yzM(z))(1+xzM(z))}.
$$

Next, by summing geometric series,  we have
$$\align
S_2&=\sum_{s=0}^{t-1}\left(\frac {z} {1-xz}\right)^{s+1}\al zM(z)
\left(zM(z)\right)^{t-s}\sum_{h=1}^{s}
\left(\al zM(z)\right)^{s-h}
\left(zM(z)\right)^{s-h}\\
&=\sum_{s=0}^{t-1}\left(\frac {z} {1-xz}\right)^{s+1}\al zM(z)
\left(zM(z)\right)^{t-s}
\frac {1-\left(\al z^2M^2(z)\right)^{s}}
{1-\al z^2M^2(z)} \\
&=\sum_{s=0}^{t-1}\left(\frac {z} {1-xz}\right)^{s+1}\al zM(z)
\left(zM(z)\right)^{t-s}
\frac {1}
{1-\al z^2M^2(z)} \\
&\kern2cm
-\sum_{s=0}^{t-1}\left(\frac {z} {1-xz}\right)^{s+1}\al zM(z)
\left(zM(z)\right)^{t-s}
\frac {\left(\al z^2M^2(z)\right)^{s}}
{1-\al z^2M^2(z)} \\
&=\left(\frac {z} {1-xz}\right)^{t}
\frac {1-\left((1-xz)M(z)\right)^t}
{1-(1-xz)M(z)} 
\cdot
\frac {\al z^2M^2(z)}
{1-\al z^2M^2(z)} \\
&\kern2cm
-\left(zM(z)\right)^{t}
\frac {1-\left(\frac {\al z^2M(z)} {1-xz}\right)^{t}}
{1-xz-\al z^2M(z)} 
\cdot
\frac {\al z^2M(z)}
{1-\al z^2M^2(z)} .
\endalign$$
Now, again by the defining equation~(\AA),
$$
1-(1-xz)M(z)=-zM(z)(y+\al zM(z)).
\tag\Aee
$$
Using this relation and (\Acc),
our last expression for $S_2$ can be rewritten as
$$\align
S_2
&=
\frac {\left(zM(z)\right)^t-
\left(\frac {z} {1-xz}\right)^{t}}
{y+\al zM(z)} 
\cdot
\frac {\al zM(z)}
{1-\al z^2M^2(z)} \\
&\kern2cm
-
\frac {\left(zM(z)\right)^{t}
-\left(\al z^2M^2(z)\right)^{t}\left(\frac {z} {1-xz}\right)^{t}}
{1+yzM(z)} 
\cdot
\frac {\al z^2M^2(z)}
{1-\al z^2M^2(z)} .
\endalign$$

Similarly, again by summing geometric series, we get
$$\align
S_3&=\sum_{h=1}^{t}
\sum_{s\ge t}\left(\frac {z} {1-xz}\right)^{s+1}\al z
\left(\al zM(z)\right)^{s-h}
M(z)\left(zM(z)\right)^{t-h}\\
&=
\sum_{h=1}^{t}\left(\frac {z} {1-xz}\right)^t
\left(\al zM(z)\right)^{t-h}
\left(zM(z)\right)^{t-h}
\frac {\al z^2M(z)}
{1-xz-\al z^2M(z)}\\
&=\left(\frac {z} {1-xz}\right)^t
\frac {1-\left(\al z^2M^2(z)\right)^{t}}
{1-\al z^2M^2(z)} 
\cdot
\frac {\al z^2M(z)}
{1-xz-\al z^2M(z)}.
\endalign$$
Using (\Acc), this expression can be rewritten as
$$
S_3
=\left(\frac {z} {1-xz}\right)^t
\frac {1-\left(\al z^2M^2(z)\right)^{t}}
{1-\al z^2M^2(z)} 
\cdot
\frac {\al z^2M^2(z)}
{1+yzM(z)}.
$$

Hence, in total we get
$$\multline
\left(\frac {z} {1-xz}\right)^t+S_1+S_2+S_3
=\left(\frac {z} {1-xz}\right)^t\\
+\frac {\al \left(zM(z)\right)^{t+2}} 
{(1+yzM(z))(1+xzM(z))}
+
\frac {\al zM(z)
\left(\left(zM(z)\right)^t-\left(\frac {z} {1-xz}\right)^{t}\right)}
 {\left(y+\al zM(z)\right)
\left(1+yzM(z)\right)},
\endmultline$$
where we have used the defining equation (\AA) for $M(z)$
another time to simplify the result.
This is exactly the expression in~(\AAa).\quad \quad \qed
\enddemo

Our next proposition provides a formula for the generating
function for the Motzkin paths that correspond to the
set-valued standard tableaux of skew shape in Theorem~\TF. More precisely,
via the bijection in Lemma~\TH, these are the two-coloured Motzkin paths
starting at $(0,f)$ and ending at height~$t$, where an umber
horizontal step does not occur before an up-step and not
at height~0,
and a denim horizontal
step does not occur before the first down-step.

\proclaim{Proposition \TJ}
Let $f$ be a positive integer and $t$ a non-negative integer.

If $t<f$,
the generating function {\rm(\AAK)} where $\Cal P$ is the above set of
two-coloured Motzkin paths corresponding to
set-valued standard tableaux of shape $(e+t,e)/(f,0)$, for some non-negative
integer~$e$, is given by
$$\align
&\left(\frac {\al z} {1-yz}\right)^{f-t}
+\frac {\al^{f+1} \left(zM(z)\right)^{t+1}} {1+xzM(z)}
\frac {\left(zM(z)\right)^f
-\left(\frac {z} {1-yz}\right)^{f}} {x+\al zM(z)}\\
&\kern1cm
+
\al^{f-t+1} zM(z)
\frac {\left(zM(z)\right)^{f-t}
-\left(\frac {z} {1-yz}\right)^{f-t}} {x+\al zM(z)}
\frac {1-\left(\al z^2M^2(z)\right)^t} {1-\al z^2M^2(z)}\\
&\kern1cm
+\left(\frac {\al z} {1-yz}\right)^{f-t}
\frac {1-\left(\frac {\al z^2M(z)} {1-yz}\right)^{t}}
{1+xzM(z)}
\frac {\al \left(zM(z)\right)^2} {1-\al z^2M^2(z)}\\
&\kern2cm
-\al^{f+1}\left(\frac { z} {1-yz}\right)^{f-t}
\frac {\left( zM(z)\right)^{t}
-\left(\frac { z} {1-yz}\right)^{t}} {x+\al zM(z)}
\frac {
\left(zM(z)\right)^{t+1}} {1-\al z^2M^2(z)}\\
&\kern1cm
+\frac {\al^{f+1} \left(zM(z)\right)^{f+t+2}} 
{(1+yzM(z))(1+xzM(z))}
+
\frac {1-\left(\al z^2M^2(z)\right)^{t}}
{1-\al z^2M^2(z)} 
\cdot
\frac {\al^{f-t+1}\left( zM(z)\right)^{f-t+2}}
{1+yzM(z)}.
\tag\Aff
\endalign$$
If $t\ge f$, then this generating function equals
$$\allowdisplaybreaks
\align
&\left(\frac {z} {1-xz}\right)^{t-f}
+\frac {\al\left(zM(z)\right)^{t-f+2}
-\al^{f+1}\left(zM(z)\right)^{f+t+2}} {(1+xzM(z))(1-\al z^2M^2(z))}\\
&\kern1cm
+\frac {\al^{f+1} \left(zM(z)\right)^{f+t+2}} 
{(1+yzM(z))(1+xzM(z))}
+
\frac {\al zM(z)
\left(\left(zM(z)\right)^{t-f}-\left(\frac {z} {1-xz}\right)^{t-f}\right)}
{(y+\al zM(z))(1+yzM(z))}\\ 
&\kern2cm
+
\frac {\al\left(zM(z)\right)^{t-f+2}}
{1+yzM(z)} 
\cdot
\frac {1-\left(\al z^2M^2(z)\right)^{f}}
{1-\al z^2M^2(z)} .
\tag\Agg
\endalign$$
\endproclaim

\demo{Proof}
We proceed as in the proof of Proposition~\TI. Therefore we shall be
briefer here.

At the beginning, there stands again a combinatorial decomposition of
the generating function in question, namely
$$\multline
\left\{\matrix 
\{uh_u^*\}^{t-f}&\text{if }t\ge f\\
\{dh_d^*\}^{f-t}&\text{if }t< f\\
\endmatrix\right\}
\,\sqcup\,
\left(
\bigsqcup_{s=1}^f\{dh_d^*\}^{s}u\left(\Cal Md\right)^{f-s+1}\Cal M_0
\left(u\Cal M\right)^t\right)\\
\,\sqcup\,
\left(
\bigsqcup_{s=1}^f\{dh_d^*\}^{s}u
\bigsqcup_{h=1}^{\min\{f-s+1,t\}}
\left(\Cal Md\right)^{f-s-h+1}\Cal M
\left(u\Cal M\right)^{t-h}\right)\\
\,\sqcup\,
\left(
\bigsqcup_{s\ge0}\{uh_u^*\}^{s+1}d\left(\Cal Md\right)^{s+f}\Cal M_0
\left(u\Cal M\right)^t\right)\\
\,\sqcup\,
\left(
\bigsqcup_{s\ge0}\{uh_u^*\}^{s+1}d
\bigsqcup_{h=1}^{\min\{s+f,t\}}
\left(\Cal Md\right)^{s+f-h}\Cal M
\left(u\Cal M\right)^{t-h}\right).
\endmultline$$
This decomposition translates into the expression
$$\allowdisplaybreaks
\align
&\chi(t\ge f)\left(\frac {z} {1-xz}\right)^{t-f}
+\chi(t<f)\left(\frac {\al z} {1-yz}\right)^{f-t}\\
&\kern1cm
+\sum_{s=1}^{f}\left(\frac {\al z} {1-yz}\right)^s
z\left(M(z)\al z\right)^{f-s+1}
M_0(z)\left(M(z)z\right)^{t}\vphantom{\sum_f^t}\\
&\kern2cm
+\sum_{s=1}^{f-t}\left(\frac {\al z} {1-yz}\right)^s
z\sum_{h=1}^{t}
\left(M(z)\al z\right)^{f-s-h+1}
M(z)\left(M(z)z\right)^{t-h}\\
&\kern2cm
+\sum_{s=\max\{1,f-t+1\}}^{f}\left(\frac {\al z} {1-yz}\right)^s
z\sum_{h=1}^{f-s+1}
\left(M(z)\al z\right)^{f-s-h+1}
M(z)\left(M(z)z\right)^{t-h}
\\
&\kern1cm
+\sum_{s\ge0}\left(\frac {z} {1-xz}\right)^{s+1}\al z
\left(M(z)\al z\right)^{s+f}
M_0(z)\left(M(z)z\right)^{t}\\
&\kern2cm
+
\sum_{s=0}^{t-f-1}\left(\frac {z} {1-xz}\right)^{s+1}\al z
\sum_{h=1}^{s+f}
\left(M(z)\al z\right)^{s+f-h}
M(z)\left(M(z)z\right)^{t-h}\\
&\kern2cm
+
\sum_{s\ge \max\{0,t-f\}}\left(\frac {z} {1-xz}\right)^{s+1}\al z
\sum_{h=1}^{t}
\left(M(z)\al z\right)^{s+f-h}
M(z)\left(M(z)z\right)^{t-h}.
\endalign$$
for the generating function in the statement of the proposition.

We denote the six sums on the right-hand side
by $S_1,S_2,S_3,S_4,S_5,S_6$, in order.
Then, by summing a geometric series and using (\Aaa) and (\Add), we have
$$\align
S_1&=\left(\frac {\al z} {1-yz}\right)^{f}
\al z^2M(z)
M_0(z)\left(M(z)z\right)^{t}\vphantom{\sum_f^t}
\frac {1-\left((1-yz)M(z)\right)^f} {1-(1-yz)M(z)}\\
&=\left(\frac {\al z} {1-yz}\right)^{f}
\frac {\al z\left(zM(z)\right)^{t+1}} {1-yz-\al z^2M(z)}
\frac {1-\left((1-yz)M(z)\right)^f} {1-(1-yz)M(z)}\\
&=
\frac {\al^{f+1} \left(zM(z)\right)^{t+1}} {1+xzM(z)}
\frac {\left(zM(z)\right)^f
-\left(\frac {z} {1-yz}\right)^{f}} {x+\al zM(z)}.
\endalign$$
By summing geometric series and using an adapted form of~(\Aee),
we have
$$\allowdisplaybreaks
\align
S_2&=\sum_{s=1}^{f-t}\left(\frac {\al z} {1-yz}\right)^s
zM(z)\left(\al zM(z)\right)^{f-t-s+1}
\sum_{h=1}^{t}
\left(\al zM(z)\right)^{t-h}
\left(M(z)z\right)^{t-h}\\
&=\sum_{s=1}^{f-t}\left(\frac {\al z} {1-yz}\right)^s
zM(z)\left(\al zM(z)\right)^{f-t-s+1}
\frac {1-\left(\al z^2M^2(z)\right)^t} {1-\al z^2M^2(z)}\\
&=\chi(t<f)\left(\frac {\al z} {1-yz}\right)^{f-t}
\al z^2M^2(z)
\frac {1-\left((1-yz)M(z)\right)^{f-t}} {1-(1-yz)M(z)}
\frac {1-\left(\al z^2M^2(z)\right)^t} {1-\al z^2M^2(z)}\\
&=\chi(t<f)
\al^{f-t+1} zM(z)
\frac {\left(zM(z)\right)^{f-t}
-\left(\frac {z} {1-yz}\right)^{f-t}} {x+\al zM(z)}
\frac {1-\left(\al z^2M^2(z)\right)^t} {1-\al z^2M^2(z)}.
\endalign$$
Summing geometric series and using (\Add) and an adapted form
of~(\Aee), we get
$$\allowdisplaybreaks
\align
S_3&=\sum_{s=F}^{f}\left(\frac {\al z} {1-yz}\right)^s
zM(z)\left(\al zM(z)\right)^{f-t-s+1}
\sum_{h=1}^{f-s+1}
\left(\al zM(z)\right)^{t-h}
\left(zM(z)\right)^{t-h}\\
&=\sum_{s=F}^{f}\left(\frac {\al z} {1-yz}\right)^s
\left(zM(z)\right)^{t-f+s}
\frac {1-\left(\al z^2M^2(z)\right)^{f-s+1}} {1-\al z^2M^2(z)}\\
&=\sum_{s=F}^{f}\left(\frac {\al z} {1-yz}\right)^s
\frac {\left(zM(z)\right)^{t-f+s}} {1-\al z^2M^2(z)}\\
&\kern2cm
-\sum_{s=F}^{f}\left(\frac {\al z} {1-yz}\right)^s
\frac {\left(zM(z)\right)^{t-f+s}
\left(\al z^2M^2(z)\right)^{f-s+1}} {1-\al z^2M^2(z)}\\
&=
\left(\frac {\al z} {1-yz}\right)^F
\frac {1-\left(\frac {\al z^2M(z)} {1-yz}\right)^{f-F+1}}
{1-yz-\al z^2M(z)}
\frac {(1-yz)\left(zM(z)\right)^{t-f+F}} {1-\al z^2M^2(z)}\\
&\kern2cm
-\left(\frac {\al z} {1-yz}\right)^{f}
\frac {1-\left((1-yz)M(z)\right)^{f-F+1}} {1-(1-yz)M(z)}
\frac {\al z^2M^2(z)
\left(zM(z)\right)^{t}} {1-\al z^2M^2(z)}\\
&=
\left(\frac {\al z} {1-yz}\right)^F
\frac {1-\left(\frac {\al z^2M(z)} {1-yz}\right)^{f-F+1}}
{1+xzM(z)}
\frac {(1-yz)M(z)\left(zM(z)\right)^{t-f+F}} {1-\al z^2M^2(z)}\\
&\kern2cm
-\al^{f+1}\left(\frac { z} {1-yz}\right)^{F-1}
\frac {\left( zM(z)\right)^{f-F+1}
-\left(\frac { z} {1-yz}\right)^{f-F+1}} {x+\al zM(z)}
\frac {
\left(zM(z)\right)^{t+1}} {1-\al z^2M^2(z)},
\endalign$$
where $F=\max\{1,f-t+1\}$.

Summing a geometric series and using (\Aaa), (\Acc) and~(\Add), we have
$$\align
S_4&=
\frac {\al^{f+1} z^2M_0(z)\left(zM(z)\right)^{f+t}} 
{1-xz-\al z^2M(z)}\\
&=\frac {\al^{f+1} z^2\left(zM(z)\right)^{f+t}} 
{(1-xz-\al z^2M(z))(1-yz-\al z^2M(z))}\\
&=\frac {\al^{f+1} \left(zM(z)\right)^{f+t+2}} 
{(1+yzM(z))(1+xzM(z))}.
\endalign$$
Furthermore, if $t<f$ the sum $S_5$ is empty and hence equal to zero.
Otherwise, by summing geometric series and using~(\Acc) and~(\Aee), we have
$$\allowdisplaybreaks
\align
S_5&=\sum_{s=0}^{t-f-1}\left(\frac {z} {1-xz}\right)^{s+1}\al zM(z)
\left(zM(z)\right)^{t-f-s}\sum_{h=1}^{s+f}
\left(\al zM(z)\right)^{s+f-h}
\left(zM(z)\right)^{s+f-h}\\
&=\sum_{s=0}^{t-f-1}\left(\frac {z} {1-xz}\right)^{s+1}\al zM(z)
\left(zM(z)\right)^{t-f-s}
\frac {1-\left(\al z^2M^2(z)\right)^{s+f}}
{1-\al z^2M^2(z)} \\
&=\sum_{s=0}^{t-f-1}\left(\frac {z} {1-xz}\right)^{s+1}\al zM(z)
\left(zM(z)\right)^{t-f-s}
\frac {1}
{1-\al z^2M^2(z)} \\
&\kern2cm
-\sum_{s=0}^{t-f-1}\left(\frac {z} {1-xz}\right)^{s+1}\al zM(z)
\left(zM(z)\right)^{t-f-s}
\frac {\left(\al z^2M^2(z)\right)^{s+f}}
{1-\al z^2M^2(z)} \\
&=\chi(t\ge f)\left(\frac {z} {1-xz}\right)^{t-f}
\frac {1-\left((1-xz)M(z)\right)^{t-f}}
{1-(1-xz)M(z)} 
\cdot
\frac {\al z^2M^2(z)}
{1-\al z^2M^2(z)} \\
&\kern2cm
-\chi(t\ge f)\left(zM(z)\right)^{t-f}
\frac {1-\left(\frac {\al z^2M(z)} {1-xz}\right)^{t-f}}
{1-xz-\al z^2M(z)} 
\cdot
\frac {\al z^2M(z)\left(\al z^2M^2(z)\right)^f}
{1-\al z^2M^2(z)} \\
&=\chi(t\ge f)
\frac {\left(zM(z)\right)^{t-f}-\left(\frac {z} {1-xz}\right)^{t-f}}
{y+\al zM(z)} 
\cdot
\frac {\al zM(z)}
{1-\al z^2M^2(z)} \\
&\kern2cm
-\chi(t\ge f)\al^{f+1}
\frac {1-\left(\frac {\al z^2M(z)} {1-xz}\right)^{t-f}}
{1+yzM(z)} 
\cdot
\frac {\left( zM(z)\right)^{t+f+2}}
{1-\al z^2M^2(z)} .
\endalign$$
Finally, by summing geometric series and using~(\Acc), we have
$$\align
S_6&=\sum_{h=1}^{t}
\sum_{s\ge\max\{0, t-f\}}\left(\frac {z} {1-xz}\right)^{s+1}\al z
\left(\al zM(z)\right)^{s+f-h}
M(z)\left(zM(z)\right)^{t-h}\\
&=
\sum_{h=1}^{t}\left(\frac {z} {1-xz}\right)^T
\left(\al zM(z)\right)^{T+f-h}
\left(zM(z)\right)^{t-h}
\frac {\al z^2M(z)}
{1-xz-\al z^2M(z)}\\
&=\left(\frac {z} {1-xz}\right)^T
\frac {1-\left(\al z^2M^2(z)\right)^{t}}
{1-\al z^2M^2(z)} 
\cdot
\frac {z\left(\al zM(z)\right)^{T-(t-f)+1}}
{1-xz-\al z^2M(z)}\\
&=\left(\frac {z} {1-xz}\right)^T
\frac {1-\left(\al z^2M^2(z)\right)^{t}}
{1-\al z^2M^2(z)} 
\cdot
\frac {zM(z)\left(\al zM(z)\right)^{T-(t-f)+1}}
{1+yzM(z)},
\endalign$$
where $T=\max\{0, t-f\}$.

In total, this gives
$$\allowdisplaybreaks
\align
&\chi(t\ge f)\left(\frac {z} {1-xz}\right)^{t-f}
+\chi(t<f)\left(\frac {\al z} {1-yz}\right)^{f-t}\\
&\kern1cm
+\frac {\al^{f+1} \left(zM(z)\right)^{t+1}} {1+xzM(z)}
\frac {\left(zM(z)\right)^f
-\left(\frac {z} {1-yz}\right)^{f}} {x+\al zM(z)}\\
&\kern1cm
+\chi(t<f)
\al^{f-t+1} zM(z)
\frac {\left(zM(z)\right)^{f-t}
-\left(\frac {z} {1-yz}\right)^{f-t}} {x+\al zM(z)}
\frac {1-\left(\al z^2M^2(z)\right)^t} {1-\al z^2M^2(z)}\\
&\kern1cm
+\left(\frac {\al z} {1-yz}\right)^F
\frac {1-\left(\frac {\al z^2M(z)} {1-yz}\right)^{f-F+1}}
{1+xzM(z)}
\frac {(1-yz)M(z)\left(zM(z)\right)^{t-f+F}} {1-\al z^2M^2(z)}\\
&\kern2cm
-\al^{f+1}\left(\frac { z} {1-yz}\right)^{F-1}
\frac {\left( zM(z)\right)^{f-F+1}
-\left(\frac { z} {1-yz}\right)^{f-F+1}} {x+\al zM(z)}
\frac {
\left(zM(z)\right)^{t+1}} {1-\al z^2M^2(z)}\\
&\kern1cm
+\frac {\al^{f+1} \left(zM(z)\right)^{f+t+2}} 
{(1+yzM(z))(1+xzM(z))}\\
&\kern1cm
+
\chi(t\ge f)
\frac {\left(zM(z)\right)^{t-f}-\left(\frac {z} {1-xz}\right)^{t-f}}
{y+\al zM(z)} 
\cdot
\frac {\al zM(z)}
{1-\al z^2M^2(z)} \\
&\kern2cm
-\chi(t\ge f)\al^{f+1}
\frac {1-\left(\frac {\al z^2M(z)} {1-xz}\right)^{t-f}}
{1+yzM(z)} 
\cdot
\frac {\left( zM(z)\right)^{t+f+2}}
{1-\al z^2M^2(z)} \\
&\kern1cm
+\left(\frac {z} {1-xz}\right)^T
\frac {1-\left(\al z^2M^2(z)\right)^{t}}
{1-\al z^2M^2(z)} 
\cdot
\frac {zM(z)\left(\al zM(z)\right)^{T-(t-f)+1}}
{1+yzM(z)}.
\endalign$$
For $t<f$ this simplifies to the expression in (\Aff),
while for $t\ge f$ this becomes
$$\align
&\left(\frac {z} {1-xz}\right)^{t-f}
+\frac {\al^{f+1} \left(zM(z)\right)^{t+1}} {1+xzM(z)}
\frac {\left(zM(z)\right)^f
-\left(\frac {z} {1-yz}\right)^{f}} {x+\al zM(z)}\\
&\kern1cm
+
\frac {1-\left(\frac {\al z^2M(z)} {1-yz}\right)^{f}}
{1+xzM(z)}
\frac {\al \left(zM(z)\right)^{t-f+2}} {1-\al z^2M^2(z)}
-\al^{f+1}
\frac {\left( zM(z)\right)^{f}
-\left(\frac { z} {1-yz}\right)^{f}} {x+\al zM(z)}
\frac {
\left(zM(z)\right)^{t+1}} {1-\al z^2M^2(z)}\\
&\kern1cm
+\frac {\al^{f+1} \left(zM(z)\right)^{f+t+2}} 
{(1+yzM(z))(1+xzM(z))}
+
\frac {\left(zM(z)\right)^{t-f}-\left(\frac {z} {1-xz}\right)^{t-f}}
{y+\al zM(z)} 
\cdot
\frac {\al zM(z)}
{1-\al z^2M^2(z)} \\
&\kern1cm
-\al^{f+1}
\frac {1-\left(\frac {\al z^2M(z)} {1-xz}\right)^{t-f}}
{1+yzM(z)} 
\cdot
\frac {\left( zM(z)\right)^{t+f+2}}
{1-\al z^2M^2(z)} \\
&\kern1cm
+\left(\frac {z} {1-xz}\right)^{t-f}
\frac {1-\left(\al z^2M^2(z)\right)^{t}}
{1-\al z^2M^2(z)} 
\cdot
\frac {\al z^2M^2(z)}
{1+yzM(z)}.
\endalign$$
The last expression can be mildly simplified, and one
finally arrives at~(\Agg).\quad \quad \qed
\enddemo

\subhead 4. Lagrange inversion\endsubhead
In this article, we make use of the second form of the
Lagrange inversion formula (cf.~\cite{\KratAC, Eq.~(1.2)}).

\proclaim{Theorem \TK\ (\smc Lagrange inversion)}
Let $f(z)$ be a formal power series with $f(0)=0$ and $f'(0)\ne0$,
and let $F(z)$ be its compositional inverse.
If $g(z)$ is a Laurent series with only finitely many terms
containing negative powers of~$z$, then
$$
\coef{z^n}g(F(z))=\coef{z^{-1}}g(z)f^{-n-1}(z)f'(z).
$$
\endproclaim

In the next section, we shall apply
this Lagrange inversion formula
with $F(z)=zM(z)$. Since, by~(\AA),
we have
$$
\frac {zM(z)} {1+(x+y)zM(z)+\al z^2M^2(z)}=z,
\tag\Ahh
$$
the series $F(z)=zM(z)$ is the compositional inverse of
$f(z)=\frac {z} {1+(x+y)z+\al z^2}$.
Furthermore, we have
$$
f'(z)=\frac {1-\al z^2} {\left(1+(x+y)z+\al z^2\right)^2}.
\tag\Aii
$$

\subhead 5. Coefficient extraction, I\endsubhead
In this section, we provide the basis for
the proof of Theorem~\TA\ given in Section~10.
We apply the Lagrange inversion formula
in Theorem~\TK\ in order to extract the coefficient of~$z^n$
from the various subexpressions in~(\AAa).

\proclaim{Lemma \TL}
We have
$$
\coef{z^n}\left(\frac {z} {1-xz}\right)^t=\binom {n-1}{t-1}x^{n-t}.
$$
\endproclaim

\demo{Proof}
This is a trivial consequence of the expansion of the binomial
series.\quad \quad \qed
\enddemo

\proclaim{Lemma \TM}
We have
$$\multline
\coef{z^n}\frac {\al \left(zM(z)\right)^{t+2}} {(1+xzM(z))(1+yzM(z))}
\\=
\underset c+d+2e+t=n\to{\sum_{c,d,e\ge0}}
x^cy^d\al^{e}
\left(\frac {(n-1)!} {(n-c-e)\,c!\,d!\,(e-1)!\,(n-c-d-e-1)!}
\vphantom{\sum_{b=d+1}^{n-e-1}}\right.\\
\left.
-\sum_{b=n-c-e+1}^{n-e+1}(-1)^{n-b-c-e+1}
\frac {(n-b)\cdot(n-1)!} {b\cdot d!\,(b-1-d)!\,
(e-1)!\,(n-b-e+1)!}\right).
\endmultline$$
\endproclaim

\demo{Proof}
With (\Ahh) and (\Aii) in mind, by the Lagrange inversion formula
in Theorem~\TK\ with $F(z)=zM(z)$,
the coefficient of $z^n$ in the statement of the lemma equals
$$\allowdisplaybreaks
\align
\coef{z^n}&\frac {\al \left(zM(z)\right)^{t+2}} {(1+xzM(z))(1+yzM(z))}
=\coef{z^{-1}}\frac {\al z^{t+2}} {(1+xz)(1+yz)}\cdot
\frac {\left(1+(x+y)z+\al z^2\right)^{n+1}} {z^{n+1}}\\
&\kern7cm
\cdot\frac {1-\al z^2} {\left(1+(x+y)z+\al z^2\right)^2}\\
&=\al\coef{z^{n-t-2}}\frac {\left(1+(x+y)z+\al z^2\right)^{n-1}(1-\al z^2)}
{(1+xz)(1+yz)}
\tag\AAb
\\
&=\al\coef{z^{n-t-2}}
\frac {\left((1+xz)+(1+yz)-(1-\al z^2)\right)^{n-1}(1-\al z^2)}
{(1+xz)(1+yz)}\\
&=\al\coef{z^{n-t-2}}
\underset a+b\le n-1\to{\sum_{a,b\ge0}}\binom {n-1}{a,b,n-1-a-b}
(1+xz)^{a-1} (1+yz)^{b-1}\\
&\kern4cm
\cdot(-1)^{n-1-a-b}(1-\al z^2)^{n-a-b}\\
&=\al\coef{z^{n-t-2}}
\underset a+b\le n-1\to{\sum_{a,b,c,d,e\ge0}}\binom {n-1}{a,b,n-1-a-b}
\binom {a-1}cx^cz^c
\binom {b-1}dy^dz^d\\
&\kern4cm
\cdot
(-1)^{n-1-a-b}\binom {n-a-b}e(-1)^e\al^ez^{2e},
\tag\AB
\endalign
$$
where, here and later, $\binom N{a,b,N-a-b}$ denotes the {\it  trinomial
coefficient} $\frac {N!} {a!\,b!\,(N-a-b)!}$.
As the attentive reader will have noticed, in the above calculation
we used the expansion of the binomial series several times.

At this point, we must distinguish four cases depending on whether
$a$ and~$b$ are zero or not.

We begin with the case where both $a$ and $b$ are greater than zero.
For those $a$ and~$b$, the coefficient in~(\AB) can be written in
the form
$$\multline
\underset c+d+2(e+1)+t=n\to{\sum_{c,d,e\ge0}}
x^cy^d\al^{e+1}
\\
\cdot
\underset a+b\le n-1\to{\sum_{a,b\ge1}}(-1)^{n-1-a-b-e}
\frac {(n-a-b)(n-1)!} {ab\,c!\,(a-1-c)!\,d!\,(b-1-d)!\,
e!\,(n-a-b-e)!}.
\endmultline
\tag\ABa
$$
We concentrate first on the sum over~$a$. Neglecting terms that are
independent of~$a$, this sum is
$$
\sum_{a=c+1}^{n-b-e}(-1)^a
\frac {(n-a-b)} {a\,(a-1-c)!\,(n-a-b-e)!}.
\tag\ACc
$$
We note that for this sum to be non-empty we must have $c+1\le n-b-e$,
or, equivalently,
$$
b+c+e\le n-1.
\tag\ACa
$$
In standard hypergeometric notation
$$
{}_r F_s\!\left[\matrix a_1,\dots,a_r\\ b_1,\dots,b_s\endmatrix;
z\right]=\sum _{l =0} ^{\infty}\frac {\po{a_1}{l }\cdots\po{a_r}{l }} 
{l !\,\po{b_1}{l }\cdots\po{b_s}{l }} z^l \ ,
$$
where the {\it Pochhammer symbol\/} $(\al)_m$ is defined by
$(\alpha)_m  =  \alpha(\hbox{$\alpha+1$})(\alpha+2)\cdots (\alpha+m-1)$ for
$m>0$, and $(\alpha)_0  = 1$, the above sum reads
$$
(-1)^{c+1}
\frac {(n-c-b-1)} {(c+1)\,(n-c-b-e-1)!}
{}_3 F_2\!\left[\matrix -n + b + c + 2, c+1, -n + b + c + e + 1\\
-n + b + c + 1, c+2\endmatrix;
1\right].
$$
We observe that, in the case of equality in~(\ACa), the $_3F_2$-series
contains an upper parameter~0, which means that it collapses to just~1,
so that the above expression reduces to
$$
(-1)^{c+1}
\frac {n-c-b-1} {c+1}.
\tag\ACb
$$
From now on, we assume that strict inequality holds in~(\ACa).
We apply the contiguous relation
$$
{} _{3} F _{2} \!\left [ \matrix A, B, C\\D, E
\endmatrix ; {\displaystyle z}\right ]
={} _{3} F _{2} \!\left [ \matrix  A-1, B, C\\D, E\endmatrix ;
  {\displaystyle z}\right ]   
+ z\frac {BC} {DE}
{} _{3} F _{2} \!\left [ \matrix A, B+1, C+1\\D+1, E+1
  \endmatrix ; {\displaystyle z}\right ] 
\tag\AC$$
with $A=-n + b + c + 2$ and $D=-n + b + c + 1$ to the last $_3F_2$-series.
Thereby we obtain
$$\multline
(-1)^{c+1}
\frac {(n-c-b-1)} {(c+1)\,(n-c-b-e-1)!}
\left(
{}_2 F_1\!\left[\matrix c+1, -n + b + c + e + 1\\
c+2\endmatrix;
1\right]\right.\\
\left.
+\frac {(c+1)  (n - b - c - e - 1)} {(c+2)  (n - b - c - 1)}
{}_2 F_1\!\left[\matrix c+2, -n + b + c + e + 2\\
c+3\endmatrix;
1\right]
\right).
\endmultline$$
Both $_2F_1$-series can be evaluated by means of the Chu--Vandermonde
summation (see \cite{\SlatAC, (1.7.7); Appendix
(III.4)}),
$$
{} _{2} F _{1} \!\left [ \matrix { A, -N}\\ { C}\endmatrix ; {\displaystyle
   1}\right ]  = {{({ \textstyle C-A}) _{N} }\over 
    {({ \textstyle C}) _{N} }},
\tag\AD
$$
where $N$ is a nonnegative integer. As a consequence,
after some simplification, we obtain
$$
(-1)^{c+1}\frac {n-b} {(c+1)_{n-b-c-e}}.
\tag\AE
$$
for the sum in (\ACc), however only in the case of strict inequality
in~(\ACa). Indeed, the expression in (\AE) differs from the one
for the case of equality in~(\ACa), namely~(\ACb).
If we substitute these findings in~(\ABa)
then we get
$$\multline
\underset c+d+2(e+1)+t=n\to{\sum_{c,d,e\ge0}}
x^cy^d\al^{e+1}
\left(\sum_{b=d+1}^{n-c-e-1}(-1)^{n-b-c-e}
\frac {(n-1)!} {b\,d!\,(b-1-d)!\,
e!}\frac {(n-b)} {(n-b-e)!}\right.\\
\left.
+
\frac {(n-1)!} {(n-c-e-1)\,c!\,d!\,e!\,(n-c-d-e-2)!}
\vphantom{\sum_{b=d+1}^{n-e}}\right).
\endmultline$$
Here, the second term between the big parentheses corresponds to
the ``correction" that is necessary when one replaces (\ACb)
by~(\AE) (with equality in~(\ACa), that is, with $n-b-c-e=1$).
It should be noted that the sum over~$b$ terminates at $n-c-e-1$
which is ``earlier" than the ``natural" upper bound $n-e$ (resulting
from the term $(n-b-e)!$ in the denominator of the summand).
Consequently, this sum can not be evaluated in closed form.
(Again using the contiguous relation~(\AC) and
the Chu--Vandermonde summation~(\AD),
it can be shown that, when $b$ would range over the full interval
$[d+1,n-e]$, then the sum would equal $(-1)^{e+1}\binom ne$; see
the following computation.)

Next we treat the expression (\AB) with $a=0$ and $b\ge1$, that is,
$$
\underset c+d+2(e+1)+t=n\to{\sum_{c,d,e\ge0}}
x^cy^d\al^{e+1}
\sum_{b=d+1}^{n-e}(-1)^{n-1-b-c-e}
\frac {(n-b)\cdot(n-1)!}
{b\cdot d!\,(b-1-d)!\,e!\,(n-b-e)!}
$$
We write the sum over $b$ in hypergeometric notation to get
$$
(-1)^{n-d-c-e}
\frac {(n-d-1)\,(n-1)!}
{(d+1)!\,e!\,(n-d-e-1)!}
{}_3 F_2\!\left[\matrix -n + d + 2, d+1, -n + d + e + 1\\
-n + d + 1, d+2\endmatrix; 1\right].
$$
From here on, we proceed in the same way as before, i.e., we apply
the contiguous relation~(\AC), thereby obtaining a (weighted) sum
of two $_2F_1$-series, both of which can be evaluated by means
of the Chu--Vandermonde summation~(\AD). After simplification, we
obtain
$$
(-1)^{e}\binom ne,
$$
so that the expression (\AB) with $a=0$ and $b\ge1$ reduces to
$$
\underset c+d+2(e+1)+t=n\to{\sum_{c,d,e\ge0}}
x^cy^d\al^{e+1}(-1)^{e}\binom ne.
\tag\AF$$

For reasons of symmetry the same happens for the
expression (\AB) with $a\ge1$ and $b=0$.

Finally, expression (\AB) with $a=b=0$ equals by its definition
the negative of~(\AF).

After replacement of $e$ by $e-1$,
all four cases together yield the formula
in the statement of the lemma.\quad \quad \qed
\enddemo

\proclaim{Lemma \TN}
We have
$$\multline
\coef{z^n}\frac {\al zM(z)
\left(\left(zM(z)\right)^t-\left(\frac {z} {1-xz}\right)^{t}\right)}
 {\left(y+\al zM(z)\right)
\left(1+yzM(z)\right)}\\
=\underset c+d+2e+t=n\to{\sum_{c,d,e\ge0}}
x^cy^d\al^{e}\frac {t\,(n-1)!} {(d+e)(n-c-e)\,(n-c-d-e-1)!\,c!\,d!\,(e-1)!}.
\endmultline$$
\endproclaim

\demo{Proof}
Again we apply  the Lagrange inversion formula
in Theorem~\TK\ with $F(z)=zM(z)$.
With (\Ahh) and (\Aii) in mind (in particular, the reader should
note that, to pass from $g(F(z))=g(zM(z))$ to $g(z)$, one must replace
$z$ by $f(z)=\frac {z} {1+(x+y)z+\al z^2}$, the compositional inverse
of $F(z)=zM(z)$, in the former expression), we obtain
$$\allowdisplaybreaks
\align
\coef{z^n}&\frac {\al zM(z)
\left(\left(zM(z)\right)^t-\left(\frac {z} {1-xz}\right)^{t}\right)}
 {\left(y+\al zM(z)\right)
\left(1+yzM(z)\right)}\\
&=\coef{z^{-1}}
\frac {\al z\left(z^t-\left(\frac {z} {1+yz+\al z^2}\right)^t\right)}
{(y+\al z)(1+yz)}\cdot
\frac {\left(1+(x+y)z+\al z^2\right)^{n+1}} {z^{n+1}}\\
&\kern7cm
\cdot\frac {1-\al z^2} {\left(1+(x+y)z+\al z^2\right)^2}\\
&=\coef{z^{n-1}}
\al \frac{\displaystyle
\left(z^t-\left(\frac {z} {1+yz+\al z^2}\right)^t\right)}
{y+\al z}\cdot
\frac {\left(1+(x+y)z+\al z^2\right)^{n-1}\left(1-\al z^2\right)}
{1+yz}.
\tag\AFa
\endalign$$
Partial fraction decomposition yields
$$
\frac {1} {(1+yz)(y+\al z)}
=\frac {1} {1-\al z^2}\left(\frac {1} {y+\al z}-\frac {z} {1+yz}\right).
$$
Substituting this in the last expression, we obtain
$$\allowdisplaybreaks
\align
\coef{z^{n-1}}&
\al 
\left(z^t-\left(\frac {z} {1+yz+\al z^2}\right)^t\right)
\left(1+(x+y)z+\al z^2\right)^{n-1}
\left(\frac {1} {y+\al z}-\frac {z} {1+yz}\right)\\
&=
\coef{z^{n-1}}
\al z^t
\left(1-\left(1+yz+\al z^2\right)^{-t}\right)\\
&\kern2cm
\cdot
\sum_{c\ge0}\binom {n-1}c x^cz^c
\left(1+yz+\al z^2\right)^{n-1-c}
\left(\frac {1} {y+\al z}-\frac {z} {1+yz}\right)\\
&=
\coef{z^{n-1}}
\al \sum_{c\ge0}\binom {n-1}c x^cz^{c+t}\\
&\kern1cm
\cdot
\left(
\sum_{k\ge0}\left(\binom {n-1-c}k-\binom {n-1-c-t}k\right) z^k(y+\al z)^{k-1}
\right.\\
&\kern3cm
\left.
-\sum_{e\ge0}\binom {n-1-c}e
\al^ez^{2e+1}(1+y z)^{n-c-e-2}
\right.\\
&\kern3cm
\left.
+\sum_{e\ge0}\binom {n-1-c-t}e
\al^ez^{2e+1}(1+y z)^{n-c-t-e-2}
\right)\\
&=
\coef{z^{n-1}}
\al \sum_{c\ge0}\binom {n-1}c x^cz^{c+t}\\
&\kern1cm
\cdot
\left(
\sum_{k\ge0}\left(\binom {n-1-c}k-\binom {n-1-c-t}k\right) z^k
\sum_{e\ge0}\binom {k-1}e y^{k-1-e}\al^ez^e
\right.\\
&\kern3cm
\left.
-\sum_{e\ge0}\binom {n-1-c}e
\al^ez^{2e+1}\sum_{d\ge0}\binom {n-c-e-2}d y^dz^d
\right.\\
&\kern3cm
\left.
+\sum_{e\ge0}\binom {n-1-c-t}e
\al^ez^{2e+1}\sum_{d\ge0}\binom {n-c-t-e-2}d y^dz^d
\right).
\endalign$$
Thus, we get
$$\allowdisplaybreaks
\align
&\underset c+d+2(e+1)+t=n\to{\sum_{c,d,e\ge0}}
x^cy^d\al^{e+1}
\binom {n-1}c \\
&\kern1cm
\cdot
\left(
\binom {n-1-c}{d+e+1}\binom {d+e}e 
-\binom {n-1-c-t}{d+e+1}\binom {d+e}e 
\right.\\
&\kern2cm
\left.
-\binom {n-1-c}e\binom {n-c-e-2}d 
+\binom {n-1-c-t}e\binom {n-c-t-e-2}d
\right)\\
&=\underset c+d+2(e+1)+t=n\to{\sum_{c,d,e\ge0}}
x^cy^d\al^{e+1}
\binom {n-1}c \\
&\kern1cm
\cdot
\left(
\frac {(n-c-1)!} {(n-c-d-e-2)!\,d!\,e!}
\left(\frac {1} {d+e+1}-\frac {1} {n-c-e-1}\right)
\right.\\
&\kern1.5cm
\left.
-\frac {(n-c-t-1)!} {(n-c-d-e-t-2)!\,d!\,e!}
\left(\frac {1} {d+e+1}-\frac {1} {n-c-e-t-1}\right)
\right)\\
&=\underset c+d+2(e+1)+t=n\to{\sum_{c,d,e\ge0}}
x^cy^d\al^{e+1}\frac {t\,(n-1)!} {(d+e+1)(n-c-e-1)\,(n-c-d-e-2)!\,c!\,d!\,e!}.
\endalign$$
After replacement of $e$ by $e-1$, this becomes the expression
in the statement of the lemma.\quad \quad \qed
\enddemo

\subhead 6. Coefficient extraction, II\endsubhead
In this section, we provide the basis for the proof of Corollary~\TD\
given in Section~10. That corollary concerns the cumulative number
of the set-valued standard tableaux in Theorem~\TA, where we do not
care about the exact shape or how many entries there are in the first
row and in the second row. In other words, what we want to compute
is the sum of~(\AAA) over all $c,d,e$ with $c+d+2e+t=n$. It turns
out that it is more convenient to extract the desired number
--- so-to-speak ``from scratch" --- directly from the (specialised)
generating function. To recall, this generating function was given
in Proposition~\TI. So we have to consider~(\AAa) with $x=y=\al=1$,
that is,
$$
\left(\frac {z} {1-z}\right)^t
+\frac {\left(zM(z)\right)^{t+2}} 
{(1+zM(z))^2}
+
\frac {zM(z)
\left(\left(zM(z)\right)^t-\left(\frac {z} {1-z}\right)^{t}\right)}
 {\left(1+zM(z)\right)^2
 }.
\tag\AG
$$
It is important to remember that the above specialisation has also
an effect on the series $M(z)$. More precisely,
in this special case, the series $zM(z)$ is the
compositional inverse
of $f(z)=\frac {z} {(1+z)^2}$;
the derivative of the latter series is
$\frac {1-z} {(1+z)^3}$, cf.~(\Ahh) and~(\Aii).

\proclaim{Lemma \TO}
With $x=y=\al=1$, we have
$$
\coef{z^n}\frac {\left(zM(z)\right)^{t+2}} 
{(1+zM(z))^2}
=\binom {2n-3}{n-t-2}-\binom {2n-3}{n-t-3}.
$$
\endproclaim

\demo{Proof}
By the Lagrange inversion formula in Theorem~\TK\ with $F(z)=zM(z)$,
and the above considerations concerning the implied $f(z)$, we have
$$\align
\coef{z^n}\frac {\left(zM(z)\right)^{t+2}} 
{(1+zM(z))^2}
&=\coef{z^{-1}}\frac {z^{t+2}} {(1+z)^2}\frac {(1+z)^{2n+2}} {z^{n+1}}
\frac {1-z} {(1+z)^3}\\
&=\coef{z^{n-t-2}}(1-z)(1+z)^{2n-3} \\
&=\binom {2n-3}{n-t-2}-\binom {2n-3}{n-t-3},
\endalign
$$
as desired.\quad \quad \qed
\enddemo

\proclaim{Lemma \TP}
With $x=y=\al=1$, we have
$$
\coef{z^n}\frac {zM(z)
\left(\left(zM(z)\right)^t-\left(\frac {z} {1-z}\right)^{t}\right)}
 {\left(1+zM(z)\right)^2}
=\binom {2n-3}{n-t-1}-\binom {2n-3}{n-t-2}
-\binom {n-2}{t-1}.
$$
\endproclaim

\demo{Proof}
By the Lagrange inversion formula in Theorem~\TK\ with $F(z)=zM(z)$,
and the earlier considerations concerning the implied $f(z)$, we have
$$\align
\coef{z^n}&\frac {zM(z)
\left(\left(zM(z)\right)^t-\left(\frac {z} {1-z}\right)^{t}\right)}
 {\left(1+zM(z)\right)^2}
=\coef{z^{-1}}\frac {z
\left(z^t-\left(\frac {z} {1+z+z^2}\right)^{t}\right)}
 {(1+z)^2}\frac {(1+z)^{2n+2}} {z^{n+1}}
\frac {1-z} {(1+z)^3}\\
&=\coef{z^{n-1-t}}
(1+z)^{2n-3}(1-z)
-\coef{z^{-1}}\left(\frac {z} {1+z+z^2}\right)^t
\frac {(1+z)^{2n-3}(1-z)} {z^n}\\
&=\binom {2n-3}{n-t-1}-\binom {2n-3}{n-t-2}
-\coef{z^{-1}}\left(\frac {z} {1+z+z^2}\right)^t
\frac {(1+z)^{2n-3}(1-z)} {z^n}.
\tag\AJ
\endalign
$$
We rewrite the last term as the complex contour integral
$$\multline
\frac {1} {2\pi i}\int_{\Cal C}\left(\frac {z} {1+z+z^2}\right)^t
\frac {(1+z)^{2n-3}(1-z)} {z^n}\,dz
\\=
\frac {1} {2\pi i}\int_{\Cal C}\left(\frac {z} {1+z+z^2}\right)^t
\frac {(1+z)^{2n-4}(1-z^2)} {z^n}\,dz,
\endmultline$$
where $\Cal C$ is a contour close to the origin, encircling the
origin once in positive direction. Subsequently,
we perform the substitution $v=\frac {z} {1+z+z^2}$.
For this substitution, we have $dv=dz\frac {1-z^2} {(1+z+z^2)^2}$, so that
the above integral becomes
$$
\frac {1} {2\pi i}\int_{\Cal C'}v^{t-2}\left(1+\frac {1} {v}\right)^{n-2}
\,dv
=\coef{v^{n-t-1}}(1+v)^{n-2}=\binom {n-2}{n-t-1}.
$$
If this is substituted in (\AJ), we obtain the claimed formula.\quad \quad \qed
\enddemo

\subhead 7. Coefficient extraction, III\endsubhead
In this section,
we concern ourselves with the coefficient extraction results
relevant for the proof of Theorem~\TE\ given in Section~10.
Since in that context we do not have to specify how many entries
there occur in the first row and how many in the second row,
what we have to do is to compute the expected number of down-steps in a
two-coloured Motzkin path from $(0,0)$ to $(n,t)$ as in Lemma~\TH\
with $f=0$, however without specifying $c$ and~$d$.
In order to do that, we have to put $x=y=1$ in the generating function
in Proposition~\TI,
differentiate the generating function
with respect to~$\al$, subsequently put $\al=1$, extract the
coefficient of~$z^n$, and finally divide
the result by the total number of these Motzkin paths, given by the
formula in Corollary~\TD. The next two lemmas provide the corresponding
coefficient extractions from the second and the third term in~(\AAa).
(The first term may be safely ignored since it does not contain~$\al$.)


\proclaim{Lemma \TQ}
With $x=y=1$, we have
$$
\multline
\frac {\partial} {\partial\al}
\coef{z^n}\frac {\al \left(zM(z)\right)^{t+2}} {(1+zM(z))^2}
\Bigg\vert_{\al=1}
=
\binom {2n-5}{n-t-2}
+\binom {2n-5}{n-t-3}\\
+(n-3)\binom {2n-5}{n-t-4}
-(n+1)\binom {2n-5}{n-t-5}.
\endmultline$$
\endproclaim

\demo{Proof}
Using~(\AAb) with $x=y=1$,
we compute
$$\allowdisplaybreaks
\align
\frac {\partial} {\partial\al}
\coef{z^n}&\frac {\al \left(zM(z)\right)^{t+2}} {(1+zM(z))^2}
\Bigg\vert_{\al=1}
=\frac {\partial} {\partial\al}
\coef{z^{n-t-2}}\frac {\left(1+2z+\al z^2\right)^{n-1}(\al-\al^2 z^2)}
{(1+z)^2}
\Bigg\vert_{\al=1}\\
&=
\coef{z^{n-t-2}}
\frac {(n-1)z^2\left(1+z\right)^{2n-4}(1-z^2)
+\left(1+z\right)^{2n-2}(1-2 z^2)}
{(1+z)^2}\\
&=
\coef{z^{n-t-2}}
 {\left(1+z\right)^{2n-6}
\left((n-1)z^2(1-z^2)+(1+z)^2(1-2z^2)\right)}\\
&=
\coef{z^{n-t-2}}
 {\left(1+z\right)^{2n-5}
 \left(1+z+(n-3)z^2-(n+1)z^3\right)}\\
&=
\binom {2n-5}{n-t-2}
+\binom {2n-5}{n-t-3}
+(n-3)\binom {2n-5}{n-t-4}
-(n+1)\binom {2n-5}{n-t-5},
\endalign$$
which is the claimed formula.\quad \quad \qed
\enddemo

\proclaim{Lemma \TR}
With $x=y=1$, we have
$$\allowdisplaybreaks
\align
\frac {\partial} {\partial\al}&\coef{z^n}\frac {\al zM(z)
\left(\left(zM(z)\right)^t-\left(\frac {z} {1-z}\right)^{t}\right)}
 {\left(1+\al zM(z)\right)
\left(1+zM(z)\right)}
\Bigg\vert_{\al=1}\\
&=\binom {2n-5}{n-t-1}+(n-3)\binom {2n-5}{n-t-3}
-n\binom {2n-5}{n-t-4}-\binom {n-3}{n-t-1}.
\endalign
$$
\endproclaim

\demo{Proof}
Using (\AFa) with $x=y=1$, we obtain
$$\allowdisplaybreaks
\align
\frac {\partial} {\partial\al}&\coef{z^n}\frac {\al zM(z)
\left(\left(zM(z)\right)^t-\left(\frac {z} {1-z}\right)^{t}\right)}
 {\left(1+\al zM(z)\right)
\left(1+zM(z)\right)}
\Bigg\vert_{\al=1}\\
&=\frac {\partial} {\partial\al}
\coef{z^{-1}}
\al 
\left(z^t-\left(\frac {z} {1+z+\al z^2}\right)^t\right)
\frac {\left(1+2z+\al z^2\right)^{n-1}} {z^n}
\frac {1-\al z^2} {(1+\al z)(1+ z)}
\Bigg\vert_{\al=1}\\
&=\frac {\partial} {\partial\al}
\coef{z^{n-t-1}}
\al {\left(1+2z+\al z^2\right)^{n-1}} 
\frac {1-\al z^2} {(1+\al z)(1+ z)}
\Bigg\vert_{\al=1}\\
&\kern1cm
-\frac {\partial} {\partial\al}
\coef{z^{-1}}
\al 
\left(\frac {z} {1+z+\al z^2}\right)^t
\frac {\left(1+2z+\al z^2\right)^{n-1}} {z^n}
\frac {1-\al z^2} {(1+\al z)(1+ z)}
\Bigg\vert_{\al=1}\\
&=\frac {\partial} {\partial\al}
\coef{z^{n-t-1}}
\al {\left(1+2z+\al z^2\right)^{n-1}} 
\frac {1-\al z^2} {(1+\al z)(1+ z)}
\Bigg\vert_{\al=1}\\
&\kern1cm
-\frac {\partial} {\partial\al}\left(\frac {1} {2\pi i}
\int_{\Cal C}
\al 
\left(\frac {z} {1+z+\al z^2}\right)^t
\frac {\left(1+2z+\al z^2\right)^{n-1}} {z^n}
\frac {1-\al z^2} {(1+\al z)(1+ z)}\,dz\right)
\Bigg\vert_{\al=1}.
\endalign$$
where $\Cal C$ is a contour close to the origin, encircling the
origin once in positive direction. 
In the integral, we perform the substitution $v=\frac {z} {1+z+\al z^2}$,
for which $dv=dz\frac {1-\al z^2} {(1+z+\al z^2)^2}$. In this manner,
we get
$$\allowdisplaybreaks
\align
\frac {\partial} {\partial\al}&\coef{z^n}\frac {\al zM(z)
\left(\left(zM(z)\right)^t-\left(\frac {z} {1-z}\right)^{t}\right)}
 {\left(1+\al zM(z)\right)
\left(1+zM(z)\right)}
\Bigg\vert_{\al=1}\\
&=\frac {\partial} {\partial\al}
\coef{z^{n-t-1}}
\al {\left(1+2z+\al z^2\right)^{n-1}} 
\frac {1-\al z^2} {(1+\al z)(1+ z)}
\Bigg\vert_{\al=1}\\
&\kern1cm
-\frac {\partial} {\partial\al}\left(\frac {1} {2\pi i}
\int_{\Cal C}
\al v^{t-2}
\left(\frac {1} {v}+1\right)^{n-1}
\frac {1} {\frac {1} {v}+\al}\,dv\right)
\Bigg\vert_{\al=1}\\
&=\frac {\partial} {\partial\al}
\coef{z^{n-t-1}}
\al {\left(1+2z+\al z^2\right)^{n-1}} 
\frac {1-\al z^2} {(1+\al z)(1+ z)}
\Bigg\vert_{\al=1}\\
&\kern1cm
-\frac {\partial} {\partial\al}
\coef{v^{n-t-1}}
\al
\frac {\left(1+v\right)^{n-1}} {1+\al v}
\Bigg\vert_{\al=1}.
\endalign$$
Now we do the differentiations and put $\al=1$, to obtain
$$\align
\frac {\partial} {\partial\al}&\coef{z^n}\frac {\al zM(z)
\left(\left(zM(z)\right)^t-\left(\frac {z} {1-z}\right)^{t}\right)}
 {\left(1+\al zM(z)\right)
\left(1+zM(z)\right)}
\Bigg\vert_{\al=1}\\
&=
\coef{z^{n-t-1}}
(1+z)^{2n-5}\left(1+(n-3) z^2-n z^3
\right)
-
\coef{v^{n-t-1}}
{\left(1+v\right)^{n-3}} \\
&=\binom {2n-5}{n-t-1}+(n-3)\binom {2n-5}{n-t-3}
-n\binom {2n-5}{n-t-4}-\binom {n-3}{n-t-1},
\endalign$$
which is the desired expression.\quad \quad \qed
\enddemo

\subhead 8. Coefficient extraction, IV\endsubhead
This section provides the basis for the proof of Theorem~\TF\
given in Section~10. Proposition~\TJ\ in Section~3
gives formulae for the generating function for the Motzkin paths
which, by Lemma~\TH, correspond bijectively to the set-valued
standard tableaux in the statement of Theorem~\TF.
The lemmas in this section provide the coefficient extractions of
the individual terms in~(\Aff) and~(\Agg).
We point out that in this section the variables $x,y,\al$ are
again unspecialised.

\proclaim{Lemma \TS}
We have
$$
\coef{z^n}\left(\frac {\al z} {1-yz}\right)^{f}
=\al^f y^{n-f}\binom {n-1} {f-1}.
$$
\endproclaim

\demo{Proof}
This is a simple consequence of the expansion of the binomial
series.\quad \quad \qed
\enddemo

\proclaim{Lemma \TT}
We have
$$\multline
\coef{z^n}
\frac {\al^{f+1} \left(zM(z)\right)^{t+1}} {1+xzM(z)}
\frac {\left(zM(z)\right)^f
-\left(\frac {z} {1-yz}\right)^{f}} {x+\al zM(z)}
=\underset c+d+2e-f+t=n\to{\sum_{c,d,e\ge0}}
x^cy^d\al^{e}\\
\cdot\left(\frac {(f+t)\,(n-1)!}
{(c+e-f)(n-d-e+f)\,(n-c-d-e+f-1)!\,c!\,d!\,(e-f-1)!}\right.\\
-\left.
\frac {t\,(n-1)!\,(n-d-f-1)!}
{(c+e-f)(n-d-e)\,(n-d-1)!\,(n-c-d-e-1)!\,c!\,d!\,(e-f-1)!}
\right).
\endmultline$$

\endproclaim

\demo{Proof}
With (\Ahh) and (\Aii) in mind, by the Lagrange inversion formula
in Theorem~\TK\ with $F(z)=zM(z)$,
the coefficient of $z^n$ in the statement of the lemma equals
$$\allowdisplaybreaks
\align
\coef{z^n}&
\frac {\al^{f+1} \left(zM(z)\right)^{t+1}} {1+xzM(z)}
\frac {\left(zM(z)\right)^f
-\left(\frac {z} {1-yz}\right)^{f}} {x+\al zM(z)}\\
&=\coef{z^{-1}}
\frac {\al^{f+1} z^{t+1}
\left(z^f-\left(\frac {z} {1+xz+\al z^2}\right)^f\right)}
{(1+xz)(x+\al z)}\cdot
\frac {\left(1+(y+x)z+\al z^2\right)^{n+1}} {z^{n+1}}\\
&\kern7cm
\cdot\frac {1-\al z^2} {\left(1+(y+x)z+\al z^2\right)^2}\\
&=\coef{z^{n-t-1}}
\al^{f+1} \frac{\displaystyle
\left(z^f-\left(\frac {z} {1+xz+\al z^2}\right)^f\right)}
{x+\al z}\cdot
\frac {\left(1+(y+x)z+\al z^2\right)^{n-1}\left(1-\al z^2\right)}
{1+xz}.
\endalign$$
Partial fraction decomposition yields
$$
\frac {1} {(1+xz)(x+\al z)}
=\frac {1} {1-\al z^2}\left(\frac {1} {x+\al z}-\frac {z} {1+xz}\right).
$$
Substituting this in the last expression, we obtain
$$\allowdisplaybreaks
\align
\coef{z^{n-t-1}}&
\al^{f+1} 
\left(z^f-\left(\frac {z} {1+xz+\al z^2}\right)^f\right)
\left(1+(y+x)z+\al z^2\right)^{n-1}
\left(\frac {1} {x+\al z}-\frac {z} {1+xz}\right)\\
&=
\coef{z^{n-t-1}}
\al^{f+1} z^f
\left(1-\left(1+xz+\al z^2\right)^{-f}\right)\\
&\kern2cm
\cdot
\sum_{c\ge0}\binom {n-1}c y^cz^c
\left(1+xz+\al z^2\right)^{n-1-c}
\left(\frac {1} {x+\al z}-\frac {z} {1+xz}\right)\\
&=
\coef{z^{n-t-1}}
\al^{f+1} \sum_{c\ge0}\binom {n-1}c y^cz^{c+f}\\
&\kern1cm
\cdot
\left(
\sum_{k\ge0}\left(\binom {n-1-c}k-\binom {n-1-c-f}k\right) z^k(x+\al z)^{k-1}
\right.\\
&\kern3cm
\left.
-\sum_{e\ge0}\binom {n-1-c}e
\al^ez^{2e+1}(1+x z)^{n-c-e-2}
\right.\\
&\kern3cm
\left.
+\sum_{e\ge0}\binom {n-1-c-f}e
\al^ez^{2e+1}(1+x z)^{n-c-f-e-2}
\right)\\
&=
\coef{z^{n-t-1}}
\al^{f+1} \sum_{c\ge0}\binom {n-1}c y^cz^{c+f}\\
&\kern1cm
\cdot
\left(
\sum_{k\ge0}\left(\binom {n-1-c}k-\binom {n-1-c-f}k\right) z^k
\sum_{e\ge0}\binom {k-1}e x^{k-1-e}\al^ez^e
\right.\\
&\kern2.5cm
\left.
-\sum_{e\ge0}\binom {n-1-c}e
\al^ez^{2e+1}\sum_{d\ge0}\binom {n-c-e-2}d x^dz^d
\right.\\
&\kern2.5cm
\left.
+\sum_{e\ge0}\binom {n-1-c-f}e
\al^ez^{2e+1}\sum_{d\ge0}\binom {n-c-f-e-2}d x^dz^d
\right).
\endalign$$
Thus, we get
$$\allowdisplaybreaks
\align
&\underset c+d+2(e+1)+f=n-t\to{\sum_{c,d,e\ge0}}
y^cx^d\al^{e+f+1}
\binom {n-1}c \\
&\kern1cm
\cdot
\left(
\binom {n-1-c}{d+e+1}\binom {d+e}e 
-\binom {n-1-c-f}{d+e+1}\binom {d+e}e 
\right.\\
&\kern2cm
\left.
-\binom {n-1-c}e\binom {n-c-e-2}d 
+\binom {n-1-c-f}e\binom {n-c-f-e-2}d
\right)\\
&=\underset c+d+2(e+1)+f+t=n\to{\sum_{c,d,e\ge0}}
y^cx^d\al^{e+f+1}
\binom {n-1}c \\
&\kern1cm
\cdot
\left(
\frac {(n-c-1)!} {(n-c-d-e-2)!\,d!\,e!}
\left(\frac {1} {d+e+1}-\frac {1} {n-c-e-1}\right)
\right.\\
&\kern1.5cm
\left.
-\frac {(n-c-f-1)!} {(n-c-d-e-f-2)!\,d!\,e!}
\left(\frac {1} {d+e+1}-\frac {1} {n-c-e-f-1}\right)
\right)\\
&=\underset c+d+2(e+1)+f+t=n\to{\sum_{c,d,e\ge0}}
y^cx^d\al^{e+f+1}\left(
\frac {(f+t)\,(n-1)!} {(d+e+1)(n-c-e-1)\,(n-c-d-e-2)!\,c!\,d!\,e!}
\right.\\
&\kern1.5cm \left.
-\frac {t\,(n-1)!\,(n-c-f-1)!}
{(d+e+1)(n-c-e-f-1)\,(n-c-1)!\,(n-c-d-e-f-2)!\,c!\,d!\,e!}
\right).
\endalign$$
After replacement of $e$ by $e-f-1$ and interchange of $c$ and~$d$,
we arrive at the expression in the statement of the lemma.\quad \quad \qed
\enddemo

\proclaim{Lemma \TU}
We have
$$\allowdisplaybreaks
\align
\coef{z^n}&
\al^{f-t+1} zM(z)
\frac {\left(zM(z)\right)^{f-t}
-\left(\frac {z} {1-yz}\right)^{f-t}} {x+\al zM(z)}
\frac {1-\left(\al z^2M^2(z)\right)^t} {1-\al z^2M^2(z)}\\
&=
\underset c+d+2e-f+t=n\to{\sum_{c,d,e\ge0}}
x^cy^d\al^{e}
\left(\frac {(n-1)!} {(c+e-f+t)\cdot c!\,d!\,(e-f+t-1)!\,
(e-1)!}\right.\\
&\kern2cm
\left.
-\frac {(n-1)!\,(n-d-f+t-1)!} {(c+e-f+t)\cdot c!\,d!\,(n-d-1)!\,(e-f+t-1)!^2}
\right.\\
&\kern1.5cm
-\left.
\frac {(n-1)!} {(c+e-f)\cdot c!\,d!\,(e-f-1)!\,(e+t-1)!}\right.\\
&\kern2cm
\left.
+\frac {(n-1)!\,(n-d-f+t-1)!} {(c+e-f)\cdot c!\,d!\,(n-d-1)!\,(e-f-1)!\,
(e-f+2t-1)!}
\right).
\endalign$$
\endproclaim

\demo{Proof}
Again with (\Ahh) and (\Aii) in mind, by the Lagrange inversion formula
in Theorem~\TK\ with $F(z)=zM(z)$,
the coefficient of $z^n$ in the statement of the lemma equals
$$\allowdisplaybreaks
\align
\coef{z^n}&
\al^{f-t+1} zM(z)
\frac {\left(zM(z)\right)^{f-t}
-\left(\frac {z} {1-yz}\right)^{f-t}} {x+\al zM(z)}
\frac {1-\left(\al z^2M^2(z)\right)^t} {1-\al z^2M^2(z)}\\
&=\coef{z^{-1}}
\al^{f-t+1} z
\frac {z^{f-t}
-\left(\frac {z} {1+xz+\al z^2}\right)^{f-t}} {x+\al z}
\frac {1-\left(\al z^2\right)^t} {1-\al z^2}\\
&\kern2cm
\cdot
\frac {(1+(x+y)z+\al z^2)^{n+1}} {z^{n+1}}
\frac {1-\al z^2} {(1+(x+y)z+\al z^2)^2}\\
&=\coef{z^{n-1}}
\al^{f-t+1} 
\left({z^{f-t}
-\left(\frac {z} {1+xz+\al z^2}\right)^{f-t}} \right)
\left(1-\left(\al z^2\right)^t\right)\\
&\kern4cm
\cdot
\frac {(1+(x+y)z+\al z^2)^{n-1}} {x+\al z}\\
&=\coef{z^{n-f+t-1}}
\al^{f-t+1} 
\left(1
-\left(1+xz+\al z^2\right)^{t-f} \right)
\frac {\left(1-\left(\al z^2\right)^t\right)} {x+\al z}\\
&\kern4cm
\cdot
\sum_{d\ge0}\binom {n-1}d y^dz^d(1+xz+\al z^2)^{n-1-d}\\
&=\coef{z^{n-f+t-1}}
\al^{f-t+1} 
\left(1-\left(\al z^2\right)^t\right) 
\sum_{d\ge0}\binom {n-1}d y^dz^d\\
&\kern2cm
\cdot
\sum_{k\ge0}\left(\binom {n-d-1}k-\binom {n-d-f+t-1}k\right)
z^k(x+\al z)^{k-1}\\
&=\coef{z^{n-f+t-1}}
\al^{f-t+1} 
\left(1-\left(\al z^2\right)^t\right) 
\\
&\kern1cm
\cdot
\sum_{d,k\ge0}\binom {n-1}d
\left(\binom {n-d-1}k-\binom {n-d-f+t-1}k\right)
\\
&\kern2cm
\cdot
\sum_{c\ge0}\binom {k-1}cx^cy^d\al^{k-1-c}z^{d+2k-1-c}\\
&=
\underset c+d+2e-f+t=n\to{\sum_{c,d,e\ge0}}
x^cy^d\al^{e}\\
&\kern1cm
\cdot
\left(\frac {(n-1)!} {(c+e-f+t)\cdot c!\,d!\,(e-f+t-1)!\,(n-c-d-e+f-t-1)!}\right.\\
&\kern2cm
\left.
-\frac {(n-1)!\,(n-d-f+t-1)!} {(c+e-f+t)\cdot c!\,d!\,(n-d-1)!\,(e-f+t-1)!\,(n-c-d-e-1)!}
\right)\\
&\kern1cm
-\underset c+d+2e-f+t=n\to{\sum_{c,d,e\ge0}}
x^cy^d\al^{e}
\left(\frac {(n-1)!} {(c+e-f)\cdot c!\,d!\,(e-f-1)!\,(n-c-d-e+f-1)!}\right.\\
&\kern2cm
\left.
-\frac {(n-1)!\,(n-d-f+t-1)!} {(c+e-f)\cdot c!\,d!\,(n-d-1)!\,(e-f-1)!\,(n-c-d+t-e-1)!}
\right)\\
&=
\underset c+d+2e-f+t=n\to{\sum_{c,d,e\ge0}}
x^cy^d\al^{e}
\left(\frac {(n-1)!} {(c+e-f+t)\cdot c!\,d!\,(e-f+t-1)!\,
(e-1)!}\right.\\
&\kern2cm
\left.
-\frac {(n-1)!\,(n-d-f+t-1)!} {(c+e-f+t)\cdot c!\,d!\,(n-d-1)!\,(e-f+t-1)!^2}
\right.\\
&\kern1.5cm
-\left.
\frac {(n-1)!} {(c+e-f)\cdot c!\,d!\,(e-f-1)!\,(e+t-1)!}\right.\\
&\kern2cm
\left.
+\frac {(n-1)!\,(n-d-f+t-1)!} {(c+e-f)\cdot c!\,d!\,(n-d-1)!\,(e-f-1)!\,
(e-f+2t-1)!}
\right),
\endalign$$
which agrees with the expression in the statement of the
lemma.\quad \quad \qed
\enddemo

\proclaim{Lemma \TV}
We have
$$\allowdisplaybreaks
\align
\coef{z^n}&
\left(\frac {\al z} {1-yz}\right)^{f-t}
\frac {1-\left(\frac {\al z^2M(z)} {1-yz}\right)^{t}}
{1+xzM(z)}
\frac {\al \left(zM(z)\right)^2} {1-\al z^2M^2(z)}\\
&=
\underset c+d+2e-f+t=n\to{\sum_{c,d,e\ge0}}
x^cy^d\al^{e}
\left(\frac {(n-1)!\,(n-d-f+t-1)!}
{(n-d-e)\,c!\,d!\,(n-d-1)!\,(e-f+t-1)!^2}
\right.\\
&\kern2cm
\left.
-\frac {(n-1)!\,(n-d-f-1)!}
{(n-d-e)\,c!\,d!\,(n-d-1)!\,(e-f-1)!\,(e-f+t-1)!}
\right).
\endalign$$
\endproclaim

\demo{Proof}
By the Lagrange inversion formula
in Theorem~\TK\ with $F(z)=zM(z)$,
the coefficient of $z^n$ in the statement of the lemma equals
$$\allowdisplaybreaks
\align
\coef{z^n}&
\left(\frac {\al z} {1-yz}\right)^{f-t}
\frac {1-\left(\frac {\al z^2M(z)} {1-yz}\right)^{t}}
{1+xzM(z)}
\frac {\al \left(zM(z)\right)^2} {1-\al z^2M^2(z)}\\
&=
\coef{z^{-1}}
\left(\frac {\al z} {1+xz+\al z^2}\right)^{f-t}
\frac {1-\left(\frac {\al z^2} {1+xz+\al z^2}\right)^{t}}
{1+xz}
\frac {\al z^2} {1-\al z^2}\\
&\kern2cm
\cdot
\frac {(1+(x+y)z+\al z^2)^{n+1}} {z^{n+1}}\frac {1-\al z^2}
{(1+(x+y)z+\al z^2)^2}\\
&=
\coef{z^{n-2}}
\al\frac {\left(\frac {\al z} {1+xz+\al z^2}\right)^{f-t}
-z^t\left(\frac {\al z} {1+xz+\al z^2}\right)^{f}}
{1+xz}
(1+(x+y)z+\al z^2)^{n-1}\\
&=
\coef{z^{n-f+t-2}}
\al^{f-t+1}\sum_{d\ge0}\binom {n-1}d y^dz^d
(1+xz+\al z^2)^{n-d-f+t-1}
(1+xz)^{-1}\\
&\kern1cm
-\coef{z^{n-f-t-2}}
\al^{f+1}\sum_{d\ge0}\binom {n-1}d y^dz^d
(1+xz+\al z^2)^{n-d-f-1}
(1+xz)^{-1}\\
&=
\coef{z^{n-f+t-2}}
\al^{f-t+1}\sum_{d,e\ge0}\binom {n-1}d y^dz^d
\binom {n-d-f+t-1}e\al^ez^{2e}\\
&\kern5cm
\cdot
(1+xz)^{n-d-e-f+t-2}\\
&\kern1cm
-\coef{z^{n-f-t-2}}
\al^{f+1}\sum_{d,e\ge0}\binom {n-1}d y^dz^d
\binom {n-d-f-1}e\al^ez^{2e}\\
&\kern5cm
\cdot
(1+xz)^{n-d-e-f-2}\\
&=
\coef{z^{n-f+t-2}}
\al^{f-t+1}\sum_{d,e\ge0}\binom {n-1}d y^dz^d
\binom {n-d-f+t-1}e\al^ez^{2e}\\
&\kern5cm
\cdot
\sum_{c\ge0}\binom {n-d-e-f+t-2}c x^c z^c\\
&\kern1cm
-\coef{z^{n-f-t-2}}
\al^{f+1}\sum_{d,e\ge0}\binom {n-1}d y^dz^d
\binom {n-d-f-1}e\al^ez^{2e}\\
&\kern5cm
\cdot
\sum_{c\ge0}\binom {n-d-e-f-2}c x^cz^c\\
&=
\underset c+d+2e=n-f+t-2\to{\sum_{c,d,e\ge0}}
x^cy^d\al^{e+f-t+1}\\
&\kern2cm
\cdot
\binom {n-1}d\binom {n-d-f+t-1}e 
\binom {n-d-e-f+t-2}c  \\
&\kern1cm
-
\underset c+d+2e=n-f-t-2\to{\sum_{c,d,e\ge0}}x^cy^d\al^{e+f+1}\\
&\kern2cm
\cdot
\binom {n-1}d 
\binom {n-d-f-1}e 
\binom {n-d-e-f-2}c \\
&=
\underset c+d+2e-f+t=n\to{\sum_{c,d,e\ge0}}
x^cy^d\al^{e}
\left(\frac {(n-1)!\,(n-d-f+t-1)!}
{(n-d-e)\,c!\,d!\,(n-d-1)!\,(e-f+t-1)!^2}
\right.\\
&\kern2cm
\left.
-\frac {(n-1)!\,(n-d-f-1)!}
{(n-d-e)\,c!\,d!\,(n-d-1)!\,(e-f-1)!\,(e-f+t-1)!}
\right),
\endalign$$
as desired.\quad \quad \qed
\enddemo

\proclaim{Lemma \TW}
We have
$$\allowdisplaybreaks
\align
\coef{z^n}&
\al^{f+1}\left(\frac { z} {1-yz}\right)^{f-t}
\frac {\left(\frac { z} {1-yz}\right)^{t}
-\left( zM(z)\right)^{t}} {x+\al zM(z)}
\frac {
\left(zM(z)\right)^{t+1}} {1-\al z^2M^2(z)}\\
&=
\underset c+d+2e-f+t=n\to{\sum_{c,d,e\ge0}}
x^cy^d\al^{e}
\left(\frac {(n-1)!\,(n-d-f-1)!}
{(c+e-f)\,c!\,d!\,(n-d-1)!\,(e-f+t-1)!\,(e-f-1)!}
\right.\\
&\kern2cm
\left.
-\frac {(n-1)!\,(n-d-f+t-1)!}
{(c+e-f)\,c!\,d!\,(n-d-1)!\,(e-f+2t-1)!\,(e-f-1)!}
\right) .
\endalign$$
\endproclaim

\demo{Proof}
By the Lagrange inversion formula
in Theorem~\TK\ with $F(z)=zM(z)$,
the coefficient of $z^n$ in the statement of the lemma equals
$$\allowdisplaybreaks
\align
\coef{z^n}&
\al^{f+1}\left(\frac { z} {1-yz}\right)^{f-t}
\frac {\left(\frac { z} {1-yz}\right)^{t}
-\left( zM(z)\right)^{t}} {x+\al zM(z)}
\frac {
\left(zM(z)\right)^{t+1}} {1-\al z^2M^2(z)}\\
&=
\coef{z^{-1}}
\al^{f+1}\left(\frac { z} {1+xz+\al z^2}\right)^{f-t}
\frac {\left(\frac { z} {1+xz+\al z^2}\right)^{t}
-z^{t}} {x+\al z}
\frac {z^{t+1}} {1-\al z^2}\\
&\kern1cm
\cdot
\frac {(1+(x+y)z+\al z^2)^{n+1}} {z^{n+1}}
\frac {1-\al z^2} {(1+(x+y)z+\al z^2)^2}\\
&=
\coef{z^{n-f-t-1}}
\al^{f+1}
\left(1+xz+\al z^2\right)^{-f}
(1+(x+y)z+\al z^2)^{n-1}(x+\al z)^{-1}\\
&\kern1cm
-\coef{z^{n-f-t-1}}
\al^{f+1}
\left(1+xz+\al z^2\right)^{t-f}
(1+(x+y)z+\al z^2)^{n-1}(x+\al z)^{-1}\\
&=
\coef{z^{n-f-t-1}}
\al^{f+1}\sum_{d\ge0}\binom {n-1}d y^dz^d
\left(1+xz+\al z^2\right)^{n-d-f-1}
(x+\al z)^{-1}\\
&\kern1cm
-\coef{z^{n-f-t-1}}
\al^{f+1}\sum_{d\ge0}\binom {n-1}d y^dz^d
\left(1+xz+\al z^2\right)^{n-d-f+t-1}
(x+\al z)^{-1}\\
&=
\coef{z^{n-f-t-1}}
\al^{f+1}\sum_{d,k\ge0}\binom {n-1}d y^dz^d
\binom {n-d-f-1}k z^k
(x+\al z)^{k-1}\\
&\kern1cm
-\coef{z^{n-f-t-1}}
\al^{f+1}\sum_{d,k\ge0}\binom {n-1}d y^dz^d
\binom {n-d-f+t-1}k z^k
(x+\al z)^{k-1}\\
&=
\coef{z^{n-f-t-1}}
\al^{f+1}\sum_{c,d,k\ge0}\binom {n-1}d y^dz^d
\binom {n-d-f-1}k z^k\\
&\kern5cm
\cdot
\binom {k-1}c x^c\al^{k-c-1}z^{k-c-1}\\
&\kern1cm
-\coef{z^{n-f-t-1}}
\al^{f+1}\sum_{c,d,k\ge0}\binom {n-1}d y^dz^d
\binom {n-d-f+t-1}k z^k\\
&\kern5cm
\cdot
\binom {k-1}c x^c\al^{k-c-1}z^{k-c-1}\\
&=
\underset c+d+2e-f+t=n\to{\sum_{c,d,e\ge0}}
x^cy^d\al^{e}
\left(\frac {(n-1)!\,(n-d-f-1)!}
{(c+e-f)\,c!\,d!\,(n-d-1)!\,(e-f+t-1)!\,(e-f-1)!}
\right.\\
&\kern2cm
\left.
-\frac {(n-1)!\,(n-d-f+t-1)!}
{(c+e-f)\,c!\,d!\,(n-d-1)!\,(e-f+2t-1)!\,(e-f-1)!}
\right) ,
\endalign$$
as we claimed.\quad \quad \qed
\enddemo

\proclaim{Lemma \TX}
We have
$$\multline
\coef{z^n}
\frac {\al^{f+1}\left(zM(z)\right)^{f+t+2}} 
{(1+yzM(z))(1+xzM(z))}
\\
=\underset c+d+2e+t-f=n\to{\sum_{c,d,e\ge0}}
x^cy^d\al^{e}
\left(\frac {(n-1)!} {(n-c-e+f)\,c!\,d!\,(e-f-1)!\,(n-c-d-e+f-1)!}
\vphantom{\sum_{b=d+1}^{n-e+f-1}}\right.\\
\left.
-\sum_{b=n-c-e+f+1}^{n-e+f+1}(-1)^{n-b-c-e+f+1}
\frac {(n-b)\cdot(n-1)!} {b\cdot d!\,(b-1-d)!\,
(e-f-1)!\,(n-b-e+f+1)!}\right).
\endmultline
$$
\endproclaim

\demo{Proof}
This follows directly from Lemma~\TM\ by replacing $t$ by $f+t$
and multiplying everything by~$\al^f$.\quad \quad \qed
\enddemo

\proclaim{Lemma \TY}
We have
$$\allowdisplaybreaks
\align
\coef{z^n}&
\frac {1-\left(\al z^2M^2(z)\right)^{t}}
{1-\al z^2M^2(z)} 
\cdot
\frac {\al^{f-t+1}\left( zM(z)\right)^{f-t+2}}
{1+yzM(z)}\\
&=
\underset c+d+2e-f+t=n\to{\sum_{c,d,e\ge0}}
x^c y^d\al^{e}
\left(\frac {(n-1)!}
{(n-c-e+f-t)\,c!\,d!\,(e-f+t-1)!\,(e-1)!}
\right.\\
&\kern4cm
\left.
-\frac {(n-1)!}
{(n-c-e+f)\,c!\,d!\,(e-f-1)!\,(e+t-1)!}
\right).
\endalign$$
\endproclaim

\demo{Proof}
By the Lagrange inversion formula
in Theorem~\TK\ with $F(z)=zM(z)$,
the coefficient of $z^n$ in the statement of the lemma equals
$$\allowdisplaybreaks
\align
\coef{z^n}&
\frac {1-\left(\al z^2M^2(z)\right)^{t}}
{1-\al z^2M^2(z)} 
\cdot
\frac {\al^{f-t+1}\left( zM(z)\right)^{f-t+2}}
{1+yzM(z)}\\
&=
\coef{z^{-1}}
\al^{f-t+1}
\frac {1-\left(\al z^2\right)^{t}} {1-\al z^2} 
\cdot
\frac {z^{f-t+2}}
{1+yz}\\
&\kern2cm
\cdot
\frac {(1+(x+y)z+\al z^2)^{n+1}} {z^{n+1}}
\frac {1-\al z^2} {(1+(x+y)z+\al z^2)^2}\\
&=
\coef{z^{n-f+t-2}}
\al^{f-t+1}
(1+(x+y)z+\al z^2)^{n-1}
(1+yz)^{-1}\\
&\kern1cm
-\coef{z^{n-f-t-2}}
\al^{f+1}
(1+(x+y)z+\al z^2)^{n-1}
(1+yz)^{-1}\\
&=
\coef{z^{n-f+t-2}}
\al^{f-t+1}
\sum_{c,e\ge0}\binom {n-1}{c,e,n-c-e-1} x^cz^c \al^ez^{2e}
(1+yz)^{n-c-e-2}\\
&\kern1cm
-\coef{z^{n-f-t-2}}
\al^{f+1}
\sum_{c,e\ge0}\binom {n-1}{c,e,n-c-e-1} x^cz^c \al^ez^{2e}
(1+yz)^{n-c-e-2}\\
&=
\coef{z^{n-f+t-2}}
\sum_{c,d,e\ge0}\binom {n-1}{c,e,n-c-e-1} x^cz^c \al^{e+f-t+1}z^{2e}
\binom {n-c-e-2}d y^dz^d\\
&\kern1cm
-\coef{z^{n-f-t-2}}
\sum_{c,d,e\ge0}\binom {n-1}{c,e,n-c-e-1} x^cz^c \al^{e+f+1}z^{2e}
\binom {n-c-e-2}d y^dz^d\\
&=
\underset c+d+2e-f+t=n\to{\sum_{c,d,e\ge0}}
x^c y^d\al^{e}
\left(\frac {(n-1)!}
{(n-c-e+f-t)\,c!\,d!\,(e-f+t-1)!\,(e-1)!}
\right.\\
&\kern4cm
\left.
-\frac {(n-1)!}
{(n-c-e+f)\,c!\,d!\,(e-f-1)!\,(e+t-1)!}
\right),
\endalign$$
which is the claimed formula.\quad \quad \qed
\enddemo

\proclaim{Lemma \TZ}
We have
$$\allowdisplaybreaks
\align
\coef{z^n}&
\frac {\al\left(zM(z)\right)^{t-f+2}-\al^{f+1}\left(zM(z)\right)^{f+t+2}}
{(1+xzM(z))(1-\al z^2M^2(z))}\\
&=
\underset c+d+2e-f+t=n\to{\sum_{c,d,e\ge0}}
x^cy^d\al^{e}
\left(\frac {(n-1)!}
{(n-d-e)\,c!\,d!\,(e-1)!\,(e-f+t-1)!}
\right.\\
&\kern4cm
\left.
-
\frac {(n-1)!}
{(n-d-e+f)\,c!\,d!\,(e-f-1)!\,(e+t-1)!}
\right).
\endalign$$
\endproclaim

\demo{Proof}
By the Lagrange inversion formula
in Theorem~\TK\ with $F(z)=zM(z)$,
the coefficient of $z^n$ in the statement of the lemma equals
$$\allowdisplaybreaks
\align
\coef{z^n}&
\frac {\al\left(zM(z)\right)^{t-f+2}-\al^{f+1}\left(zM(z)\right)^{f+t+2}}
{(1+xzM(z))(1-\al z^2M^2(z))}\\
&=
\coef{z^{-1}}
\frac {\al z^{t-f+2}-\al^{f+1}z^{f+t+2}}
{(1+xz)(1-\al z^2)}
\cdot
\frac {(1+(x+y)z+\al z^2)^{n+1}} {z^{n+1}}
\frac {1-\al z^2} {(1+(x+y)z+\al z^2)^2}\\
&=
\coef{z^{n+f-t-2}}
\al
(1+(x+y)z+\al z^2)^{n-1}({1+xz})^{-1}\\
&\kern1cm
-
\coef{z^{n-f-t-2}}
\al^{f+1} 
(1+(x+y)z+\al z^2)^{n-1}({1+xz})^{-1}\\
&=
\coef{z^{n+f-t-2}}
\al
\sum_{d,e\ge0}\binom {n-1}{d,e,n-d-e-1}y^dz^d\al^ez^{2e}
(1+xz)^{n-d-e-2}\\
&\kern1cm
-
\coef{z^{n-f-t-2}}
\al^{f+1} 
\sum_{d,e\ge0}\binom {n-1}{d,e,n-d-e-1}y^dz^d\al^e
z^{2e}
(1+xz)^{n-d-e-2}\\
&=
\underset c+d+2e=n+f-t-2\to{\sum_{c,d,e\ge0}}
x^cy^d\al^{e+1}
\binom {n-1}{d,e,n-d-e-1}
\binom {n-d-e-2}c \\
&\kern1cm
-
\underset c+d+2e=n-f-t-2\to{\sum_{c,d,e\ge0}}
x^cy^d\al^{e+f+1}
\binom {n-1}{d,e,n-d-e-1}
\binom {n-d-e-2}c \\
&=
\underset c+d+2e-f+t=n\to{\sum_{c,d,e\ge0}}
x^cy^d\al^{e}
\left(\frac {(n-1)!}
{(n-d-e)\,c!\,d!\,(e-1)!\,(e-f+t-1)!}
\right.\\
&\kern4cm
\left.
-
\frac {(n-1)!}
{(n-d-e+f)\,c!\,d!\,(e-f-1)!\,(e+t-1)!}
\right),
\endalign$$
as required.\quad \quad \qed
\enddemo

\proclaim{Lemma \TTA}
We have
$$\multline
\coef{z^n}\frac {\al zM(z)\left(\left(zM(z)\right)^{t-f}
-\left(\frac {z} {1-xz}\right)^{t-f}\right)}
{(y+\al zM(z))(1+yzM(z))} \\
=\underset c+d+2e+t-f=n\to{\sum_{c,d,e\ge0}}
x^cy^d\al^{e}\frac {(t-f)\,(n-1)!} {(d+e)(n-c-e)\,(n-c-d-e-1)!\,c!\,d!\,(e-1)!}.
\endmultline
$$
\endproclaim

\demo{Proof}
This follows from Lemma~\TN\ by replacing $t$ by $t-f$.\quad \quad \qed
\enddemo

\proclaim{Lemma \TTB}
We have
$$\allowdisplaybreaks
\align
\coef{z^n}&\frac {\al \left(zM(z)\right)^{t-f+2}}
{1+yzM(z)} 
\cdot
\frac {1-\left(\al z^2M^2(z)\right)^f}
{1-\al z^2M^2(z)} \\
&=
\underset c+d+2e-f+t=n\to
{\sum_{c,d,e\ge0}}x^cy^d\al^{e}
\left(\frac {(n-1)!}
{(d+e-f+t)\,c!\,d!\,(e-1)!\,(e-f+t-1)!}
\right.\\
&\kern4cm
-\left.
\frac {(n-1)!}
{(d+e+t)\,c!\,d!\,(e-f-1)!\,(e+t-1)!}
\right).
\endalign$$
\endproclaim

\demo{Proof}
This follows easily from Lemma~\TY\ by interchanging the roles
of $f$ and~$t$ and multiplying the result by $\al^{f-t}$.\quad \quad \qed
\enddemo

\subhead 9. Coefficient extraction, V\endsubhead
In this section, we apply the Lagrange inversion formula
in Theorem~\TK\ in order to provide the basis for the proof of
Theorem~\TG\ given in Section~10.
In principle, we could sum~(\AAF) respectively~(\AAG)
over all $c,d,e$ with $c+d+2e-f+t=n$. However again,
it is more convenient to extract the desired number
directly from the (specialised) generating function.
This generating function was given
in Proposition~\TJ.

\medskip
We begin with the case where $t<f$.
So we have to consider~(\Aff) with $x=y=\al=1$,
that is,
$$\allowdisplaybreaks
\align
&\chi(t<f)\left(\frac { z} {1-z}\right)^{f-t}
+ \left(zM(z)\right)^{t+1}
\frac {\left(zM(z)\right)^f
-\left(\frac {z} {1-z}\right)^{f}} {(1+ zM(z))^2}\\
&\kern1cm
+
 zM(z)
\frac {\left(zM(z)\right)^{f-t}
-\left(\frac {z} {1-z}\right)^{f-t}} {1+zM(z)}
\frac {1-\left( z^2M^2(z)\right)^t} {1- z^2M^2(z)}\\
&\kern1cm
+\left(\frac {z} {1-z}\right)^{f-t}
\frac {1-\left(\frac {z^2M(z)} {1-z}\right)^{t}}
{1+zM(z)}
\frac {\left(zM(z)\right)^2} {1- z^2M^2(z)}\\
&\kern2cm
-\left(\frac { z} {1-z}\right)^{f-t}
\frac {\left( zM(z)\right)^{t}
-\left(\frac { z} {1-z}\right)^{t}} {1+ zM(z)}
\frac {
\left(zM(z)\right)^{t+1}} {1-z^2M^2(z)}\\
&\kern1cm
+\frac { \left(zM(z)\right)^{f+t+2}} 
{(1+zM(z))^2}
+
\frac {1-\left(z^2M^2(z)\right)^{t}}
{1-z^2M^2(z)} 
\cdot
\frac {\left( zM(z)\right)^{f-t+2}}
{1+zM(z)}.
\tag\AL
\endalign$$
Lemmas~\TTC--\TTH\ below afford the coefficient extraction from the
individual terms in~(\AL).

\proclaim{Lemma \TTC}
With $x=y=\al=1$, we have
$$\align
\coef{z^n}&
\left(zM(z)\right)^{t+1}
\frac {\left(zM(z)\right)^f
-\left(\frac {z} {1-z}\right)^{f}} {(1+ zM(z))^2}\\
&=
\sum_{k\ge f+1}(-1)^{k-f+1}\binom {k-1}{f-1}
\left(\binom {2n-3+k-f}{n-k-t-1}-\binom {2n-3+k-f}{n-k-t-2}\right).
\endalign$$
\endproclaim

\demo{Proof}
By the Lagrange inversion formula in Theorem~\TK\ with $F(z)=zM(z)$,
and the earlier considerations concerning the implied $f(z)$
under the specialisation $x=y=\al=1$
(see the paragraph above Lemma~\TO), we have
$$\align
\coef{z^n}&
\left(zM(z)\right)^{t+1}
\frac {\left(zM(z)\right)^f
-\left(\frac {z} {1-z}\right)^{f}} {(1+ zM(z))^2}\\
&=
\coef{z^{-1}}
z^{t+1}
\frac {z^f
-\left(\frac {z} {1+z+z^2}\right)^{f}} {(1+ z)^2}
\frac {(1+z)^{2(n+1)}} {z^{n+1}}\frac {1-z} {(1+z)^3}\\
&=
\coef{z^{n-t-1}}
\left(z^f
-\left(\frac {z} {1+z+z^2}\right)^{f} \right)
(1+z)^{2n-3}(1-z)\\
&=
\coef{z^{n-t-1}}
\sum_{k\ge f+1}(-1)^{k-f+1}\binom {k-1}{f-1}z^k
(1+z)^{2n-3+k-f}(1-z)\\
&=
\sum_{k\ge f+1}(-1)^{k-f+1}\binom {k-1}{f-1}
\left(\binom {2n-3+k-f}{n-k-t-1}-\binom {2n-3+k-f}{n-k-t-2}\right),
\endalign$$
as desired.\quad \quad \qed
\enddemo

\proclaim{Lemma \TTD}
With $x=y=\al=1$, we have
$$
\multline
\coef{z^n}
 zM(z)
\frac {\left(zM(z)\right)^{f-t}
-\left(\frac {z} {1-z}\right)^{f-t}} {1+zM(z)}
\frac {1-\left( z^2M^2(z)\right)^t} {1- z^2M^2(z)}\\
=
\sum_{k\ge f-t+1}(-1)^{k-f+t+1}\binom {k-1}{f-t-1}
\left(\binom {2n+k-f+t-3}{n-k-1}-\binom {2n+k-f+t-3}{n-k-2t-1}\right).
\endmultline$$
\endproclaim

\demo{Proof}
By the Lagrange inversion formula in Theorem~\TK\ with $F(z)=zM(z)$,
and the earlier considerations concerning the implied $f(z)$
under the specialisation $x=y=\al=1$
(see the paragraph above Lemma~\TO), we have
$$\allowdisplaybreaks
\align
\coef{z^n}&
 zM(z)
\frac {\left(zM(z)\right)^{f-t}
-\left(\frac {z} {1-z}\right)^{f-t}} {1+zM(z)}
\frac {1-\left( z^2M^2(z)\right)^t} {1- z^2M^2(z)}\\
&=\coef{z^{-1}}
 z
\frac {z^{f-t}
-\left(\frac {z} {1+z+z^2}\right)^{f-t}} {1+z}
\frac {1-z^{2t}} {1- z^2}
\frac {(1+z)^{2(n+1)}} {z^{n+1}}\frac {1-z} {(1+z)^3}\\
&=\coef{z^{n-1}}
\left(z^{f-t}
-\left(\frac {z} {1+z+z^2}\right)^{f-t}\right)
(1-z^{2t})
(1+z)^{2n-3}\\
&=\coef{z^{n-1}}
\sum_{k\ge f-t+1}(-1)^{k-f+t+1}\binom {k-1}{f-t-1}
z^k
(1+z)^{2n+k-f+t-3}
(1-z^{2t})\\
&=
\sum_{k\ge f-t+1}(-1)^{k-f+t+1}\binom {k-1}{f-t-1}\\
&\kern3cm
\cdot
\left(\binom {2n+k-f+t-3}{n-k-1}-\binom {2n+k-f+t-3}{n-k-2t-1}\right),
\endalign$$
confirming our claim.\quad \quad \qed
\enddemo

\proclaim{Lemma \TTE}
With $x=y=\al=1$, we have
$$
\align
\coef{z^n}&
\left(\frac {z} {1-z}\right)^{f-t}
\frac {1-\left(\frac {z^2M(z)} {1-z}\right)^{t}}
{1+zM(z)}
\frac {\left(zM(z)\right)^2} {1- z^2M^2(z)}\\
&=
\sum_{k\ge f-t}(-1)^{k-f+t}\binom {k-1}{f-t-1}
\binom {2n+k-f+t-3}{n-k-2}\\
&\kern1cm
-
\sum_{k\ge f}(-1)^{k-f}\binom {k-1}{f-1}
\binom {2n+k-f-3}{n-k-t-2}.
\endalign$$
\endproclaim

\demo{Proof}
By the Lagrange inversion formula in Theorem~\TK\ with $F(z)=zM(z)$,
and the earlier considerations concerning the implied $f(z)$
under the specialisation $x=y=\al=1$
(see the paragraph above Lemma~\TO), we have
$$\allowdisplaybreaks
\align
\coef{z^n}&
\left(\frac {z} {1-z}\right)^{f-t}
\frac {1-\left(\frac {z^2M(z)} {1-z}\right)^{t}}
{1+zM(z)}
\frac {\left(zM(z)\right)^2} {1- z^2M^2(z)}\\
&=\coef{z^{-1}}
\frac {\left(\frac {z} {1+z+z^2}\right)^{f-t}
-z^{t}\left(\frac {z} {1+z+z^2}\right)^{f}}
{1+z}
\frac {z^2} {1- z^2}
\frac {(1+z)^{2(n+1)}} {z^{n+1}}\frac {1-z} {(1+z)^3}\\
&=\coef{z^{n-2}}
\left(\left(\frac {z} {1+z+z^2}\right)^{f-t}
-z^{t}\left(\frac {z} {1+z+z^2}\right)^{f}\right)
(1+z)^{2n-3}\\
&=\coef{z^{n-2}}
\sum_{k\ge f-t}(-1)^{k-f+t}\binom {k-1}{f-t-1}z^k
(1+z)^{2n+k-f+t-3}\\
&\kern1cm
-
\coef{z^{n-2}}
\sum_{k\ge f}(-1)^{k-f}\binom {k-1}{f-1}
z^{t+k}
(1+z)^{2n+k-f-3}\\
&=
\sum_{k\ge f-t}(-1)^{k-f+t}\binom {k-1}{f-t-1}
\binom {2n+k-f+t-3}{n-k-2}\\
&\kern1cm
-
\sum_{k\ge f}(-1)^{k-f}\binom {k-1}{f-1}
\binom {2n+k-f-3}{n-k-t-2},
\endalign$$
which is exactly the expression in the statement of the
lemma.\quad \quad \qed
\enddemo

\proclaim{Lemma \TTF}
With $x=y=\al=1$, we have
$$
\align
-\coef{z^n}&
\left(\frac { z} {1-z}\right)^{f-t}
\frac {\left( zM(z)\right)^{t}
-\left(\frac { z} {1-z}\right)^{t}} {1+ zM(z)}
\frac {
\left(zM(z)\right)^{t+1}} {1-z^2M^2(z)}\\
&=
\sum_{k\ge f}(-1)^{k-f}\binom {k-1}{f-1}
\binom {2n+k-f-3}{n-k-t-1}\\
&\kern1cm
-
\sum_{k\ge f-t}(-1)^{k-f+t}\binom {k-1}{f-t-1}
\binom {2n+k-f+t-3}{n-k-2t-1}.
\endalign$$
\endproclaim

\demo{Proof}
By the Lagrange inversion formula in Theorem~\TK\ with $F(z)=zM(z)$,
and the earlier considerations concerning the implied $f(z)$
under the specialisation $x=y=\al=1$
(see the paragraph above Lemma~\TO), we have
$$\allowdisplaybreaks
\align
-\coef{z^n}&
\left(\frac { z} {1-z}\right)^{f-t}
\frac {\left( zM(z)\right)^{t}
-\left(\frac { z} {1-z}\right)^{t}} {1+ zM(z)}
\frac {
\left(zM(z)\right)^{t+1}} {1-z^2M^2(z)}\\
&=
\coef{z^{-1}}
\left(\frac { z} {1+z+z^2}\right)^{f-t}
\frac {\left(\frac { z} {1+z+z^2}\right)^{t}
-z^{t}} {1+ z}
\frac {z^{t+1}} {1-z^2}
\frac {(1+z)^{2(n+1)}} {z^{n+1}}\frac {1-z} {(1+z)^3}\\
&=
\coef{z^{n-t-1}}
\left(\left(\frac { z} {1+z+z^2}\right)^{f}
-z^{t}\left(\frac { z} {1+z+z^2}\right)^{f-t}
\right)
(1+z)^{2n-3}\\
&=
\coef{z^{n-t-1}}
\sum_{k\ge f}(-1)^{k-f}\binom {k-1}{f-1}z^k
(1+z)^{2n+k-f-3}\\
&\kern1cm
-
\coef{z^{n-t-1}}
\sum_{k\ge f-t}(-1)^{k-f+t}\binom {k-1}{f-t-1}
z^{k+t}
(1+z)^{2n+k-f+t-3}\\
&=
\sum_{k\ge f}(-1)^{k-f}\binom {k-1}{f-1}
\binom {2n+k-f-3}{n-k-t-1}\\
&\kern1cm
-
\sum_{k\ge f-t}(-1)^{k-f+t}\binom {k-1}{f-t-1}
\binom {2n+k-f+t-3}{n-k-2t-1},
\endalign$$
as we claimed.\quad \quad \qed
\enddemo

\proclaim{Lemma \TTG}
With $x=y=\al=1$, we have
$$
\coef{z^n}
\frac { \left(zM(z)\right)^{f+t+2}} 
{(1+zM(z))^2}
=\binom {2n-3}{n-f-t-2}-\binom {2n-3}{n-f-t-3}.
$$
\endproclaim

\demo{Proof}
By the Lagrange inversion formula in Theorem~\TK\ with $F(z)=zM(z)$,
and the earlier considerations concerning the implied $f(z)$
under the specialisation $x=y=\al=1$
(see the paragraph above Lemma~\TO), we have
$$\align
\coef{z^n}
\frac { \left(zM(z)\right)^{f+t+2}} 
{(1+zM(z))^2}
&=
\coef{z^{-1}}
\frac { z^{f+t+2}} 
{(1+z)^2}\frac {(1+z)^{2(n+1)}} {z^{n+1}}\frac {1-z} {(1+z)^3}\\
&=
\coef{z^{n-f-t-2}}
(1+z)^{2n-3}(1-z)\\
&=\binom {2n-3}{n-f-t-2}-\binom {2n-3}{n-f-t-3},
\endalign$$
as required.\quad \quad \qed
\enddemo

\proclaim{Lemma \TTH}
With $x=y=\al=1$, we have
$$
\coef{z^n}
\frac {1-\left(z^2M^2(z)\right)^{t}}
{1-z^2M^2(z)} 
\cdot
\frac {\left( zM(z)\right)^{f-t+2}}
{1+zM(z)}
=\binom {2n-3}{n-f+t-2}-\binom {2n-3}{n-f-t-2}.
$$
\endproclaim

\demo{Proof}
By the Lagrange inversion formula in Theorem~\TK\ with $F(z)=zM(z)$,
and the earlier considerations concerning the implied $f(z)$
under the specialisation $x=y=\al=1$
(see the paragraph above Lemma~\TO), we have
$$\align
\coef{z^n}
\frac {1-\left(z^2M^2(z)\right)^{t}}
{1-z^2M^2(z)} 
\cdot
\frac {\left( zM(z)\right)^{f-t+2}}
{1+zM(z)}
&=
\coef{z^{-1}}
\frac {1-z^{2t}}
{1-z^2} 
\cdot
\frac {z^{f-t+2}}
{1+z}\frac {(1+z)^{2(n+1)}} {z^{n+1}}\frac {1-z} {(1+z)^3}\\
&=
\coef{z^{n-f+t-2}}
(1-z^{2t})
(1+z)^{2n-3}\\
&=\binom {2n-3}{n-f+t-2}-\binom {2n-3}{n-f-t-2},
\endalign$$
establishing the assertion of the lemma.\quad \quad \qed
\enddemo

\medskip
On the other hand, if $t\ge f$ then we
have to consider~(\Agg) with $x=y=\al=1$,
that is, after little simplification,
$$\allowdisplaybreaks
\align
&\left(\frac {z} {1-z}\right)^{t-f}
+\frac {2 } {1+zM(z)}
\frac {\left(zM(z)\right)^{t-f+2}
-\left(zM(z)\right)^{f+t+2}} {1-z^2M^2(z)}\\
&\kern1cm
+\frac { \left(zM(z)\right)^{f+t+2}} 
{(1+zM(z))^2}
+
\frac {zM(z)
\left(\left(zM(z)\right)^{t-f}-\left(\frac {z} {1-z}\right)^{t-f}\right)}
{(1+zM(z))^2}.
\tag\AN
\endalign$$
Lemmas~\TTI\ and~\TTJ\ below afford the coefficient extraction from the
individual terms in~(\AN) (together with Lemma~\TTG).

\proclaim{Lemma \TTI}
With $x=y=\al=1$, we have
$$
\coef{z^n}
\frac {1 } {1+zM(z)}
\frac {\left(zM(z)\right)^{t-f+2}
-\left(zM(z)\right)^{f+t+2}} {1-z^2M^2(z)}
=\binom {2n-3}{n+f-t-2}-\binom {2n-3}{n-f-t-2}.
$$
\endproclaim

\demo{Proof}
By the Lagrange inversion formula in Theorem~\TK\ with $F(z)=zM(z)$,
and the earlier considerations concerning the implied $f(z)$
under the specialisation $x=y=\al=1$
(see the paragraph above Lemma~\TO), we have
$$\align
\coef{z^n}
\frac {1 } {1+zM(z)}&
\frac {\left(zM(z)\right)^{t-f+2}
-\left(zM(z)\right)^{f+t+2}} {1-z^2M^2(z)}\\
&=
\coef{z^{-1}}
\frac {1 } {1+z}
\frac {z^{t-f+2}
-z^{f+t+2}} {1-z^2}
\frac {(1+z)^{2(n+1)}} {z^{n+1}}\frac {1-z} {(1+z)^3}\\
&=
\coef{z^{n}}
\left(z^{t-f+2}
-z^{f+t+2}\right)
(1+z)^{2n-3} \\
&=\binom {2n-3}{n+f-t-2}-\binom {2n-3}{n-f-t-2},
\endalign$$
in accordance with the assertion of the lemma.\quad \quad \qed
\enddemo

\proclaim{Lemma \TTJ}
With $x=y=\al=1$, we have
$$\align
\coef{z^n}&
\frac {zM(z)
\left(\left(zM(z)\right)^{t-f}-\left(\frac {z} {1-z}\right)^{t-f}\right)}
{(1+zM(z))^2}
=
\sum_{k\ge t-f+1}(-1)^{k-t-f-1}\binom {k-1}{t-f-1}\\
&\kern3cm
\cdot
\left(\binom {2n+k+f-t-3}{n-k-1}-\binom {2n+k+f-t-3}{n-k-2}\right).
\endalign$$
\endproclaim

\demo{Proof}
By the Lagrange inversion formula in Theorem~\TK\ with $F(z)=zM(z)$,
and the earlier considerations concerning the implied $f(z)$
under the specialisation $x=y=\al=1$
(see the paragraph above Lemma~\TO), we have
$$\allowdisplaybreaks
\align
\coef{z^n}&
\frac {zM(z)
\left(\left(zM(z)\right)^{t-f}-\left(\frac {z} {1-z}\right)^{t-f}\right)}
{(1+zM(z))^2}\\ 
&=\coef{z^{-1}}
\frac {z
\left(z^{t-f}-\left(\frac {z} {1+z+z^2}\right)^{t-f}\right)}
{(1+z)^2}\frac {(1+z)^{2(n+1)}} {z^{n+1}}\frac {1-z} {(1+z)^3}\\ 
&=\coef{z^{n-1}}
\left(z^{t-f}-\left(\frac {z} {1+z+z^2}\right)^{t-f}\right)
(1+z)^{2n-3}(1-z)\\ 
&=\coef{z^{n-1}}
\sum_{k\ge t-f+1}(-1)^{k-t-f-1}\binom {k-1}{t-f-1}z^k
(1+z)^{2n-3+k-t+f}(1-z)\\ 
&=
\sum_{k\ge t-f+1}(-1)^{k-t-f-1}\binom {k-1}{t-f-1}\\
&\kern3cm
\cdot
\left(\binom {2n+k+f-t-3}{n-k-1}-\binom {2n+k+f-t-3}{n-k-2}\right),
\endalign$$
as we claimed.\quad \quad \qed
\enddemo

\subhead 10. Proofs\endsubhead
Here we prove the enumeration results that we stated in
Section~1.

\demo{Proof of Theorem \TA}
We combine Proposition \TI\ in Section~3
with Lemmas~\TL--\TN\ in Section~5.\quad \quad \qed
\enddemo

\demo{Proof of Corollary \TB}
We sum the expression (\AAA) over all possible $c$ and $d$ such that
$c+d=n-2e-t$. This means that we have to compute
$$\multline
\chi(e=0)\binom {n-1}{t-1}
+\sum_{d\ge0}
\frac {(n-1)!} {(d+e)\cdot (n-d-2e-t)!\,d!\,(e-1)!\,(e+t-1)!}\\
-\sum_{b=e+t+1}^{n-e+1}(-1)^{b+e+t-1}
\frac {(n-b)\cdot(n-1)!} {b!\,
(e-1)!\,(n-b-e+1)!}\sum_{d=0}^ {b-e-t-1}(-1)^d\binom {b-1}d.
\endmultline
$$
The second sum over $d$ can be evaluated by using the simple formula
$$
\sum _{L=0} ^{K}(-1)^ {K-L}\binom ML=\binom {M-1}K.
\tag\AH$$
After replacing $b$ by $d+e+t+1$ in the result, we arrive exactly
at~(\AAB).\quad \quad \qed
\enddemo

\demo{Proof of Corollary \TC}
We sum the expression (\AAA) over all possible $c,d,e$ such that
$c+e+t=m$ and $c+d+2e+t=n$. This leads to the expression
$$\align
\chi(m=n)&\binom {n-1}{t-1}
+\sum_{e\ge1}
\frac {(n-1)!} {(n-m)\cdot(m-e-t)!\,(n-m-e)!\,(e-1)!\,(e+t-1)!}\\
&\kern1cm
-\sum_{b=n-m+t+1}^{n}(-1)^{n-b-m+t-1}\\
&\kern2cm
\cdot
\sum_{e\ge1}\frac {(n-b)\cdot(n-1)!} {b\cdot(n-m-e)!\,(b-1-n+m+e)!\,
(e-1)!\,(n-b-e+1)!}\\
&=
\chi(m=n)\binom {n-1}{t-1}
+\binom {n-1}{m-1}\sum_{e\ge1}
\binom {n-m-1}{e-1}\binom {m-1}{m-e-t}\\
&\kern1cm
-\sum_{b=n-m+t+1}^{n}(-1)^{n-b-m+t-1}\binom {n-1}{b}
\sum_{e\ge1}\binom {n-b}{e-1} \binom {b-1}{n-m-e}.
\endalign
$$
Both sums over $e$ can be evaluated by the Chu--Vandermonde summation formula
in binomial form (see e.g\. \cite{\GrKPAA, Sec.~5.1, Eq.~(5.27)})
$$
\sum _{L=0} ^{K}\binom ML\binom N{K-L}=\binom {M+N}K.
$$
Thus the above expression simplifies to
$$\align
\chi(m=n)&\binom {n-1}{t-1}
+
\binom {n-1}{m-1}
\binom {n-2}{m-t-1}\\
&\kern2cm
-\sum_{b=n-m+t+1}^{n}(-1)^{n-b-m+t-1}
\binom {n-1}{b}
\binom {n-1}{n-m-1}.
\endalign
$$
By (\AH), the sum over $b$ can be evaluated.
This turns the above expression into
$$\align
\chi(m=n)&\binom {n-1}{t-1}
+
\binom {n-1}{m-1}
\binom {n-2}{m-t-1}
-\binom {n-2}{m-t-2}
\binom {n-1}{m}\\
&=
\chi(m=n)\binom {n-1}{t-1}
+\frac {(nt+n-m)\,(n-1)!\,(n-2)!} {m!\,(n-m)!\,(m-t-1)!\,(n-m+t)!}\\
&=
\chi(m=n)\binom {n-1}{t-1}
+\chi(m<n)\frac {t} {n-1}\binom nm\binom {n-1}{m-t-1}\\
&\kern4cm
+\frac {1} {n-1}\binom {n-1}m\binom {n-1}{m-t-1},
\endalign$$
which is what we wanted to prove.\quad \quad \qed
\enddemo

\demo{Proof of Corollary \TD}
By combining (\AG)
with Lemma~\TL\ in Section~5 and Lemmas~\TO\ and~\TP\ in Section~6, we get
$$\multline
\binom {n-1}{t-1}
+\binom {2n-3}{n-t-2}-\binom {2n-3}{n-t-3}
+\binom {2n-3}{n-t-1}-\binom {2n-3}{n-t-2}-\binom {n-2}{t-1}\\
=\binom {2n-3}{n-t-1}-\binom {2n-3}{n-t-3}+\binom {n-2}{t-2}
=\binom {2n-2}{n-t-1}-\binom {2n-2}{n-t-2}+\binom {n-2}{t-2}
\endmultline
$$
for the desired number, as required.\quad \quad \qed
\enddemo

\demo{Proof of Theorem \TE}
As we explained in the paragraph at the beginning of Section~7,
what we have to do
is to put $x=y=1$ in the generating function
in Proposition~\TI,
differentiate the generating function
with respect to~$\al$, subsequently put $\al=1$, extract the
coefficient of~$z^n$, and finally divide
the result by the total number of relevant set-valued standard
tableaux given in Corollary~\TD.

For the first part, we have to add the corresponding
coefficient extractions from the second and the third term in~(\AAa)
(while the first term may be safely ignored since it does not
contain~$\al$). Hence, we have to sum the expressions given in
Lemmas~\TQ\ and~\TR\ in Section~7. This sum simplifies to
$$
\binom {2n-4}{n-t-1}
  + (n - 2) \binom {2 n - 4} {  n - t - 3}
  - (n + 1) \binom {2 n - 4} {n - t - 4} - 
 \binom {n - 3} {t-2}.
$$
Finally, this expression must be divided by the expression given
in Corollary~\TD. This quotient indeed agrees with the claimed
formula in the statement of Theorem~\TE.\quad \quad \qed
\enddemo

\demo{Proof of Theorem \TF}
For $t<f$,
we combine Proposition \TJ\ in Section~3 with $t<f$
with Lemmas~\TS--\TY\ in Section~8.
Subsequently, the resulting expression is considerably
simplified.

On the other hand, if $t\ge f$, then
we combine Proposition \TJ\ with $t\ge f$
with Lemma~\TL\ in Section~5 and Lemmas~\TZ, \TX, \TTA, and~\TTB\
in Section~8.
Again, the resulting expression is considerably
simplified.\quad \quad \qed
\enddemo

\demo{Proof of Theorem \TG}
For $t<f$, we use 
Lemma~\TL\ in Section~5 with $t$ replaced by~$f-t$ and
Lemmas~\TTC--\TTH\ in Section~9 to extract the coefficient
of~$z^n$ in~(\AL).
Subsequently, the resulting expression is considerably
simplified.

On the other hand, if $t\ge f$, then we use
Lemma~\TL\ with $t$ replaced by $t-f$ and
Lemmas~\TTI, \TTG, and~\TTJ\ in Section~9 to extract the coefficient
of~$z^n$ in~(\AN).
Also here, the resulting expression can be considerably
simplified.\quad \quad \qed
\enddemo

\centerline{\smc AI Declaration}

\medskip
AI has been used for proofreading, but for nothing else.

\enddocument